\documentclass[reqno]{amsart}
\usepackage{eurosym}
\usepackage{amsmath}
\usepackage{amsfonts}
\usepackage{amssymb}
\usepackage{amsfonts}
\usepackage{epsfig}
\usepackage{caption}
\usepackage{subcaption}
\usepackage{multirow}
\usepackage{enumitem}
\usepackage{stackengine}
\stackMath

\newtheorem{theorem}{Theorem}

\newtheorem{remark}[theorem]{Remark}

\newtheorem{Thm}{Theorem}[section]
\newtheorem{lem}[Thm]{Lemma}
\newtheorem{pro}[Thm]{Proposition}
\newtheorem{de}[Thm]{Definition}

\begin{document}
	
\title[Minimization of free energy in a PB problem]
{Minimization of free energy in a Poisson--Boltzmann problem}

\author[S. Ali]{S. Ali}
\address{Dipartimento di Matematica, Informatica e Geoscienze,
	Universit\`a degli Studi di Trieste,
	Alfonso Valerio 12/1, 34127 Trieste, Italy}
\email{shahid.ali@phd.units.it}

\author[K. Bohinc]{K. Bohinc}
\address{Faculty of Health Sciences, University of Ljubljana,
	Zdravstvena 5, SI-1000 Ljubljana, Slovenia}
\email{klemen.bohinc@zf.uni-lj.si}

\author[A. Mart\'{\i}nez]{A. Mart\'{\i}nez}
\address{Dipartimento di Matematica, Informatica e Geoscienze,
	Universit\`a degli Studi di Trieste,
	Alfonso Valerio 12/1, 34127 Trieste, Italy}
\email{amartinez@units.it}

\author[S. Maset]{S. Maset}
\address{Dipartimento di Matematica, Informatica e Geoscienze,
	Universit\`a degli Studi di Trieste,
	Alfonso Valerio 12/1, 34127 Trieste, Italy}
\email{maset@units.it}

\date{}



\begin{abstract}
	In this paper, we study a mathematical model in which an ionic
	solution composed of several types of spherical ions is confined
	between two uniformly charged parallel planes. The ion concentrations
	at equilibrium, determined by the interplay between thermal and
	electrostatic effects, are obtained by minimizing the free energy of
	the system. The resulting minimization problem is infinite-dimensional,
	and we prove the existence and uniqueness of its solution by means of
	a norm-convergence direct method based on the Banach--Saks property of
	the space in which the problem is posed. The existence proof is
	technically challenging because the functional to be minimized is not
	coercive in general. Three illustrative problems are considered and
	solved numerically.
\end{abstract}

\maketitle

\section{Introduction}

The study of electrostatic interactions between charged structures and ions is crucial for understanding the behavior of soft matter, particularly living matter (see \cite{6a}). Biological systems typically contain a multitude of charged structures. Membranes, cellular components, globular proteins, polyelectrolytes, and DNA are examples of biological structures carrying electric charges and, consequently, exhibiting electrostatic behavior (see \cite{2a,4a}). These charged structures are immersed in biological media, such as blood, containing free ions in addition to water molecules. Understanding the physical processes occurring between two negatively or positively charged structures is important from the biological point of view (see \cite{6a}).

We model two negatively or positively charged structures as two uniformly charged infinite parallel planes, with an ionic solution confined between them. The ions in the solution may be either \emph{counter-ions}, i.e., ions whose charge has the opposite sign to that of the planes, or \emph{co-ions}, i.e., ions whose charge has the same sign. The ions undergo thermal motion and, at the same time, are subject to the electrostatic forces generated by the charged planes and by the other ions. Of particular interest is the determination of the ion concentrations at thermal and electrostatic equilibrium. This is known as a \emph{Poisson--Boltzmann (PB) problem}. At equilibrium, the ion concentrations are those that minimize the \emph{free energy} of the system.

The aim of this paper is to present such a model for the case of spherical ions and to prove the existence and uniqueness of the solution to the PB problem, i.e., the existence and uniqueness of the ion concentrations minimizing the free energy. The existence proof is based on a direct method, but it does not employ standard weak-compactness arguments, such as those based on dual spaces and the Banach--Alaoglu theorem. Instead, it relies on the Banach--Saks property of the infinite-dimensional space in which the PB problem is formulated and uses norm convergence. Moreover, because the functional to be minimized may fail to be coercive, the proof requires nonstandard and technically challenging arguments, which also yield an essential boundedness result. Existence and uniqueness results for free-energy minimization in PB problems different from the one considered in the present paper, and based on different approaches, can be found in \cite{Li} and \cite{Huang}.

By means of a variational procedure based on the Euler--Lagrange equations, the free-energy minimization problem can be reformulated as a boundary value problem with Neumann conditions for a second-order differential equation whose solution is the electrostatic potential (see \cite{Armenta}). This potential determines the ion concentrations minimizing the free energy. When the ions are modeled as point particles, the resulting differential equation is an ordinary differential equation. More realistic models can be obtained by considering ions with a geometrical structure, such as rod-like or spherical ions. In these cases, the resulting differential equation is an integro-differential equation. Such integro-differential equations, and their numerical treatments, are presented in \cite{Maset,May,Grime} for rod-like ions and in \cite{3a,S,S1,Bohinc2023} for spherical ions. These papers do not address the question of existence and uniqueness of the solution of the boundary value problem for the second-order integro-differential equation. 

The present paper includes three illustrative examples of PB problems for the spherical ions model considered here, in which the ion concentrations minimizing the free energy are computed numerically. Unlike the approaches adopted in the papers cited above, these concentrations are not obtained by reformulating the minimization problem as a boundary value problem for the potential. Instead, they are computed directly by minimizing the free energy. This approach avoids the lengthy mathematical calculations required by the variational procedure used to derive the boundary value problem. Since the minimization problem is infinite-dimensional, it must first be discretized into a finite-dimensional minimization problem, after which a standard optimization solver can be used to compute its minimum. The description and analysis of the discretization of the infinite-dimensional minimization problem are not included here and will be the subject of a companion paper.

The present paper is organized as follows. Section~2 describes the spherical ions model. Section~3 provides the expression of the free energy of the system. Section~4 formulates the infinite-dimensional minimization problem and its dimensionless version. Section~5 proves that the infinite-dimensional minimization problem has a unique solution. Section~6 presents the three illustrative examples of PB problems together with their numerical solutions. Conclusions are in Section~7. Finally, the appendix contains the results on which the Banach--Saks direct method of proof is based.

\section{The spherical ions model}

The spherical ions model presented here is more general than those considered in
\cite{S,Bohinc2023}, because it allows for an arbitrary number of different
types of spherical ions.

\begin{figure}[h]
\centering
\includegraphics[height=6cm,width=7cm]{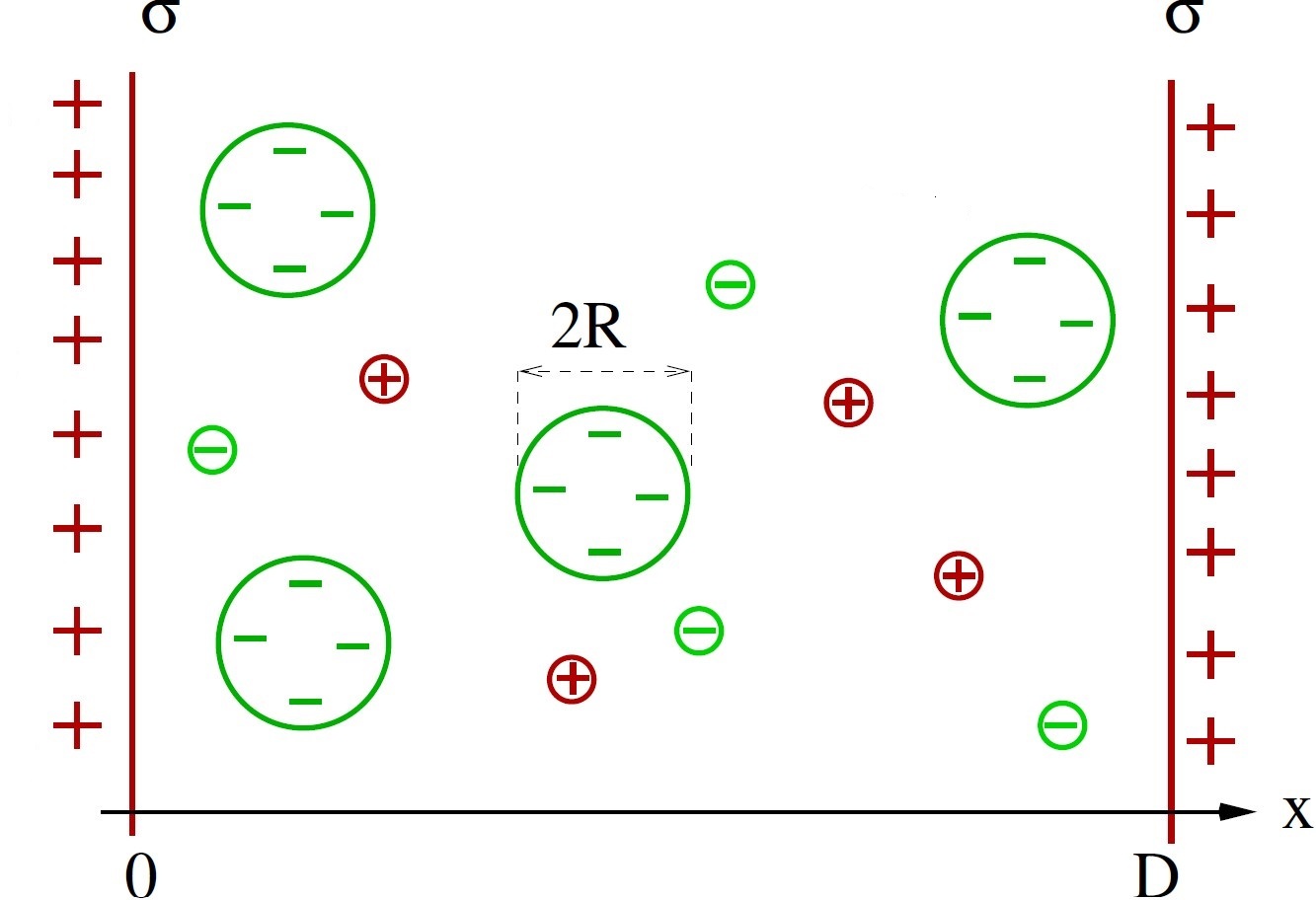}
\caption{Schematic illustration of two like-charged infinite parallel planes (located
at positions $x=0$ and $x=D$ along the $x$ axis) with surface charge density $%
\protect\sigma$. The ionic solution between the planes consists of three different types of ions: two types of counter-ions and one type of co-ions.}
\label{Fig1}
\end{figure}

As remarked above in the introduction, the two negatively or positively charged structures
are modeled by two infinite parallel planes. We consider the two planes at a distance $%
D $ and uniformly charged with surface charge density $\sigma$. Between the
planes, there is the ionic solution consisting of $N>1$ different types of
spherical ions. An ion of type $i$, $i\in\{1,\ldots,N\}$, is a sphere of radius $%
R_i$ and has (positive or negative) \emph{valence} $Z_i$. The charge on the $i$-th type of ion is therefore $Z_ie$, where $e$ is the (positive) elementary charge. The charge on the spherical ions is assumed to be uniformly distributed on their
surface. Figure \ref{Fig1} illustrates the situation.

We assume that the spherical ions cannot penetrate the planes. Therefore, we have the
constraints 
\begin{equation}
2R_{i}<D,\ i\in\{1,\ldots,N\},  \label{costraints}
\end{equation}
for the radii of the ions. The limiting case $2R_i=D$ is not addressed, since in this case the
ions cannot move in the direction normal to the planes.  

The abscissa $x$ denotes the distance from the left plane of a point between the two planes. For $i\in\{1,\ldots ,N\}$, $n_{i}(x)$ denotes the \emph{%
concentration} of the $i$-th type of ions at a point of abscissa $x$, i.e., $%
n_{i}(x)$ is the number of spherical ions of type $i$ per unit volume whose centers have abscissa $x$. Since the ions cannot penetrate the
planes, we have 
\begin{equation}
n_{i}(x)=0\text{ for $x\notin \left[ R_{i},D-R_{i}\right]$.}  \label{zero}
\end{equation}

Note that we are dealing with a one-dimensional model, as the ion concentrations at equilibrium 
remain constant in any plane parallel to the two charged planes. This appears
evident from translational invariance: all points in a plane at a distance $ x $ from the
left plane perceive the same configuration, namely, two uniformly charged
infinite parallel planes, with one plane on the left at a distance $x$ and the other on the right at a distance
$D-x$.

\section{The free energy}

The \emph{free energy} of a thermodynamic system is a potential that attains its minimum value at equilibrium under the appropriate constraints. In this section, we present the expression for the free energy of the system consisting of the charged planes with the ionic solution between them.

In the region consisting of a right prism with bases of area $A$ on the
two planes, the free energy $F$ is given by 
\begin{equation}
\frac{F}{Ak_{B}T}=\frac{F^{\mathrm{en}}}{Ak_{B}T}+\frac{F^{\mathrm{el}}}{Ak_{B}T}=\sum\limits_{i=1}^{N}\int\limits_{R_i}^{D-R_i}S_{i}\left( x\right) dx+\int\limits_{0}^{D}\frac{E\left( x\right) ^{2}}{8\pi l_{B}}
dx.
\label{F}
\end{equation}%
Here $F$ is divided by the area $A$ and $k_BT$, where $k_B$ is the Boltzmann
constant and $T$ is the absolute temperature of the ionic solution. Observe that $%
\frac{F}{Ak_{B}T}$ has dimensions $L^{-2}$ ($L$ denotes length). In (\ref{F}), the free
energy (divided by $Ak_BT$) is split into the sum of two terms: the \emph{entropic term} 
\begin{equation}
	\frac{F^{\mathrm{en}}}{Ak_{B}T}=\sum\limits_{i=1}^{N}\int\limits_{R_i}^{D-R_i}S_{i}\left( x\right) dx  \label{ent}
\end{equation}
associated with the thermal motion, and the \emph{%
electrostatic term} 
\begin{equation}
\frac{F^{\mathrm{el}}}{Ak_{B}T}=\int\limits_{0}^{D}\frac{E\left( x\right) ^{2}}{8\pi l_{B}}dx  \label{elt}
\end{equation}%
associated with the electrostatic interactions. These two terms are described in the next two subsections. 

\subsection{The entropic term}\label{et}

In the entropic term (\ref{ent}), where the integrals are from $R_i$ to $D-R_i$ since the ions cannot penetrate the planes, we have
\begin{equation}
S_{i}\left( x\right)=
n_i\left( x\right) \left( \log \frac{n_{i}\left(
x\right) }{n_{i0}}-1\right) +n_{i0},\ x\in[R_i,D-R_i]\text{\ and\ }i\in\{1,\ldots,N\}.
\label{Ai}
\end{equation}
The quantity $S_i(x)$ is the volume density at the point of abscissa $x$ of the entropic relative free energy $\frac{F^{\mathrm{en}}}{k_{B}T}$, for the $i$-th type of ion. The positive constant concentrations 
$n_{i0}$, $i\in\{1,\ldots,N\}$, are known and called \emph{bulk
concentrations}. The bulk concentrations are the concentrations of the different types of ions in the solution, outside the region between the charged structures modeled by the two parallel infinite planes. Therefore, they are the concentrations at equilibrium when the charged planes are removed. Further details on bulk concentrations are given in Subsection \ref{bulk}.

Observe that the region between the planes is treated as being in contact
with a reservoir (the bulk), in which the different ion species have their
prescribed bulk concentrations. Therefore, the region is not a closed
system, and the total number of ions within it is not fixed.

In (\ref{Ai}) and in the following, we assume $0\log 0=0$.

The expression \((\ref{ent})\), with \((\ref{Ai})\), for the entropic term is obtained by treating the thermal motion of the spherical ions as that of an ideal gas, under the
assumption that the ions are allowed to overlap.
For a derivation of (\ref{ent}), with (\ref{Ai}), see \cite{B}.

\subsection{The electrostatic term}

In the electrostatic term (\ref{elt}), $l_{B}$ is the \emph{Bjerrum length},
that is, the distance at which the electrostatic interaction energy between two elementary charges is equal to the thermal energy $k_BT$, and 
\begin{equation}
E\left( x\right) =4\pi l_{B}\left( \frac{\sigma }{e}+\int\limits_{0}^{x} 
\frac{\rho \left( z\right) }{e}dz\right) ,\ x\in \left[ 0,D\right],
\label{eq1}
\end{equation}%
is the electric field at the point of abscissa $x$. 
This electric field is the negative gradient of a dimensionless potential, and so its
dimension is $L^{-1}$, in accordance with the dimension $L^{-2}$ of both
sides in (\ref{F}).

In (\ref{eq1}), $\rho (z)$ is the \emph{volume
charge density} at the point of abscissa $z\in[0,x]$ and it is given by 
\begin{equation}
 \frac{\rho \left( z\right) }{e}=\sum\limits_{i=1}^{N}Z_i\cdot\frac{1}{2R_{i}}
 \int\limits_{z-R_{i}}^{z+R_{i}}n_{i}\left(s\right) ds,\ z\in \left[ 0,D%
 \right].  \label{eq1.5}
\end{equation}%
Observe that (\ref{zero}) holds for the integrand functions $n_i$, $i\in\{1,\ldots,N\}$, on the right-hand
side of (\ref{eq1.5}).

As with the entropic term, the expression (\ref{elt}), with (\ref{eq1}) and (\ref{eq1.5}), is derived under the assumption that the spherical ions may overlap. For a derivation of (\ref{elt}), with (\ref{eq1}), see \cite{Jackson}, and for a derivation of (\ref{eq1.5}), see \cite{3a}.

Observe that (\ref{eq1.5}) is not defined if $R_i=0$ for some $i\in\{1,\ldots,N\}$. This case will be considered below in Subsection \ref{Facts}.

\section{The minimization problem}

We are interested in determining the concentrations of the different types of
ions at thermal and electrostatic equilibrium. These concentrations are those that minimize the free energy $F$, once we
consider $F$ as a function of the concentrations.

Hence, we have to minimize $F$ with respect to the non-negative functions $%
n_{i}:\mathbb{R}\rightarrow \mathbb{R}$, $i\in\{1,\ldots ,N\}$, which are
considered zero outside $\left[ R_{i}, D-R_i\right]$ by recalling (\ref{zero}). Of course, the
non-negativity of the functions $n_i$ is imposed since these
functions represent concentrations. This minimization is accomplished under
the constraint of the \emph{electroneutrality condition} 
\begin{equation}
2\sigma+\int\limits_{0}^{D}\rho(z)dz=0,  \label{electron}
\end{equation}
where $\rho$ is the volume charge density given in (\ref{eq1.5}). This
condition states that the total charge is zero in any right prism with bases on
the planes. 

Along with determining the concentrations $n_{i}^{\min}$, $i\in\{1,\ldots,N\}$,
minimizing $F$, we are interested in the minimal free energy $F^{\min }$. It
is important to understand how $F^{\mathrm{min}}$ depends on the distance $D$
between the two planes, i.e., to determine the function $F^{\min }(D)$. In particular, we are interested in determining whether this function of $D$
has a minimum. If this happens, the
charged structures, modeled as the two parallel planes, tend to move towards
the distance $D^{\min}$ at which we have the minimum of the function $%
F^{\min}(D)$. Thus, if the charged structures are at a distance \(D\) in a right neighborhood of \(D^{\min}\) in which \(F^{\min}\) is increasing, they attract each other, although they are charged with the same sign. This phenomenon, which has also been observed experimentally in other
physical settings, is called \emph{like-charge attraction} (see \cite{Wang}). In the
present model, it is mediated by the counter-ions in the solution between
the charged structures; see, for example, \cite{3a,S,S1}. In this context, it is important to look at the 
\emph{pressure} 
\begin{equation}
	P(D)=-\frac{\left(F^{\min}\right)^\prime(D)}{A}, \label{pressure}
\end{equation}
which is the pressure exerted by the two charged planes towards the exterior,
when they are at a distance $D$. Here, $A$ is the same area as in (\ref{F}). Like-charge attraction is characterized by a
negative pressure.

\subsection{Bulk concentrations}\label{bulk}
The constant bulk concentrations $n_{i0}$, $i\in\{1,\ldots,N\}$, have been defined in Subsection \ref{et}. They are the ion concentrations at equilibrium when the charged planes are removed. This means that they are the concentrations minimizing the free energy (\ref{F}), where now in the electrostatic term (\ref{elt}) we have
\begin{equation*}
	E\left( x\right) =4\pi l_{B}\int\limits_{0}^{x} 
	\frac{\rho \left( z\right) }{e}dz,\ x\in \left[ 0,D\right],
\end{equation*}%
with $\rho(z)$ given in (\ref{eq1.5}), but the concentrations $n_i$, $i\in\{1,\ldots,N\}$, no longer have to satisfy (\ref{zero}), since the charged planes are removed. For the same reason, now we have
$$
\frac{F^{\mathrm{en}}}{Ak_{B}T}=\sum\limits_{i=1}^{N}\int\limits_{0}^{D}S_{i}\left( x\right) dx
$$
for the entropic term.

In this minimization problem without the charged planes, the electroneutrality condition (\ref{electron}) becomes
\begin{equation*}
	\int\limits_{0}^{D}\rho(z)dz=0.
\end{equation*}
For the bulk concentrations, this electroneutrality condition reads 
\begin{equation}
\sum\limits_{i=1}^N Z_in_{i0}=0.  \label{electronb}
\end{equation}
 In the following, we assume that the known valences $Z_i$  and bulk concentrations $n_{i0}$, $i\in\{1,\ldots,N\}$, satisfy (\ref{electronb}). Under this assumption, it follows immediately that the bulk concentrations have zero free energy and therefore minimize it, in the case in which the charged planes are removed.

\subsection{Nondimensionalization} \label{nondim}

Now, we introduce an
equivalent version of the minimization problem, in which all the quantities involved are dimensionless. This reduces the number of parameters involved and yields a cleaner problem formulation. 

We write the bulk concentrations as 
\begin{equation*}
n_{i0}=c_{i}n_{0},\ i\in\{1,\ldots ,N\}.
\end{equation*}%
Here, $n_{0}$ is some reference value for bulk concentrations, such as their average, maximum value or minimum value. The numbers $c_i$ can be considered
nondimensional bulk concentrations. Moreover, we introduce the
nondimensional abscissa 
\begin{equation*}
\widehat{x}=\frac{x}{D}\in[0,1],\ x\in[0,D],
\end{equation*}%
the nondimensional concentrations 
\begin{equation*}
\widehat{n}_{i}\left( \widehat{x}\right) =\frac{n_{i}\left( D\widehat{x}
\right) }{n_{i0}},\ \widehat{x}\in[0,1]\text{\ and\ }i\in\{1,\ldots ,N\},
\end{equation*}%
and the nondimensional radii%
\begin{equation*}
\widehat{R}_{i}=\frac{R_{i}}{D},\ i\in\{1,\ldots ,N\}.
\end{equation*}
Observe that the nondimensional concentration $\widehat{n}_i$, $i\in\{1,\ldots,N\}$%
, is considered to be zero outside $[\widehat{R}_i,1-\widehat{R}_i]$. Moreover, we have $%
\widehat{R}_i<\frac{1}{2}$ by (\ref{costraints}).

The next proposition expresses the nondimensional free energy 
\begin{equation*}
\widehat{F}:=\frac{F}{Dn_{0}Ak_{B}T}
\end{equation*}
in terms of the nondimensional quantities introduced above. Moreover, the proposition
also provides a nondimensional version of the electroneutrality condition (%
\ref{electron}).

\begin{pro}
Let
\begin{equation}
	\widehat{S}_{i}\left( \widehat{x}\right):=c_i\left( \widehat{n}_{i}\left( \widehat{x}\right)
	\left( \log \widehat{n}_{i}\left(\widehat{x}\right) -1\right) +1\right) ,\ \widehat{x}\in [\widehat{R}_i,1-\widehat{R}_i] 
	\text{\ and\ }i\in\{1,\ldots,N\},  \label{Aihat}
\end{equation}
be the nondimensional entropic volume densities of the relative free energy. Let 
\begin{equation}
\widehat{\sigma}:=\frac{\sigma }{eDn_{0}}  \label{sigmahat}
\end{equation}
be the nondimensional surface charge density of the plates; let 
\begin{equation}
\widehat{\rho}(\widehat{z}):=\sum\limits_{i=1}^{N}Z_{i}c_{i}\cdot \frac{1}{2\widehat{R}_i}%
\int\limits_{\widehat{z}-\widehat{R}_{i}}^{\widehat{z}+\widehat{R}_{i}}\widehat{n}%
_{i}\left(\widehat{s}\right)d\widehat{s},\ \widehat{z}\in[0,1],  \label{rhohat}
\end{equation}
be the nondimensional volume charge density, and let 
\begin{equation}
\widehat{E}\left( \widehat{x}\right):=\widehat{\sigma} +\int\limits_{0}^{\widehat{x}}\widehat{%
\rho} \left( \widehat{z}\right)d\widehat{z},\ \widehat{x}\in \left[ 0,1\right],  \label{Ehat}
\end{equation}
be the nondimensional electric field.

The nondimensional free energy
is given by 
\begin{equation}
\widehat{F}:=\widehat{F}^{\mathrm{en}}+\widehat{F}^{\mathrm{el}}=\sum\limits_{i=1}^{N}\int\limits_{\widehat{R}_i}^{1-\widehat{R}_i}\widehat{S}_{i}\left(
\widehat{x}\right) d\widehat{x}+K\int\limits_{0}^{1}\widehat{E}\left( \widehat{x}\right)
^{2}d\widehat{x},  \label{Fhat}
\end{equation}%
where
\begin{equation*}
K:=2\pi l_{B}D^{2}n_{0}. \label{Kappa}
\end{equation*} 

The nondimensional electroneutrality condition is 
\begin{equation}
2\widehat{\sigma}+\int\limits_{0}^{1}\widehat{\rho}(\widehat{z})d\widehat{z}=0
\label{electronhat}
\end{equation}
and the nondimensional electroneutrality condition for bulk concentrations is
\begin{equation}
	\sum\limits_{i=1}^NZ_ic_i=0. \label{bulken}
\end{equation}

\end{pro}

\begin{proof}
We write the entropic term (\ref{ent}) by using nondimensional quantities. For $i\in\{1,\ldots ,N\}$, we have
\begin{eqnarray*}
	\int\limits_{R_i}^{D-R_i}S_{i}\left( x\right) dx
	&=&\int\limits_{R_i}^{D-R_i}\left( n_{i}\left( x\right) \left( \log 
	\frac{n_{i}\left( x\right) }{n_{i0}}-1\right) +n_{i0}\right) dx \\
	&=&D\int\limits_{\widehat{R}_i}^{1-\widehat{R}_i}\left( n_{i}\left( D%
	\widehat{x}\right) \left( \log \frac{n_{i}\left( D\widehat{x}\right) }{n_{i0}%
	}-1\right) +n_{i0}\right) d\widehat{x} \\
	&=&Dn_{0}\int\limits_{\widehat{R}_i}^{1-\widehat{R}_i}c_{i}\left( 
	\widehat{n}_{i}\left( \widehat{x}\right) \left( \log \widehat{n}_{i}\left( 
	\widehat{x}\right) -1\right) +1\right) d\widehat{x} \\
	&=&Dn_{0}\int\limits_{\widehat{R}_i}^{1-\widehat{R}_i}\widehat{S}%
	_{i}\left( \widehat{x}\right) d\widehat{x},
\end{eqnarray*}%
where $\widehat{S}_{i}\left( \widehat{x}\right)$ is given in (\ref{Aihat}). Thus, the entropic term is
\begin{equation*}
	\frac{F^{\mathrm{en}}}{Ak_{B}T}=\sum\limits_{i=1}^{N}\int\limits_{R_i}^{D-R_i}S_{i}\left( x\right)
	dx=Dn_{0}\sum\limits_{i=1}^{N}\int\limits_{\widehat{R}_i}^{1-\widehat{R}_i}\widehat{S}_{i}\left( \widehat{x}\right) d\widehat{x}.
\end{equation*}

Now, we write the electrostatic term (\ref{elt}) by using nondimensional quantities. For $\widehat{z}\in \left[ 0,1\right] 
$, we have
\begin{eqnarray}
	\frac{\rho \left( D\widehat{z}\right) }{e} &=&\sum\limits_{i=1}^{N}Z_i\cdot\frac{1}{2R_{i}}%
	\int\limits_{D\widehat{z}-R_{i}}^{D\widehat{z}+R_{i}}n_{i}\left(s\right) ds=\sum\limits_{i=1}^{N}Z_i\cdot \frac{1}{2\widehat{R}_{i}}%
	\int\limits_{\widehat{z}-\widehat{R}_{i}}^{\widehat{z}+\widehat{R}_{i}}n_{i}\left(D\widehat{s}\right) d\widehat{s}\notag \\
	&=&n_0\sum\limits_{i=1}^{N}Z_ic_i\cdot \frac{1}{2\widehat{R}_{i}}%
	\int\limits_{\widehat{z}-\widehat{R}_{i}}^{\widehat{z}+\widehat{R}_{i}}\widehat{n}_{i}\left(\widehat{s}\right) d\widehat{s}=n_{0}\widehat{\rho }\left( \widehat{z}\right), \notag \\
	\label{Dn0}
\end{eqnarray}%
where $\widehat{\rho }\left( \widehat{z}\right)$ is given in (\ref{rhohat}). Hence, for $\widehat{x}\in \left[ 0,1\right] $, we have
\begin{eqnarray}
	E\left( D\widehat{x}\right) &=&4\pi l_{B}\left( \frac{\sigma }{e}%
	+\int\limits_{0}^{D\widehat{x}}\frac{\rho \left( z\right) }{e}dz\right)=4\pi l_{B}\left( \frac{\sigma }{e}%
	+D\int\limits_{0}^{\widehat{x}}\frac{\rho \left( D\widehat{z}\right) }{e}d\widehat{z}\right)\notag \\
	&=&4\pi l_{B}\left( \frac{\sigma }{e}+Dn_{0}\int\limits_{0}^{\widehat{x}}%
	\widehat{\rho} \left( \widehat{z}\right)d\widehat{z}\right) \text{\ \ (see (\ref{Dn0}))}\notag \\
	&=&4\pi l_{B}Dn_{0}\left( \frac{\sigma }{eDn_{0}}+\int\limits_{0}^{\widehat{%
			x}}\widehat{\rho }\left( \widehat{z}\right)d\widehat{z}\right)=4\pi l_{B}Dn_{0}\widehat{E}\left( \widehat{x}\right), \notag\\
		\label{Dn02}
\end{eqnarray}%
where $\widehat{E}(\widehat{x})$ is given in (\ref{Ehat}) (recall (\ref{sigmahat})). Thus, the electrostatic term becomes%
\begin{eqnarray*}
	\frac{F^{\mathrm{el}}}{Ak_{B}T}&=&\int\limits_{0}^{D}\frac{E\left( x\right) ^{2}}{8\pi l_{B}}%
	dx=D\int\limits_{0}^{1}\frac{E\left( D\widehat{x}\right)^2 }{8\pi l_{B}}d%
	\widehat{x}\\
	&=&D\int\limits_{0}^{1}\left( 4\pi l_{B}Dn_{0}\right) ^{2}\frac{\widehat{E}%
		\left( \widehat{x}\right) ^{2}}{8\pi l_{B}}d\widehat{x}  \text{\ \ (see (\ref{Dn02}))}\\
	&=&2\pi l_{B}D^{3}n_{0}^{2}\int\limits_{0}^{1}\widehat{E}\left( \widehat{x}%
	\right) ^{2}d\widehat{x}.
\end{eqnarray*}
The free energy $F$ (divided by $Ak_BT$) can now be rewritten as%
\begin{eqnarray*}
	\frac{F}{Ak_{B}T}=\frac{F^{\mathrm{en}}}{Ak_{B}T}+\frac{F^{\mathrm{el}}}{Ak_{B}T}=Dn_{0}\sum\limits_{i=1}^{N}\int%
	\limits_{\widehat{R}_i}^{1-\widehat{R}_i}\widehat{S}_{i}\left( \widehat{x%
	}\right) d\widehat{x}+2\pi l_{B}D^{3}n_{0}^{2}\int\limits_{0}^{1}\widehat{E}%
	\left( \widehat{x}\right) ^{2}d\widehat{x}.
\end{eqnarray*}%
By dividing by $Dn_0$, we obtain the nondimensional free energy (\ref{Fhat}).

Regarding the nondimensionalization of the electroneutrality condition (\ref{electron}), note that this condition can be written as
\begin{eqnarray*}
	0=2\sigma+\int\limits_{0}^{D} \rho(x)dx=2\sigma+eD\int\limits_{0}^{1}\frac{\rho \left( D\widehat{z}\right) }{e}d\widehat{z}=2\sigma+eDn_0\int\limits_{0}^{1}\widehat{\rho}(\widehat{z})d\widehat{z}
\end{eqnarray*}
by using (\ref{Dn0}). A division by $eDn_0$ provides the nondimensional electroneutrality condition (\ref{electronhat}) (recall (\ref{sigmahat})). Finally, (\ref{bulken}) follows immediately from (\ref{electronb}) by writing $n_{i0}=c_in_0$, $i\in\{1,\ldots,N\}$.
\end{proof}

\subsection{Formulation as the  minimization problem MP on the infinite-\\ dimensional space $X$}\label {MPsection}

To avoid cumbersome notation, from now on, the nondimensional abscissas 
$\widehat{x}$, $\widehat{z}$ and $\widehat{s}$ are written without hats. 

We consider the spaces $L^p(\widehat{R}_i,1-\widehat{R%
}_i)$, $i\in\{1,\ldots,N\}$, where $p\in(1,+\infty)$. The functions in $L^p(\widehat{R}_i,1-\widehat{R%
}_i)$ are extended by zero values
outside the interval $[\widehat{R}_i,1-\widehat{R}_i]$.

We minimize the nondimensional free energy $\widehat{F}$ with respect
to nondimensional concentrations $\widehat{n}_i\in L^p(\widehat{R}_i,1-\widehat{R%
}_i)$, $i\in\{1,\ldots,N\}$.

Regarding the entropic part $\widehat{F}^{\mathrm{en}}$, we have in (\ref{Fhat}), for $i\in\{1,\ldots ,N\}$, 
\begin{eqnarray*}
	\int\limits_{\widehat{R}_i}^{1-\widehat{R}_i}\widehat{S}_{i}\left( x\right) dx
	&=&c_i\int\limits_{\widehat{R}_i}^{1-\widehat{R}_i}\nu(\widehat{n}_{i}(x))dx=c_i\left\langle \nu\circ \widehat{n}_{i}\right\rangle_i,
\end{eqnarray*}%
where $\nu:[0,+\infty)\rightarrow \mathbb{R}$ is the function
\begin{equation}
	\nu(u)=u(\log u-1)+1,\ u\geq 0, \label{nu}
\end{equation}
and $\left\langle\ \cdot\ \right\rangle_i$ denotes the integral 
\begin{equation*}
	\left\langle f\right\rangle_i =\int\limits_{\widehat{R}_i}^{1-\widehat{R}_i}f\left( x\right)dx,\ f\in  L^1(\widehat{R}_i,1-\widehat{R%
	}_i).
\end{equation*}
The fact that the function $x\mapsto \nu(\widehat{n}_i(x))$, $x\in[\widehat{R}_i,1-\widehat{R}_i]$, belongs to $L^1(\widehat{R}_i,1-\widehat{R%
}_i)$ is proved in Proposition \ref{Lp-L1} below. Therefore, for the entropic part, we have
$$
\widehat{F}^{\mathrm{en}}=\sum\limits_{i=1}^{N}c_{i}\left\langle \nu\circ \widehat{n}%
_i\right\rangle_i.
$$

Regarding the electrostatic part $\widehat{F}^{\mathrm{el}}$, we have

\begin{equation}
\int\limits_{0}^{x}\widehat{\rho} \left( z\right)dz=\sum\limits_{i=1}^{N} 
Z_ic_i\int\limits_{0}^{x}\frac{1}{2\widehat{R}_i}\int\limits_{z-\widehat{R%
}_{i}}^{z+\widehat{R}_{i}} \widehat{n}_{i}\left(s\right)dsdz,\ x\in \lbrack
0,1],  \label{integral}
\end{equation}%
in (\ref{Ehat}). Thus, we can write the nondimensional electric field as 
\begin{equation*}
\widehat{E}=\widehat{\sigma} +\sum\limits_{i=1}^{N}Z_ic_{i}\mathcal{L}_{i}%
(\widehat{n}_{i}),
\end{equation*}%
where, for $i\in\{1,\ldots ,N\}$, $\mathcal{L}_{i}:L^p(\widehat{R}_i,1-\widehat{R%
}_i)\rightarrow C\left( \left[
0,1\right] \right) $ is the linear operator given by%
\begin{equation}
[\mathcal{L}_{i}(u)] \left( x\right) =%
\int\limits_{0}^{x}\frac{1}{2\widehat{R}_i}\int\limits_{z-\widehat{R}_{i}}^{z+\widehat{R}_{i}}
u\left(s\right)dsdz,\ x\in \left[ 0,1\right] \ \text{and\ }u\in L^p(\widehat{R}_i,1-\widehat{R%
}_i).
\label{def}
\end{equation}%
Here and in the following, we denote a constant function by its value. Therefore, in the expression for $\widehat{E}$ above, $\widehat{\sigma}$ denotes the function in $C([0,1])$ with constant value  $\widehat{\sigma}$. Finally, in (\ref{Fhat}) we can write
\begin{eqnarray*}
\int\limits_{0}^{1}\widehat{E}\left( x\right) ^{2}dx &=&\left\langle %
\widehat{\sigma} +\sum\limits_{i=1}^{N}Z_ic_{i}\mathcal{L}_{i}(\widehat{n}%
_{i}),\widehat{\sigma} +\sum\limits_{i=1}^{N}Z_ic_{i}\mathcal{L}_{i}(\widehat{%
n}_{i})\right\rangle,
\end{eqnarray*}%
where $\left\langle \cdot,\cdot\right\rangle$ denotes the scalar product 
\begin{equation*}
\left\langle f,g\right\rangle :=\int\limits_{0}^{1}f\left( x\right) g\left(
x\right) dx,\ f,g\in C\left([0,1]\right),
\end{equation*}
and then
$$
\widehat{F}^{\mathrm{el}}=K\left\langle %
\widehat{\sigma} +\sum\limits_{i=1}^{N}Z_ic_{i}\mathcal{L}_{i}(\widehat{n}%
_{i}),\widehat{\sigma} +\sum\limits_{i=1}^{N}Z_ic_{i}\mathcal{L}_{i}(\widehat{n}_i)\right\rangle .
$$
In the following, the space $C([0,1])$ is equipped with its usual maximum norm
$$
\Vert f\Vert_\infty=\max\limits_{x\in[0,1]}\vert f(x) \vert,\ f\in C([0,1]).
$$ 

We conclude that 
\begin{eqnarray*}
\widehat{F}=\widehat{F}^{\mathrm{en}}+\widehat{F}^{\mathrm{el}}=\sum\limits_{i=1}^{N}c_{i}\left\langle \nu\circ \widehat{n}%
_i\right\rangle_i+K\left\langle %
\widehat{\sigma} +\sum\limits_{i=1}^{N}Z_ic_{i}\mathcal{L}_{i}(\widehat{n}%
_{i}),\widehat{\sigma} +\sum\limits_{i=1}^{N}Z_ic_{i}\mathcal{L}_{i}(\widehat{n}_i)\right\rangle .
\end{eqnarray*}
 Moreover, the electroneutrality condition (\ref{electronhat}) reads 
\begin{equation*}
2\widehat{\sigma}+\sum\limits_{i=1}^{N}Z_ic_{i}[\mathcal{L}_{i}(\widehat{n}%
_{i})](1)=0
\end{equation*}
by using (\ref{integral}).

Now, we introduce:

\begin{itemize}
\item the product space $$X:=\prod\limits_{i=1}^N L^p(\widehat{R}_i,1-\widehat{R%
}_i);$$
each element $U\in X$ has $N$ components $U_i\in L^p(\widehat{R}_i,1-\widehat{R%
}_i)$, $i\in\{1,\ldots,N\}$.
\item the non-linear functional
$$
\widehat{F}:X_+:=\{U\in X:U_i\geq 0\text{\ for all $i\in\{1,\ldots,N\}$}\}\rightarrow \mathbb{R}
$$
given by 
\begin{equation}
\widehat{F}(U)=\sum\limits_{i=1}^{N}c_{i}\left\langle \nu\circ U
_i\right\rangle_i+K\left\langle %
\widehat{\sigma} +\mathcal{T}(U),\widehat{\sigma} +\mathcal{T}(U)\right\rangle,\ U\in X_+, \label{FhatU}
\end{equation}
where $\mathcal{T}:X\rightarrow C([0,1])$ is the linear operator
\begin{equation}
	\mathcal{T}(U)=\sum\limits_{i=1}^{N}Z_ic_{i}\mathcal{L}_{i}(U_i),\ U\in X; \label{loT}
\end{equation}
\item the affine functional $E:X\rightarrow \mathbb{R}$ given by%
\begin{equation*}
E(U)=2\widehat{\sigma }+[\mathcal{T}(U)](1),\ U\in X.
\end{equation*}
\end{itemize}

The minimization problem is:

\begin{itemize}
\item[MP] \textit{minimize $\widehat{F}\left( U\right)$ with respect to $U\in X$, subject to the \emph{non-negativity condition}%
\begin{equation}
U_{i}\geq 0,\ i\in\{1,\ldots ,N\}, \label{nncond}
\end{equation}
and to the \emph{electroneutrality condition}
\begin{equation}
E(U)=0.  \label{encond0}
\end{equation}}
\end{itemize}

Note that, here, the nondimensional concentrations are denoted by $U_i$ rather
than by $\widehat{n}_i$. Moreover, observe that the non-negativity condition (\ref{nncond}) means
\begin{equation}
U_{i}(x)\geq 0\text{ for almost every $x\in[\widehat{R}_i,1-\widehat{R}_i]$},\ i\in\{1,\ldots ,N\}, \label{nncond2}
\end{equation}
since $U_i$ is defined up to a set of measure zero. Similar shorthand from statements like (\ref{nncond2}) to the form (\ref{nncond}) will be used below. Finally, observe that the nondimensional free energy $\widehat{F}$ in (\ref{FhatU}) is split into the entropic part
\begin{equation}
	\widehat{F}^{\mathrm{en}}(U)=\sum\limits_{i=1}^N c_i\langle\nu\circ U_i\rangle_i,\ U\in X_+, \label{FenhatU}
\end{equation}
and the electrostatic part
\begin{equation}
	\widehat{F}^{\mathrm{el}}(U)=K\left\langle %
	\widehat{\sigma} +\mathcal{T}(U),\widehat{\sigma} +\mathcal{T}(U)\right\rangle,\ U\in X_+. \label{FelhatU}
\end{equation}

\subsection{Facts about the problem MP}\label{Facts}

Concerning the problem MP, the following facts, reported in the next propositions, hold.

The next proposition describes properties of the function $\nu$ defined in (\ref{nu}) and appearing in the entropic part (\ref{FenhatU}) of the nondimensional free energy $\widehat{F}$.

\begin{pro}\label{firstfact}
	$\nu$ is a non-negative strictly convex function, strictly decreasing on the interval $[0,1]$ from the value $1$ at 0 to its minimum value $0$ at $1$, and strictly increasing and diverging to $+\infty$ on the interval $ [1,+\infty)$. Moreover, it has an infinite negative derivative at $0$. Finally, for any $\alpha>1$, there exists $M_\alpha\geq 1$ such that $\nu$ is strictly increasing in $[M_\alpha,+\infty)$ and 
	$$
	\nu(u)\leq u^\alpha,\ u\geq M_\alpha.
	$$ 
\end{pro}

\begin{proof}
We have $\nu(0)=1$ and $\lim\limits_{u\rightarrow +\infty}\nu(u)=+\infty$, since
$$
\lim\limits_{u\rightarrow 0^+}u\log u=0\text{\ and\ }\lim\limits_{u\rightarrow +\infty}u\log u=+\infty.
$$
For $u>0$, we have 
\begin{equation*}
	\nu^{\prime }\left( u\right) =\log u-1+u\cdot\frac{1}{u}=\log u.
\end{equation*}%
Therefore, $\nu$ is strictly decreasing on $[0,1]$ and strictly increasing on $[1,+\infty)$. Moreover, for $u>0$, we have
\begin{equation*}
	\nu^{\prime \prime }\left( u\right) =\frac{1}{u}>0
\end{equation*}
and then $\nu$ is strictly convex. Finally, for any $\alpha>1$, we have
$$
\lim\limits_{u\rightarrow +\infty}\frac{\nu(u)}{u^\alpha}=0
$$
and then the last part of the proposition follows.
\end{proof}

The next proposition shows that the entropic term (\ref{FenhatU}) is well defined and explains why we cannot take $p=1$ in defining the space $X$.

\begin{pro} \label{Lp-L1}
	
For any $i\in\{1,\ldots,N\}$, $u\mapsto \nu\circ u$ is a map
	$$
	\left\{u\in L^p(\widehat{R}_i,1-\widehat{R}_i):u\geq 0\right\}\rightarrow L^1(\widehat{R}_i,1-\widehat{R}_i),
	$$
	but it is not a map
	$$
		\left\{u\in L^1(\widehat{R}_i,1-\widehat{R}_i):u\geq 0\right\}\rightarrow L^1(\widehat{R}_i,1-\widehat{R}_i).
	$$
\end{pro}
\begin{proof}
	Let $M_p$ be the number $M_{\alpha}$, for $\alpha=p$, defined in Proposition \ref{firstfact}. For $u\in L^p(\widehat{R}_i,1-\widehat{R}_i)$ such that $u\geq 0$, we have
	\begin{eqnarray*}
	&&\int\limits_{\widehat{R}_i}^{1-\widehat{R}_i}\nu(u(x))dx\\
	&&=\int\limits_{\{x\in[\widehat{R}_i,1-\widehat{R}_i]:u(x)\leq M_p\}}\nu(u(x))dx+\int\limits_{\{x\in[\widehat{R}_i,1-\widehat{R}_i]:u(x)> M_p\}}\nu(u(x))dx\\
	\\
	&&\leq (1-2\widehat{R}_i)\max\limits_{u\in[0,M_p]}\nu(u)+\int\limits_{\{x\in[\widehat{R}_i,1-\widehat{R}_i]:u(x)> M_p\}}u(x)^pdx\\
	&&\leq (1-2\widehat{R}_i)\max\limits_{u\in[0,M_p]}\nu(u)+\int\limits_{\widehat{R}_i}^{1-\widehat{R}_i}u(x)^pdx<+\infty.
	\end{eqnarray*}
Thus, $\nu\circ u\in L^1(\widehat{R}_i,1-\widehat{R}_i)$.

For the second part, consider a sequence $(I_n)_{n\geq n_0}$ of disjoint subintervals of $[\widehat{R}_i,1-\widehat{R}_i]$ of Lebesgue measure $m( I_n)=\frac{1}{\mathrm{e}^n n^2}$, where $n_0\in\mathbb{N}^+$ is such that
$$
\sum\limits_{n=n_0}^{\infty}m(I_n)=\sum\limits_{n=n_0}^{\infty}\frac{1}{\mathrm{e}^n n^2}<1-2\widehat{R}_i,
$$
and let 
$$
u=\sum\limits_{n=n_0}^\infty \mathrm{e}^n\chi_{I_n},
$$
where $\chi_{I_n}$ denotes the characteristic function of $I_n$. We have
$$
\int\limits_{\widehat{R}_i}^{1-\widehat{R}_i}u(x)dx=\sum\limits_{n=n_0}^\infty \int\limits_{x\in I_n}u(x)dx=\sum\limits_{n=n_0}^\infty \mathrm{e}^n m(I_n)=\sum\limits_{n=n_0}^\infty \frac{1}{n^2}<+\infty
$$
and then $u\in L^1(\widehat{R}_i,1-\widehat{R}_i)$.  On the other hand, we have
$$
\int\limits_{\widehat{R}_i}^{1-\widehat{R}_i}\nu(u(x))dx=\int\limits_{\widehat{R}_i}^{1-\widehat{R}_i}u(x)\log u(x)dx-\int\limits_{\widehat{R}_i}^{1-\widehat{R}_i}u(x)dx+1-2\widehat{R}_i
$$
with
\begin{eqnarray*}
\int\limits_{\widehat{R}_i}^{1-\widehat{R}_i}u(x)\log u(x)dx&=&\sum\limits_{n=n_0}^\infty \int\limits_{x\in I_n}u(x)\log u(x)dx=\sum\limits_{n=n_0}^\infty \mathrm{e}^n\log \mathrm{e}^n m(I_n)\\
&=&\sum\limits_{n=n_0}^\infty \frac{1}{n}=+\infty
\end{eqnarray*}
and then $\nu\circ u\notin L^1(\widehat{R}_i,1-\widehat{R}_i)$.

\end{proof}

 The next proposition provides a single integral representation for the operators $\mathcal{L}_i$, $i\in\{1,\ldots,N\}$, appearing in the electrostatic part (\ref{FelhatU}), with (\ref{loT}), of $\widehat{F}$.

 \begin{pro} \label{uniqueint}
 	For $i\in\{1,\ldots,N\}$, we have
 	\begin{eqnarray}
 		&&[\mathcal{L}_i(u)](x)=\frac{1}{2\widehat{R}_i}\int\limits_{\widehat{R}_i}^{\min\{x+\widehat{R}_i,1-\widehat{R}_i\}}u(s)\min\{x+\widehat{R}_i-s,2\widehat{R}_i\}ds,\notag \\
 		&&x\in \left[ 0,1\right] \ \text{and\ }u\in L^p(\widehat{R}_i,1-\widehat{R%
 		}_i). \notag\\
 		\label{def2}
 	\end{eqnarray}
 \end{pro}
 
 \begin{proof}
 	In (\ref{def}), we integrate over the region in the plane $(z,s)$ where $z$ varies from $0$ to $x$ and $s$ is between the two lines $s=z-\widehat{R}_i$ and $s=z+\widehat{R}_i$. This means that $s$ varies from $-\widehat{R}_i$ to $x+\widehat{R}_i$ and $z$ is between the lines $z=s-\widehat{R}_i$ and $z=s+\widehat{R}_i$ with the constraint $z\in [0,x]$.
 	
 	Thus, by swapping the two integrals in (\ref{def}), we obtain, since $u$ is zero outside $[\widehat{R}_i,1-\widehat{R}_i]$, 
 	\begin{eqnarray*}
 	\int\limits_{0}^{x}\frac{1}{2\widehat{R}_i}\int\limits_{z-\widehat{R}_{i}}^{z+\widehat{R}_{i}}u(s)dsdz&=&\frac{1}{2\widehat{R}_i}%
 	\int\limits_{-\widehat{R}_i}^{x+\widehat{R}_i}\int\limits_{\max\{0,s-\widehat{R}_i\}}^{\min\{x,s+\widehat{R}_i\}}u(s)dzds\\
 	&=&\frac{1}{2\widehat{R}_i}%
 	\int\limits_{\widehat{R}_i}^{\min\{x+\widehat{R}_i,1-\widehat{R}_i\}}\int\limits_{s-\widehat{R}_{i}}^{\min\{x,s+\widehat{R}_i\}}u\left(s\right)dzds\\
 &=&\frac{1}{2\widehat{R}_i}\int\limits_{\widehat{R}_i}^{\min\{x+\widehat{R}_i,1-\widehat{R}_i\}}u(s)\left(\min\{x,s+\widehat{R}_i\}-(s-\widehat{R}_i)\right)ds\\
 &=&\frac{1}{2\widehat{R}_i}\int\limits_{\widehat{R}_i}^{\min\{x+\widehat{R}_i,1-\widehat{R}_i\}}u(s)\min\{x+\widehat{R}_i-s,2\widehat{R}_i\}ds.
 	\end{eqnarray*}  
 \end{proof}
In particular, by Proposition \ref{uniqueint}, we have, for $i\in\{1,\ldots,N\}$,
\begin{eqnarray*}
	&&[\mathcal{L}_i(u)](1)=\int\limits_{\widehat{R}_i}^{1-\widehat{R}_i}u(s)ds,\ 
u\in L^p(\widehat{R}_i,1-\widehat{R%
	}_i)
\end{eqnarray*}
and then
\begin{equation}
[\mathcal{T}(U)](1)=\sum\limits_{i=1}^N Z_i c_i \int\limits_{\widehat{R}_i}^{1-\widehat{R}_i}U_i(s)ds,\ U\in X.  \label{Tu1}
\end{equation}
Therefore, the electroneutrality condition (\ref{encond0}) reads
\begin{equation}
\sum\limits_{i=1}^N Z_i c_i \int\limits_{\widehat{R}_i}^{1-\widehat{R}_i}U_i(s)ds=-2\widehat{\sigma}. \label{encond}
\end{equation}
Moreover, observe that
\begin{equation}
\Vert \mathcal{L}_i(u)\Vert_\infty \leq \Vert u\Vert_1,\ u\in L^p(\widehat{R}_i,1-\widehat{R%
}_i).  \label{operatornorm}
\end{equation}

The expressions (\ref{def}) and (\ref{def2}) for $[\mathcal{L}_i(u)](x)$ are not defined for $\widehat{R}_i=0$, since an indeterminate form $\frac{0}{0}$ appears in such expressions. In the case $\widehat{R}_i=0$, we set
\begin{equation}
[\mathcal{L}_i(u)](x)=\int\limits_{0}^{x}u(s)ds,\ x\in \left[ 0,1\right] \ \text{and\ }u\in L^p(0,1).  \label{R=0}
\end{equation}
The next proposition shows that (\ref{R=0}) is the limiting case of (\ref{def}) and   (\ref{def2}), as $\widehat{R}_i\rightarrow 0$.  Here, we denote by $\mathcal{L}_i^{\widehat{R}_i}$ and $\mathcal{L}_i^{0}$ the linear operators $\mathcal{L}_i$ corresponding to the radii $\widehat{R}_i$ and $\widehat{R}_i=0$, respectively.

\begin{pro} \label{tozero}
For any $u\in L^p(0,1)$, we have
$$
\left\Vert \mathcal{L}_i^{\widehat{R}_i}\left(u|_{\left[\widehat{R}_i,1-\widehat{R}_i\right]}\right)-\mathcal{L}_i^{0}u\right\Vert_\infty \leq C(\widehat{R}_i)\left\Vert u\right\Vert_p,  
$$
with $C(\widehat{R}_i)\rightarrow 0$ as $\widehat{R}_i\rightarrow 0$.
\end{pro}

\begin{proof}
	Let $x\in[0,1]$. We have
	\begin{eqnarray*}
	&&\left[\mathcal{L}_i^{\widehat{R}_i}\left(u|_{\left[\widehat{R}_i,1-\widehat{R}_i\right]}\right)\right](x)-[\mathcal{L}_i^{0}(u)](x)\\
	&&=\frac{1}{2\widehat{R}_i}\int\limits_{\widehat{R}_i}^{\min\{x+\widehat{R}_i,1-\widehat{R}_i\}}u(s)\min\{x+\widehat{R}_i-s,2\widehat{R}_i\}ds-\int\limits_{0}^{x}u(s)ds.
	\end{eqnarray*}
	For $x\leq 2\widehat{R}_i$, we have
		\begin{eqnarray*}
		&&\left\vert \left[\mathcal{L}_i^{\widehat{R}_i}\left(u|_{\left[\widehat{R}_i,1-\widehat{R}_i\right]}\right)\right](x)-[\mathcal{L}_i^{0}(u)](x)\right\vert\\
		&&\leq \left\vert \left[\mathcal{L}_i^{\widehat{R}_i}\left(u|_{\left[\widehat{R}_i,1-\widehat{R}_i\right]}\right)\right](x)\right\vert+\left\vert[\mathcal{L}_i^{0}(u)](x)\right\vert\\
		&&\leq\int\limits_{\widehat{R}_i}^{\min\{x+\widehat{R}_i,1-\widehat{R}_i\}}\vert u(s)\vert ds+\int\limits_{0}^{x}\vert u(s)\vert ds\\
		&&\leq \int\limits_{\widehat{R}_i}^{\min\{3\widehat{R}_i,1-\widehat{R}_i\}}\vert u(s)\vert ds+\int\limits_{0}^{2\widehat{R}_i}\vert u(s)\vert ds\\
		&&\leq \left(\min\{2\widehat{R}_i,1-2\widehat{R}_i\}\right)^{\frac{1}{q}}\Vert u\Vert_p+\left(2\widehat{R}_i\right)^{\frac{1}{q}}\Vert u\Vert_p\leq 2 \left(2\widehat{R}_i\right)^{\frac{1}{q}}\Vert u\Vert_p,
	\end{eqnarray*}
where \(q\) is the conjugate exponent of \(p\), namely
	\[
	\frac1p+\frac1q=1.
	\]
	For $2\widehat{R}_i\leq x\leq 1-2\widehat{R}_i$, we have
	\begin{eqnarray*}
		&&\left[\mathcal{L}_i^{\widehat{R}_i}\left(u|_{\left[\widehat{R}_i,1-\widehat{R}_i\right]}\right)\right](x)-[\mathcal{L}_i^{0}(u)](x)\\
		&&=\frac{1}{2\widehat{R}_i}\int\limits_{\widehat{R}_i}^{x+\widehat{R}_i}u(s)\min\{x+\widehat{R}_i-s,2\widehat{R}_i\}ds-\int\limits_{0}^{x}u(s)ds\\
		&&=\frac{1}{2\widehat{R}_i}\int\limits_{\widehat{R}_i}^{x-\widehat{R}_i}u(s)2\widehat{R}_ids+\frac{1}{2\widehat{R}_i}\int\limits_{x-\widehat{R}_i}^{x+\widehat{R}_i}u(s)(x+\widehat{R}_i-s)ds-\int\limits_{0}^{x}u(s)ds\\
		&&=\frac{1}{2\widehat{R}_i}\int\limits_{x-\widehat{R}_i}^{x+\widehat{R}_i}u(s)(x+\widehat{R}_i-s)ds-\int\limits_{0}^{\widehat{R}_i}u(s)ds-\int\limits_{x-\widehat{R}_i}^{x}u(s)ds.\\
	\end{eqnarray*}
	Therefore,
	\begin{eqnarray*}
		&&\left\vert \left[\mathcal{L}_i^{\widehat{R}_i}\left(u|_{\left[\widehat{R}_i,1-\widehat{R}_i\right]}\right)\right](x)-[\mathcal{L}_i^{0}(u)](x)\right\vert \\
		&&\leq \int\limits_{x-\widehat{R}_i}^{x+\widehat{R}_i}\vert u(s)\vert ds+\int\limits_{0}^{\widehat{R}_i}\vert u(s)\vert ds+\int\limits_{x-\widehat{R}_i}^{x}\vert u(s)\vert ds\\
		&&\leq (2\widehat{R}_i)^{\frac{1}{q}}\Vert u\Vert_p+(\widehat{R}_i)^{\frac{1}{q}}\Vert u\Vert_p+(\widehat{R}_i)^{\frac{1}{q}}\Vert u\Vert_p\\
		&&\leq 3(2\widehat{R}_i)^{\frac{1}{q}}\Vert u\Vert_p.
	\end{eqnarray*}
	For $x\geq \max\{2\widehat{R}_i,1-2\widehat{R}_i\}$, we have
	\begin{eqnarray*}
		&&\left[\mathcal{L}_i^{\widehat{R}_i}\left(u|_{\left[\widehat{R}_i,1-\widehat{R}_i\right]}\right)\right](x)-[\mathcal{L}_i^{0}(u)](x)\\
		&&=\frac{1}{2\widehat{R}_i}\int\limits_{\widehat{R}_i}^{1-\widehat{R}_i}u(s)\min\{x+\widehat{R}_i-s,2\widehat{R}_i\}ds-\int\limits_{0}^{x}u(s)ds\\
		&&=\frac{1}{2\widehat{R}_i}\int\limits_{\widehat{R}_i}^{x-\widehat{R}_i}u(s)2\widehat{R}_ids+\frac{1}{2\widehat{R}_i}\int\limits_{x-\widehat{R}_i}^{1-\widehat{R}_i}u(s)(x+\widehat{R}_i-s)ds-\int\limits_{0}^{x}u(s)ds\\
		&&=\frac{1}{2\widehat{R}_i}\int\limits_{x-\widehat{R}_i}^{1-\widehat{R}_i}u(s)(x+\widehat{R}_i-s)ds-\int\limits_{0}^{\widehat{R}_i}u(s)ds-\int\limits_{x-\widehat{R}_i}^{x}u(s)ds.\\
	\end{eqnarray*}
	Therefore,
	\begin{eqnarray*}
		&&\left\vert \left[\mathcal{L}_i^{\widehat{R}_i}\left(u|_{\left[\widehat{R}_i,1-\widehat{R}_i\right]}\right)\right](x)-[\mathcal{L}_i^{0}(u)](x)\right\vert \\
		&&\leq \int\limits_{x-\widehat{R}_i}^{1-\widehat{R}_i}\vert u(s)\vert ds+\int\limits_{0}^{\widehat{R}_i}\vert u(s)\vert ds+\int\limits_{x-\widehat{R}_i}^{x}\vert u(s)\vert ds\\
		&&\leq (1-x)^{\frac{1}{q}}\Vert u\Vert_p+(\widehat{R}_i)^{\frac{1}{q}}\Vert u\Vert_p+(\widehat{R}_i)^{\frac{1}{q}}\Vert u\Vert_p\\
		&&\leq (2\widehat{R}_i)^{\frac{1}{q}}\Vert u\Vert_p+(\widehat{R}_i)^{\frac{1}{q}}\Vert u\Vert_p+(\widehat{R}_i)^{\frac{1}{q}}\Vert u\Vert_p\\
		&&\leq 3(2\widehat{R}_i)^{\frac{1}{q}}\Vert u\Vert_p.
	\end{eqnarray*}
	
	We conclude that
	$$
	\left\Vert \mathcal{L}_i^{\widehat{R}_i}\left(u|_{\left[\widehat{R}_i,1-\widehat{R}_i\right]}\right)-\mathcal{L}_i^{0}u\right\Vert_\infty \leq C(\widehat{R}_i)\left\Vert u\right\Vert_p,  
	$$
	with $C(\widehat{R}_i)=3(2\widehat{R}_i)^{\frac{1}{q}}$.

\end{proof}

Observe that the linear operator $\mathcal{L}_i$, $i\in\{1,\ldots,N\}$,
and hence the electrostatic term (\ref{FelhatU}) can also be defined as acting on functions in
$L^1(\widehat{R}_i,1-\widehat{R}_i)$. However, when the electrostatic term is combined with the entropic term (\ref{FenhatU}) in the free energy $\widehat{F}$, it is necessary to restrict its action to functions
in $L^p(\widehat{R}_i,1-\widehat{R}_i)$, with $p>1$. Moreover, the proof of
the continuity property in the previous Proposition \ref{tozero} also requires working in a space $L^p(\widehat{R}_i,1-\widehat{R}_i)$ with $p>1$.

\section{Existence and uniqueness of the minimum}\label{EUsection}

In this section, we show that the problem MP has a unique minimum $\overline{U}=(\overline{U}_1,\ldots,\overline{U}_N)$. The background material for this section is in Appendix A. 

We consider the space $X$ equipped with the norm
$$
\Vert U\Vert=\max\{\Vert U_1\Vert_p,\ldots,\Vert U_N\Vert_p\},\ U\in X,
$$
and, with reference to the appendix, we consider 
\begin{equation}
A=\{U\in X: U_i\geq 0,\ i\in\{1,\ldots,N\},\text{ and }E(U)=0\}  \label{A}
\end{equation}
as the domain of the functional $\widehat{F}$.

Note that $A\neq\emptyset$. Indeed, for $\widehat{\sigma}\neq 0$, let $U\in X$ have all components equal to zero except for one constant component, of value $k$, corresponding to a counter-ion type, where $k$ is chosen to satisfy the electroneutrality condition $E(U)=0$. Then $U$ belongs to $A$. For $\widehat{\sigma}=0$, $U=0$ belongs to $A$.

The reason we consider the functional $\widehat{F}$ as acting on elements of the space $X$ is that such a space possesses the Banach–Saks property (see Theorem \ref{BSC} in Appendix A) and, as a consequence, Theorem \ref{existence2} in Appendix A for the existence of a minimum of $\widehat{F}$ can be applied. Such a theorem requires closedness and convexity for the domain $A$ and lower semicontinuity and convexity for the functional $\widehat{F}$.  Moreover, for the uniqueness of the minimum, the strict convexity of $\widehat{F}$ is required (see Theorem \ref{uniqueness} in Appendix A). 

\subsection{The hypotheses for the existence and uniqueness}

The next three propositions state that $A$ is closed and convex and $\widehat{F}$ is continuous and strictly convex. In the following, $m(B)$ will denote the Lebesgue measure of a measurable subset $B$ of $\mathbb{R}$. 

\begin{pro}
	$A$ is closed and convex.
\end{pro}

\begin{proof}
We prove that $A$ is closed. Let $(U^{(n)})_{n\in\mathbb{N}^+}$ be a sequence in $A$ such that $U^{(n)}\rightarrow U$, $n\rightarrow \infty$, in the space $X$, i.e.,
\begin{equation*}
	\lim\limits_{n\rightarrow \infty}\Vert U^{(n)}_i-U_i\Vert_p=0,\ i\in\{1,\ldots,N\}.	
\end{equation*}
For any $i\in\{1,\ldots,N\}$, the $L^p$-convergence of $U_i^{(n)}$ to $U_i$ implies the convergence in measure
$$
\forall \gamma >0: \lim\limits_{n\rightarrow \infty}m\left(\{x\in[\widehat{R}_i,1-\widehat{R}_i]:\vert U_i^{(n)}(x)-U_i(x)\vert>\gamma \}\right)=0.
$$
Fix $i\in\{1,\ldots,N\}$ and $\gamma>0$. We have, for all $n\in\mathbb{N}^+$, 
\begin{eqnarray*}
&&\{x\in[\widehat{R}_i,1-\widehat{R}_i]:U_i(x)<-\gamma\text{\ and\ }U^{(n)}_i(x)\geq 0 \}\\
&&\subseteq \{x\in[\widehat{R}_i,1-\widehat{R}_i]:\vert U_i^{(n)}(x)-U_i(x)\vert>\gamma \}.
\end{eqnarray*} 
In fact, for $x\in[\widehat{R}_i,1-\widehat{R}_i]$ such that $U_i(x)<-\gamma$ and $U^{(n)}_i(x)\geq 0$, we have
$$
\vert U_i^{(n)}(x)-U_i(x)\vert\geq U_i^{(n)}(x)-U_i(x)>0-(-\gamma)=\gamma.
$$
Thus, for all $n\in\mathbb{N}^+$,
\begin{eqnarray*}
&&m(\{x\in[\widehat{R}_i,1-\widehat{R}_i]:U_i(x)<-\gamma\})\\
&&=m(\{x\in[\widehat{R}_i,1-\widehat{R}_i]:U_i(x)<-\gamma\text{\ and\ }U^{(n)}_i(x)\geq 0 \})\\
&&\leq m(\{x\in[\widehat{R}_i,1-\widehat{R}_i]:\vert U_i^{(n)}(x)-U_i(x)\vert>\gamma \}).
\end{eqnarray*}
By passing to the limit, as $n\rightarrow \infty$, we obtain
$$
m(\{x\in[\widehat{R}_i,1-\widehat{R}_i]:U_i(x)<-\gamma\})=0.
$$
Since
$$
\forall \gamma >0: m(\{x\in[\widehat{R}_i,1-\widehat{R}_i]:U_i(x)<-\gamma\})=0,
$$
 we have $U_i\geq 0$.
 
 Moreover, for all $n\in\mathbb{N}^+$, we have
 $$
 E(U^{(n)})=2\widehat{\sigma}+[\mathcal{T}(U^{(n)})](1)=0.
 $$
 As $n\rightarrow \infty$, by the boundedness of the linear operator $\mathcal{T}$, we obtain 
 $$
 E(U)=2\widehat{\sigma}+[\mathcal{T}(U)](1)=0.
 $$
We conclude that $U\in A$ and then $A$ is closed.
 
 Now, we prove that $A$ is convex. For $U,V\in A$ and $\alpha\in[0,1]$, we have
 $$
 (1-\alpha)U+\alpha V\geq (1-\alpha)\cdot 0+\alpha \cdot 0=0. 
 $$ 
 Moreover, by the linearity of $\mathcal{T}$, we obtain
\begin{eqnarray*}
 E((1-\alpha)U+\alpha V)&=&2\widehat{\sigma}+[\mathcal{T}((1-\alpha)U+\alpha V)](1)\\
 &=&2\widehat{\sigma}+(1-\alpha)[\mathcal{T}(U)](1)+\alpha[\mathcal{T}(V)](1)\\
 &=&(1-\alpha)(2\widehat{\sigma}+[\mathcal{T}(U)](1))+\alpha(2\widehat{\sigma}+[\mathcal{T}(V)](1))\\
 &=&(1-\alpha)E(U)+\alpha E(V)=0.
 \end{eqnarray*}
 We conclude that $(1-\alpha)U+\alpha V\in A$. Therefore, $A$ is convex.
\end{proof}

\begin{pro}
	$\widehat{F}$ is continuous.
\end{pro}

\begin{proof}
	Let $U\in A$ and let $(U^{(n)})_{n\in\mathbb{N}^+}$ be a sequence in $A$ such that $U^{(n)}\rightarrow U$, $n\rightarrow \infty$, in the space $X$, i.e.,
	\begin{equation*}
		\lim\limits_{n\rightarrow \infty}\Vert U^{(n)}_i-U_i\Vert_p=0,\ i\in\{1,\ldots,N\}.	
	\end{equation*}
	We have, for any $n\in\mathbb{N}^+$, 
	\begin{eqnarray*}
		\widehat{F}(U^{(n)})-\widehat{F}(U)&=&\widehat{F}^{\mathrm{en}}(U^{(n)})-\widehat{F}^{\mathrm{en}}(U)\\
		&&+\widehat{F}^{\mathrm{el}}(U^{(n)})-\widehat{F}^{\mathrm{el}}(U)\\ &=&\sum\limits_{i=1}^{N}c_{i}\int\limits_{\widehat{R}_i}^{1-\widehat{R}_i}\left(\nu(U^{(n)}_i(x))-\nu(U_i(x))\right)dx\\ &&+K\left(\left\Vert 	\widehat{\sigma} +\mathcal{T}(U^{(n)})\right\Vert_2^2-\left\Vert 	\widehat{\sigma} +\mathcal{T}(U)\right\Vert_2^2\right) 
	\end{eqnarray*}
	(recall (\ref{FhatU}), (\ref{FenhatU}) and (\ref{FelhatU})). Here, $\Vert\ \cdot\ \Vert_2$ denotes the $L^2(0,1)$-norm 
	$$
	\Vert f\Vert^2_2=\left\langle f,f\right\rangle,\ f\in C([0,1]).
	$$
	We show that
	\begin{equation*}
		\widehat{F}^{\mathrm{en}}(U^{(n)})\rightarrow \widehat{F}^{\mathrm{en}}(U)\text{\ and\ }\widehat{F}^{\mathrm{el}}(U^{(n)})\rightarrow\widehat{F}^{\mathrm{el}}(U), n\rightarrow \infty.
	\end{equation*}

	Regarding the electrostatic term $\widehat{F}^{\mathrm{el}}$, the boundedness of the linear operator $\mathcal{T}$ immediately implies
	$$
	\lim\limits_{n\rightarrow \infty}\left\Vert 	\widehat{\sigma} +\mathcal{T}(U^{(n)})\right\Vert_2=\left\Vert 	\widehat{\sigma} +\mathcal{T}(U)\right\Vert_2.
	$$
	Hence, we have
	$
	\widehat{F}^{\mathrm{el}}(U^{(n)})\rightarrow \widehat{F}^{\mathrm{el}}(U),\ n\rightarrow +\infty.
	$

	Now, we consider the entropic term $\widehat{F}^{\mathrm{en}}$. We prove that, for $i\in\{1,\ldots,N\}$, we have
	\begin{equation}
		\lim\limits_{n\rightarrow \infty}\int\limits_{\widehat{R}_i}^{1-\widehat{R}_i}\left\vert \nu(U^{(n)}_i(x))-\nu(U_i(x))\right\vert dx=0 \label{conv}.
	\end{equation}
	In other words, with reference to Proposition  \ref{Lp-L1}, we prove that the map
	$$
	\left\{u\in L^p(\widehat{R}_i,1-\widehat{R}_i):u\geq 0\right\}\rightarrow L^1(\widehat{R}_i,1-\widehat{R}_i),\ u\mapsto \nu\circ u,
	$$
	is continuous.
	
	Let $\alpha\in (1,p)$ and let $M_\alpha\geq 1$ be as in Proposition \ref{firstfact}. We have $\nu$  increasing on $[M_\alpha,+\infty)$ and 
	$$
	\nu(u)\leq u^\alpha,\ u\geq M_\alpha.
	$$
	Let $M\geq \max\{M_\alpha,\mathrm{e}\}$. For any $n\in\mathbb{N}^+$, we partition the interval $[\widehat{R}_i,1-\widehat{R}_i]$ as $[\widehat{R}_i,1-\widehat{R}_i]=B^{(n)}\cup C^{(n)}$, where
	\begin{equation*}
		B^{(n)}=\left\{x\in[\widehat{R}_i,1-\widehat{R}_i]:\max\left\{U_i^{(n)}(x),U_i(x)\right\}>M\right\}\\
	\end{equation*}
	and
	\begin{equation*}
		C^{(n)}=\left\{x\in[\widehat{R}_i,1-\widehat{R}_i]:\max\left\{U_i^{(n)}(x),U_i(x)\right\}\leq M\right\}\\
	\end{equation*}
	and write
	\begin{eqnarray*}
		\int\limits_{\widehat{R}_i}^{1-\widehat{R}_i}\left\vert \nu(U^{(n)}_i(x))-\nu(U_i(x))\right\vert dx &=&\int\limits_{x\in B^{(n)}}\left\vert \nu(U^{(n)}_i(x))-\nu(U_i(x))\right\vert dx\\
		&&+\int\limits_{x\in C^{(n)}}\left\vert \nu(U^{(n)}_i(x))-\nu(U_i(x))\right\vert dx.
	\end{eqnarray*}
	We have
	\begin{eqnarray*}
		\int\limits_{x\in B^{(n)}}\left\vert \nu(U^{(n)}_i(x))-\nu(U_i(x))\right\vert dx&\leq& \int\limits_{x\in B^{(n)}}\max\left\{\nu(U^{(n)}_i(x)),\nu\left(U_i(x)\right)\right\}dx\\
		&=&\int\limits_{x\in B^{(n)}}\nu\left(\max\left\{U^{(n)}_i(x),U_i(x)\right\}\right)dx\\
		&\leq &\int\limits_{x\in B^{(n)}}\max\left\{U^{(n)}_i(x),U_i(x)\right\}^\alpha dx\\
		&=&\int\limits_{x\in B^{(n)}}\frac{\max\left\{U^{(n)}_i(x),U_i(x)\right\}^p}{\max\left\{U^{(n)}_i(x),U_i(x)\right\}^{p-\alpha}} dx\\
		&\leq&\frac{1}{M^{p-\alpha}}\int\limits_{x\in B^{(n)}}\max\left\{U^{(n)}_i(x),U_i(x)\right\}^pdx\\
		&=&\frac{1}{M^{p-\alpha}}\int\limits_{x\in B^{(n)}}\max\left\{U^{(n)}_i(x)^p,U_i(x)^p\right\}dx \\
		&\leq&\frac{1}{M^{p-\alpha}}\int\limits_{x\in B^{(n)}}\left(U^{(n)}_i(x)^p+U_i(x)^p\right)dx\\
		&\leq&\frac{1}{M^{p-\alpha}}\left(\sup\limits_{n\in\mathbb{N}^+} \left\Vert U^{(n)}_i\right\Vert_p^p+\left\Vert U_i \right\Vert_p^p \right).
	\end{eqnarray*}

	Now, fix $\varepsilon>0$ and let $M\geq\max\{M_\alpha,\mathrm{e}\}$ be such that
	$$
	\frac{1}{M^{p-\alpha}}\left(\sup\limits_{n\in\mathbb{N}^+} \left\Vert U^{(n)}_i\right\Vert_p^p+\left\Vert U_i \right\Vert_p^p \right)\leq \varepsilon.
	$$
	Since $\nu$ is uniformly continuous on $[0,M]$, there exists $\delta>0$ such that, for all $n\in\mathbb{N}^+$,  
	$$
	x\in C^{(n)}\text{\ and\ }\left\vert U^{(n)}_i(x)-U_i(x)\right\vert\leq \delta\ \ \Rightarrow\ \  \left\vert \nu(U^{(n)}_i(x))-\nu(U_i(x))\right\vert \leq \varepsilon. 
	$$
	We partition the complementary set  $C^{(n)}$ of $B^{(n)}$ as $C^{(n)}=D^{(n)}\cup E^{(n)}$, where
	$$
	D^{(n)}=\left\{x\in C^{(n)}:\left\vert U^{(n)}_i(x)-U_i(x)\right\vert\leq \delta\right\}
	$$ 
	and
	$$
	E^{(n)}=\left\{x\in C^{(n)}:\left\vert U^{(n)}_i(x)-U_i(x)\right\vert>\delta\right\}
	$$
	and write
	\begin{eqnarray*}
		\int\limits_{x\in C^{(n)}}\left\vert \nu(U^{(n)}_i(x))-\nu(U_i(x))\right\vert dx&=&\int\limits_{x\in D^{(n)}}\left\vert \nu(U^{(n)}_i(x))-\nu(U_i(x))\right\vert dx\\
		&&+\int\limits_{x\in E^{(n)}}\left\vert \nu(U^{(n)}_i(x))-\nu(U_i(x))\right\vert dx
	\end{eqnarray*}
	with
	$$
	\int\limits_{x\in D^{(n)}}\left\vert \nu(U^{(n)}_i(x))-\nu(U_i(x))\right\vert dx\leq m\left(D^{(n)}\right)\varepsilon \leq (1-2\widehat{R}_i)\varepsilon
	$$
	and
	\begin{eqnarray*}
		\int\limits_{x\in E^{(n)}}\left\vert \nu(U^{(n)}_i(x))-\nu(U_i(x))\right\vert dx&\leq& \int\limits_{x\in E^{(n)}}\max\left\{\nu(U^{(n)}_i(x)),\nu\left(U_i(x)\right)\right\}dx\\
		&\leq &\max\limits_{u\in[0,M]} \nu(u)\cdot m\left(E^{(n)}\right).
	\end{eqnarray*}
	The $L^p$-convergence of $U_i^{(n)}$ to $U_i$ as $n\rightarrow \infty$ implies the convergence in measure
	$$
	\forall \gamma >0: \lim\limits_{n\rightarrow \infty}m\left(\{x\in[\widehat{R}_i,1-\widehat{R}_i]:\vert U_i^{(n)}(x)-U_i(x)\vert>\gamma \}\right)=0.
	$$
	Therefore, taking $\gamma=\delta$, there exists $n_0\in\mathbb{N}^+$ such that, for $n\geq n_0$, we have
	$$
	m\left(E^{(n)}\right)\leq \frac{\varepsilon}{\max\limits_{u\in[0,M]} \nu(u)}.
	$$
	Thus, for $n\geq n_0$, we have
	\begin{eqnarray*}
		\int\limits_{\widehat{R}_i}^{1-\widehat{R}_i}\left\vert \nu(U^{(n)}_i(x))-\nu(U_i(x))\right\vert dx &=&\int\limits_{x\in B^{(n)}}\left\vert \nu(U^{(n)}_i(x))-\nu(U_i(x))\right\vert dx\\
		&&+\int\limits_{x\in D^{(n)}}\left\vert \nu(U^{(n)}_i(x))-\nu(U_i(x))\right\vert dx\\
		&&+\int\limits_{x\in E^{(n)}}\left\vert \nu(U^{(n)}_i(x))-\nu(U_i(x))\right\vert dx\\
		&\leq& \varepsilon+(1-2\widehat{R}_i)\varepsilon+\varepsilon\leq 3\varepsilon.
	\end{eqnarray*}
	This shows that (\ref{conv}) holds. Hence, we have $\widehat{F}^{\mathrm{en}}\left(U^{(n)}\right)\rightarrow \widehat{F}^{\mathrm{en}}(U)$, $n\rightarrow \infty$. 
\end{proof}

\begin{pro}\label{Thconvex}
	$\widehat{F}$ is strictly convex. 
\end{pro} 

\begin{proof}
	By observing that the sum of two
	convex functionals $A\rightarrow \mathbb{R}$ is convex and, if one of them is strictly convex, then the sum is strictly convex, we prove that $\widehat{F}=\widehat{F}^{\mathrm{en}}+\widehat{F}^{\mathrm{el}}$ (recall (\ref{FhatU}), (\ref{FenhatU}) and (\ref{FelhatU})) is strictly convex by showing that $\widehat{F}^{\mathrm{en}}$ is strictly convex and $\widehat{F}^{\mathrm{el}}$ is convex.
	
	We prove the strict convexity of $\widehat{F}^{\mathrm{en}}$. For $U,V\in A$ and $\alpha \in \left[ 0,1\right] $, we have%
	\begin{eqnarray}
		&& \widehat{F}^{\mathrm{en}}\left( \left( 1-\alpha \right) U+\alpha V\right)
		=\sum\limits_{i=1}^{N}c_{i}\int\limits_{\widehat{R}_{i}}^{1-\widehat{R}%
			_{i}}\nu\left( \left( 1-\alpha \right) U_{i}(x)+\alpha V_{i}(x)\right)dx \notag \\
		&&\quad\quad\quad\quad\quad\quad\quad\quad\quad \leq \sum\limits_{i=1}^{N}c_{i}\int\limits_{\widehat{R}_{i}}^{1-\widehat{R}%
			_{i}}\left( \left( 1-\alpha \right) \nu\left( U_{i}\left( x\right) \right)
		+\alpha \nu\left( V_{i}\left( x\right) \right) \right) dx \label{sc}\\
		&&\quad\quad\quad\quad\quad\quad\quad\quad\quad=\left( 1-\alpha \right) \sum\limits_{i=1}^{N}c_{i}\int\limits_{\widehat{R}%
			_{i}}^{1-\widehat{R}_{i}}\nu\left( U_{i}\left( x\right) \right) dx\notag\\
			&&\quad\quad\quad\quad\quad\quad\quad\quad\quad\quad\quad+\alpha
		\sum\limits_{i=1}^{N}c_{i}\int\limits_{\widehat{R}_{i}}^{1-\widehat{R}%
			_{i}}\nu\left( V_{i}\left( x\right) \right) dx \notag \\
		&&\quad\quad\quad\quad\quad\quad\quad\quad\quad=\left( 1-\alpha \right) \widehat{F}^{\mathrm{en}}\left( U\right) +\alpha \widehat{F}^{\mathrm{en}}\left( V\right)  \notag
	\end{eqnarray}
	by the convexity of $\nu$ (used in (\ref{sc})). This shows that $\widehat{F}^{\mathrm{en}}$ is convex.
	
	Moreover, suppose $U\neq V$, i.e., there exists $i\in\{1,\ldots,N\}$ such that $U_{i}\neq
	V_{i}$, i.e., $U_{i}\left( x\right) \neq V_{i}\left( x\right) $ for $x$ in a
	set $E$ of positive measure. For $x\in E$ and $\alpha\in(0,1)$, we have 
	\begin{equation*}
		\nu\left( \left( 1-\alpha \right) U_{i}\left( x\right) +\alpha V_{i}\left(
		x\right) \right) <\left( 1-\alpha \right) \nu\left( U_{i}\left( x\right)
		\right) +\alpha \nu\left( V_{i}\left( x\right) \right),
	\end{equation*}%
	since $\nu$ is strictly convex. Thus 
	\begin{equation*}
		\int\limits_{x\in E}\nu\left( \left( 1-\alpha
		\right) U_{i}\left( x\right) +\alpha V_{i}\left( x\right) \right)
		dx<\int\limits_{x\in E}\left( \left( 1-\alpha
		\right) \nu\left( U_{i}\left( x\right) \right) +\alpha \nu\left( V_{i}\left(
		x\right) \right) \right) dx,
	\end{equation*}%
	since $E$ has positive measure, and then we have $<$ instead of $\leq$ in (\ref{sc}) for $\alpha\in(0,1)$. This shows that $\widehat{F}^{\mathrm{en}}$ is strictly convex.

	Now, we prove that $\widehat{F}^{\mathrm{el}}$ is convex. We write $\widehat{F}^{\mathrm{el}}=KL$, where
	\begin{equation*}
		L\left( U\right) =\left\langle \widehat{\sigma }+\mathcal{T}\left( U\right) ,%
	\widehat{\sigma }+\mathcal{T}\left( U\right) \right\rangle,\ U\in A,
	\end{equation*}
	and prove that $L$ is convex (if $L$ is convex, then $\widehat{F}^{\mathrm{el}}=KL$ is convex).

	For $U,V\in A$ and $\alpha \in \left[ 0,1\right] $, we have%
	\begin{eqnarray*}
		&&L\left( \left( 1-\alpha \right) U+\alpha V\right)\\
		&&=\left\langle \widehat{%
			\sigma }+\mathcal{T}\left( \left( 1-\alpha \right) U+\alpha V\right) ,%
		\widehat{\sigma }+\mathcal{T}\left( \left( 1-\alpha \right) U+\alpha
		V\right) \right\rangle \\
		&&=\left\langle \widehat{\sigma }+\left( 1-\alpha \right) \mathcal{T}\left(
		U\right) +\alpha \mathcal{T}\left( V\right) ,\widehat{\sigma }+\left(
		1-\alpha \right) \mathcal{T}\left( U\right) +\alpha \mathcal{T}\left(
		V\right) \right\rangle \\
		&&=\left\langle \widehat{\sigma },\widehat{\sigma }\right\rangle +\left(
		1-\alpha \right) \left\langle \widehat{\sigma },\mathcal{T}\left( U\right)
		\right\rangle +\alpha \left\langle \widehat{\sigma },\mathcal{T}\left(
		V\right) \right\rangle \\
		&&\quad+\left( 1-\alpha \right) \left\langle \mathcal{T}\left( U\right) ,\widehat{%
			\sigma }\right\rangle +\left( 1-\alpha \right) ^{2}\left\langle \mathcal{T}%
		\left( U\right) ,\mathcal{T}\left( U\right) \right\rangle +\left( 1-\alpha
		\right) \alpha \left\langle \mathcal{T}\left( U\right) ,\mathcal{T}\left(
		V\right) \right\rangle \\
		&&\quad+\alpha \left\langle \mathcal{T}\left( V\right) ,\widehat{\sigma }%
		\right\rangle +\alpha \left( 1-\alpha \right) \left\langle \mathcal{T}\left(
		V\right) ,\mathcal{T}\left( U\right) \right\rangle +\alpha ^{2}\left\langle 
		\mathcal{T}\left( V\right) ,\mathcal{T}\left( V\right) \right\rangle \\
		&&=\left\langle \widehat{\sigma },\widehat{\sigma }\right\rangle+2\left( 1-\alpha \right) \left\langle \mathcal{T}\left( U\right) ,%
		\widehat{\sigma }\right\rangle +2\alpha \left\langle \mathcal{T}\left(
		V\right) ,\widehat{\sigma }\right\rangle \\
		&&\quad+\left(
		1-\alpha \right) ^{2}\left\langle \mathcal{T}\left( U\right) ,\mathcal{T}%
		\left( U\right) \right\rangle +2\left( 1-\alpha \right) \alpha
		\left\langle \mathcal{T}\left( U\right) ,\mathcal{T}\left( V\right)
		\right\rangle+\alpha ^{2}\left\langle \mathcal{T}\left(
		V\right) ,\mathcal{T}\left( V\right) \right\rangle 
	\end{eqnarray*}%
	and%
	\begin{eqnarray*}
		&&\left( 1-\alpha \right) L\left( U\right) +\alpha L\left( V\right)\\
		&&=\left(
		1-\alpha \right) \left\langle \widehat{\sigma }+\mathcal{T}\left( U\right) ,%
		\widehat{\sigma }+\mathcal{T}\left( U\right) \right\rangle +\alpha
		\left\langle \widehat{\sigma }+\mathcal{T}\left( V\right) ,\widehat{\sigma }+%
		\mathcal{T}\left( V\right) \right\rangle \\
		&&=\left( 1-\alpha \right) \left( \left\langle \widehat{\sigma },\widehat{%
			\sigma }\right\rangle +2\left\langle \mathcal{T}\left( U\right) ,\widehat{%
			\sigma }\right\rangle +\left\langle \mathcal{T}\left( U\right) ,\mathcal{T}%
		\left( U\right) \right\rangle \right) \\
		&&\quad+\alpha \left( \left\langle \widehat{\sigma },\widehat{\sigma }%
		\right\rangle +2\left\langle \mathcal{T}\left( V\right) ,\widehat{\sigma }%
		\right\rangle +\left\langle \mathcal{T}\left( V\right) ,\mathcal{T}\left(
		V\right) \right\rangle \right) \\
		&&=\left\langle \widehat{\sigma },\widehat{\sigma }\right\rangle
		+2\left( 1-\alpha \right) \left\langle \mathcal{T}\left( U\right) ,%
		\widehat{\sigma }\right\rangle +2\alpha \left\langle \mathcal{T}\left(
		V\right) ,\widehat{\sigma }\right\rangle \\
		&&\quad+(1-\alpha)\left\langle \mathcal{T}\left( U\right) ,\mathcal{T}\left( U\right)
		\right\rangle +\alpha\left\langle \mathcal{T}\left( V\right) ,\mathcal{T}\left(
		V\right) \right\rangle .
	\end{eqnarray*}%
	Thus%
	\begin{equation*}
		L\left( \left( 1-\alpha \right) U+\alpha V\right) \leq \left( 1-\alpha
		\right) L\left( U\right) +\alpha L\left( V\right)
	\end{equation*}%
	if and only if%
	\begin{eqnarray*}
		&&\left\langle \widehat{\sigma },\widehat{\sigma }\right\rangle+2\left( 1-\alpha \right) \left\langle \mathcal{T}\left( U\right) ,%
		\widehat{\sigma }\right\rangle +2\alpha \left\langle \mathcal{T}\left(
		V\right) ,\widehat{\sigma }\right\rangle \\
		&&\quad+\left(
		1-\alpha \right) ^{2}\left\langle \mathcal{T}\left( U\right) ,\mathcal{T}%
		\left( U\right) \right\rangle +2\left( 1-\alpha \right) \alpha
		\left\langle \mathcal{T}\left( U\right) ,\mathcal{T}\left( V\right)
		\right\rangle+\alpha ^{2}\left\langle \mathcal{T}\left(
		V\right) ,\mathcal{T}\left( V\right) \right\rangle \\
		&&\leq \left\langle \widehat{\sigma },\widehat{\sigma }\right\rangle
		+2\left( 1-\alpha \right) \left\langle \mathcal{T}\left( U\right) ,%
		\widehat{\sigma }\right\rangle +2\alpha \left\langle \mathcal{T}\left(
		V\right) ,\widehat{\sigma }\right\rangle\\
		&&\ \ \ +(1-\alpha)\left\langle \mathcal{T}\left( U\right) ,\mathcal{T}\left( U\right)
		\right\rangle +\alpha\left\langle \mathcal{T}\left( V\right) ,\mathcal{T}\left(
		V\right) \right\rangle 
	\end{eqnarray*}%
	if and only if%
	\begin{eqnarray*}
		&&0 \leq (1-\alpha)\alpha\left\langle \mathcal{T}\left( U\right) ,\mathcal{T}%
		\left( U\right) \right\rangle-2\left( 1-\alpha \right) \alpha \left\langle \mathcal{T}\left( U\right) ,%
		\mathcal{T}\left( V\right) \right\rangle\\
		&&\quad\quad +\left( 1-\alpha\right)\alpha \left\langle 
		\mathcal{T}\left( V\right) ,\mathcal{T}\left( V\right) \right\rangle
		\\
		&&\ \  =(1-\alpha)\alpha \left(\left\langle \mathcal{T}\left( U\right) ,\mathcal{T}%
		\left( U\right) \right\rangle-2\left\langle \mathcal{T}\left( U\right) ,%
		\mathcal{T}\left( V\right) \right\rangle  +\left\langle 
		\mathcal{T}\left( V\right) ,\mathcal{T}\left( V\right) \right\rangle\right)\\
		&&\ \  =(1-\alpha)\alpha \left\langle \mathcal{T}\left( U\right) -\mathcal{T}\left( V\right) ,\mathcal{T}\left( U\right) -\mathcal{T}\left( V\right) \right\rangle\\
		&&\ \ =(1-\alpha)\alpha \left\langle \mathcal{T}\left( U-V\right) ,\mathcal{T}\left( U-V\right)\right\rangle.
	\end{eqnarray*}
	Since
	$$
	(1-\alpha)\alpha \left\langle \mathcal{T}\left( U-V\right),\mathcal{T}\left( U-V\right)\right\rangle\geq 0
	$$
	holds, the convexity of $L$ is proved.  
\end{proof}

\subsection{Non-coercivity of $\widehat{F}$}

To apply Theorem \ref{existence2} in Appendix A for the existence of a minimum of $\widehat{F}$, one needs, besides the hypotheses on $A$ and $\widehat{F}$ stated in the previous subsection, the existence of a bounded minimizing sequence. This condition is satisfied if \(\widehat{F}\) is coercive. However, \(\widehat{F}\) may fail to be coercive.
\begin{pro}
In the case
\begin{equation}
\widehat{\sigma}\geq 0, N=2,\ \widehat{R}_{1}=\widehat{R}_{2}=r,\ c_{1}=c_{2}=1,\ Z_{1}=1\text{ and }Z_{2}=-1, \label{case}
\end{equation}
the functional $\widehat{F}$ is not coercive.
\end{pro}
\begin{proof}
In the case (\ref{case}), we have, for $U\in X$,
$$
\mathcal{T}(U)=\mathcal{L}U_{1}-\mathcal{L}%
U_{2}=\mathcal{L}(U_1-U_2),
$$
where
\begin{equation*}
	[\mathcal{L}u]\left( x\right)=\frac{1}{2r}\int\limits_{r}^{\min\{x+r,1-r\}}u\left( s\right)\min\{x+r-s,2r\}ds,\ x\in \left[ 0,1\right] \ \text{%
		and\ }u\in L^{p}\left( r,1-r\right),
\end{equation*}
(for the case $r=0$, see (\ref{R=0})). Then, for $U\in X_+$, 
\begin{equation*}
	\widehat{F}\left( U\right)=\int\limits_{r}^{1-r}\nu\left( U_{1}\left( x\right) \right)
	dx+\int\limits_{r}^{1-r}\nu\left( U_{2}\left( x\right) \right) dx+K\left\langle \widehat{\sigma }+\mathcal{L}\left( U_{1}-U_{2}\right) ,%
	\widehat{\sigma }+\mathcal{L}\left( U_{1}-U_{2}\right) \right\rangle.
\end{equation*}%
Moreover, we have
\begin{equation*}
E(U)=2\widehat{\sigma}+\left[\mathcal{T}U\right](1)=2\widehat{\sigma}+\left[\mathcal{L}\left(U_{1}-U_{2}\right) \right](1).
\end{equation*}%

Let $\alpha,\beta>0$ be such that  $1<\frac{\beta}{\alpha}<p$. Consider, for any $n\in\mathbb{N}^+$ sufficiently large, a subinterval $I_n$ of $\left[ r,1-r\right]$ of measure $\frac{1}{n^{\beta}}$ and $U^{(n)}\in A$ with components
\begin{equation*}
U_{1}^{\left( n\right) }\left( x\right) =\left\{ 
\begin{array}{l}
	n^{\alpha}\text{ if }x\in I_n\\ 
	0\ \text{otherwise}%
\end{array}%
\right. ,\ x\in \left[ r,1-r\right],
\end{equation*}
and
\begin{equation*}
U_{2}^{\left( n\right) }\left( x\right) =\gamma +U_{1}^{\left( n\right)
}\left( x\right) ,\ x\in \left[ r,1-r\right],
\end{equation*}%
where $\gamma\geq 0$ is such that the electroneutrality condition
\begin{equation*}
E(U^{(n)})=2\widehat{\sigma }+\left[ \mathcal{L}\left( U_{1}^{(n)}-U_{2}^{(n)}\right) \right]
\left( 1\right) =2\widehat{\sigma }+\left[ \mathcal{L}\left( -\gamma \right) %
\right] \left( 1\right) =2\widehat{\sigma }-\gamma(1-2r)=0
\end{equation*}%
is satisfied,  i.e., $\gamma =\frac{2\widehat{\sigma }}{1-2r}$.

We have 
\begin{equation*}
\left\Vert U_{1}^{\left( n\right) }\right\Vert
_{p}^{p}=\int\limits_{r}^{1-r}\left( U_{1}^{\left( n\right) }\left( x\right)
\right) ^{p}dx=\int\limits_{x\in I_n}n^{\alpha p}dx=n^{\alpha p}\cdot \frac{1%
}{n^{\beta}}=n^{\alpha p-\beta}\rightarrow \infty ,\ n\rightarrow \infty,
\end{equation*}%
and%
\begin{equation*}
\left\Vert U_{1}^{\left( n\right) }\right\Vert
_{p}^{p}=\int\limits_{r}^{1-r}\left( U_{1}^{\left( n\right) }\left( x\right)
\right) ^{p}dx\leq \int\limits_{r}^{1-r}\left( \gamma +U_{1}^{\left( n\right)
}\left( x\right) \right) ^{p}dx=\left\Vert U_{2}^{\left( n\right)
}\right\Vert _{p}^p
\end{equation*}%
and then%
\begin{equation*}
\left\Vert U_{2}^{\left( n\right) }\right\Vert _{p}^{p}\rightarrow \infty ,\
n\rightarrow \infty,
\end{equation*}%
as well. Moreover, we have 
\begin{eqnarray*}
\widehat{F}\left( U^{\left( n\right) }\right)&=&\int\limits_{r}^{1-r}\nu\left( U_{1}^{(n)}\left(
x\right) \right) dx+\int\limits_{r}^{1-r}\nu\left( U_{2}^{(n)}\left( x\right)
\right) dx\\
&&+K\left\langle \widehat{\sigma }+\mathcal{L}\left( U_{1}^{(n)}-U_{2}^{(n)}\right) ,\widehat{\sigma }+\mathcal{L}\left(
U_{1}^{(n)}-U_{2}^{(n)}\right) \right\rangle
\end{eqnarray*}%
with 
\begin{eqnarray*}
\int\limits_{r}^{1-r}\nu\left( U_{1}^{(n)}\left( x\right) \right) dx
&=&\int\limits_{r}^{1-r}\left( U_{1}^{(n)}\left( x\right) \left( \log U_{1}^{(n)}\left(
x\right) -1\right) +1\right) dx \\
&=&\int\limits_{r}^{1-r}U_{1}^{(n)}\left( x\right) \left( \log U_{1}^{(n)}\left(
x\right) -1\right) dx+1-2r \\
&=&\int\limits_{x\in I_n}n^{\alpha}\left( \log n^{\alpha}-1\right) dx+1-2r
\\
&=&n^{\alpha}\left( \log n^{\alpha}-1\right) \cdot \frac{1}{n^{\beta}}+1-2r \\
&=&\frac{\alpha\log n-1}{n^{\beta-\alpha}}+1-2r
\end{eqnarray*}%
and%
\begin{eqnarray*}
\int\limits_{r}^{1-r}\nu\left( U_{2}^{(n)}\left( x\right) \right) dx
&=&\int\limits_{r}^{1-r}\left( U_{2}^{(n)}\left( x\right) \left( \log U_{2}^{(n)}\left(
x\right) -1\right) +1\right) dx \\
&=&\int\limits_{r}^{1-r}U_{2}^{(n)}\left( x\right) \left( \log U_{2}^{(n)}\left(
x\right) -1\right) dx+1-2r \\
&=&\int\limits_{x\in I_n}\left( \gamma +n^{\alpha}\right) \left( \log
\left( \gamma +n^{\alpha}\right) -1\right) dx\\
&&\quad +\int\limits_{x\in [r,1-r]\setminus I_n}\gamma \left( \log
\gamma -1\right) dx+1-2r \\
&=&\left( \gamma +n^{\alpha}\right) \left( \log \left( \gamma +n^{\alpha}\right)
-1\right) \cdot \frac{1}{n^{\beta}}\\
&&\quad +\gamma(\log \gamma-1)(1-2r-\frac{1}{n^\beta})+1-2r \\
&=&\left( \frac{\gamma }{n^{\alpha}}+1\right) \left( \alpha\log n+\log \left( \frac{%
	\gamma }{n^{\alpha}}+1\right) -1\right) \cdot \frac{1}{n^{\beta-\alpha}}\\
	&&\quad +(1-2r)(\gamma(\log\gamma-1)+1)-\frac{\gamma(\log\gamma-1)}{n^{\beta}}
	\end{eqnarray*}
and
\begin{eqnarray*}
	K\left\langle \widehat{\sigma }+\mathcal{L}\left( U_{1}^{(n)}-U_{2}^{(n)}\right) ,\widehat{\sigma }+\mathcal{L}\left(
	U_{1}^{(n)}-U_{2}^{(n)}\right) \right\rangle&=&K\left\langle \widehat{\sigma }+\mathcal{L}\left( -\gamma \right) ,%
\widehat{\sigma }+\mathcal{L}\left( -\gamma \right) \right\rangle\\
&=&K\left\langle \widehat{\sigma }-\gamma\mathcal{L}(1),%
\widehat{\sigma }-\gamma\mathcal{L}(1) \right\rangle.
\end{eqnarray*}
Therefore, as $n\rightarrow \infty$,  
	\begin{equation*}
\left\Vert U^{\left( n\right) }\right\Vert\rightarrow \infty
\end{equation*}
and
\begin{equation*}
	\widehat{F}\left( U^{\left( n\right) }\right) \rightarrow \left( 1-2r\right)(\gamma(\log\gamma-1)+2)+K\left\langle \widehat{\sigma }-\gamma\mathcal{L}(1),%
\widehat{\sigma }-\gamma\mathcal{L}(1) \right\rangle.
\end{equation*}
This shows that $\widehat{F}$ is not coercive.

\end{proof}

The lack of coercivity of $\widehat{F}$ forces us to apply Theorem \ref{existence2} to the restriction of $\widehat{F}$ to a bounded subset of the domain $A$, as will be illustrated in the next subsection.

\subsection {The problem MP$(r)$}

Since the functional $\widehat{F}$ is, in general, non-coercive, we cannot ensure the existence of a bounded minimizing sequence for MP. Consequently, Theorem \ref{existence2} in Appendix A cannot be invoked to conclude that MP admits a minimum, despite all other hypotheses being satisfied.

For this reason, instead of MP, we consider the following problem MP$(r)$, where $r>0$.

\begin{itemize}
	\item[MP$(r)$] \textit{minimize $\widehat{F}\left( U\right)$ with respect to $U\in X$, subject to the \emph{non-negativity condition}%
		$$
		U_{i}\geq 0,\ i\in\{1,\ldots ,N\},
		$$
		to the \emph{electroneutrality condition}
		$$
		E(U)=0
		$$
		and to the \emph{boundedness condition}}
		$$
		\Vert U\Vert\leq r.
		$$
\end{itemize}

For the problem MP$(r)$, any minimizing sequence is trivially bounded, and then we can apply Theorems \ref{existence2} and \ref{uniqueness} in Appendix A to conclude that MP$(r)$ has a unique minimum. Observe that, for the problem MP$(r)$, the functional $\widehat{F}$ is considered to have the domain
$$
A(r)=\{U\in X: U_i\geq 0,\ i\in\{1,\ldots,N\},E(U)=0\text{ and }\Vert U\Vert\leq r\},
$$
which is a subset of the domain (\ref{A}) for the problem MP. The domain $A(r)$ is non-empty for $r$ sufficiently large since $A\neq\emptyset$, and it is closed and convex since it is the intersection between the closed and convex set $A$ and the closed ball in the space $X$ centered at $0$ with radius $r$, which is convex. Moreover, the functional $\widehat{F}$ with domain $A(r)$ is continuous and strictly convex, since it is continuous and strictly convex on the larger domain (\ref{A}).

In the following, we denote the minimum of the problem MP$(r)$ by $\overline{U}(r)$. Moreover, let $r_0>0$ be such that $A(r)\neq \emptyset$, for $r\geq r_0$.

\subsection{From MP$(r)$ to MP}

Although the problem MP$(r)$ is not MP and the minimum of MP$(r)$ may not be the minimum of MP, the following result relates the two problems. Observe that MP has at most one minimum, since $\widehat{F}$ is strictly convex.

\begin{pro} \label{twoalt}
	Let $(r_n)_{n\in\mathbb{N}^+}$ be an increasing divergent sequence of positive numbers such that $r_n\geq r_0$, $n\in\mathbb{N}^+$. One and only one of the following two alternatives holds:
	\begin{itemize}
		\item [a)] there exists $n_0\in\mathbb{N}^+$ such that, for $n\geq n_0$, $\overline{U}(r_n)=\overline{U}(r_{n_0})$;
		\item [b)] there exists a subsequence $(r_{n_k})_{k\in\mathbb{N}^+}$ such that $\left\Vert \overline{U}(r_{n_k})\right\Vert\rightarrow +\infty$, $k\rightarrow \infty$.
	\end{itemize}
	If a) holds, then $\overline{U}(r_{n_0})$ is the unique minimum of the problem MP.  
\end{pro}

\begin{proof}
	Suppose that a) does not hold, i.e., for any $n_0\in\mathbb{N}^+$ there exists $n>n_0$ such that $\overline{U}(r_n)\neq \overline{U}(r_{n_0})$. Therefore, we can obtain a subsequence $(r_{n_k})_{k\in\mathbb{N}^+}$ of $(r_n)_{n\in\mathbb{N}^+}$ such that  $\overline{U}(r_{n_{k+1}})\neq \overline{U}(r_{n_k})$ for all $k\in\mathbb{N}^+$. For any $k\in\mathbb{N}^+$, we have $\Vert \overline{U}(r_{n_{k+1}})\Vert>r_{n_k}$; otherwise, the problem MP$(r_{n_k})$ would have
	two minima instead of only one. In fact, since $r_{n_{k+1}}>r_{n_k}$, and thus $A(r_{n_k})\subseteq A(r_{n_{k+1}})$, it turns out that $\widehat{F}(\overline{U}(r_{n_{{k+1}}}))\leq \widehat{F}(U)$ for all $U\in A(r_{n_k})$. We conclude that $\left\Vert \overline{U}(r_{n_k})\right\Vert\rightarrow +\infty$.
	
	If a) holds, then we have $\widehat{F}(\overline{U}(r_{n_0}))\leq \widehat{F}(U)$ for all $U\in A(r_n)$ and for all $n\geq n_0$, i.e.,  $\widehat{F}(\overline{U}(r_{n_0}))\leq \widehat{F}(U)$ for all $U\in A=\bigcup\limits_{n=n_0}^\infty A(r_n)$.
\end{proof}

\subsection{The three lemmas}
Now, we show that the alternative b) in Proposition \ref{twoalt} cannot hold. To this aim, we introduce three lemmas.

Although $\nu^\prime$ coincides with $\log$, we will continue to write $\nu'$ below, since the next three lemmas remain valid under the assumption that $\nu$ is a non-negative continuous function defined on $[0,+\infty)$ and $\nu'$ is an arbitrary strictly increasing continuous function defined on $(0,+\infty)$ with range $\nu^{\prime}(0,+\infty)=\mathbb{R}$. Observe that $\nu^\prime(0)=-\infty$ and $\nu^\prime(+\infty)=+\infty$, and then $\left(\nu^\prime\right)^{-1}(-\infty)=0$ and $\left(\nu^\prime\right)^{-1}(+\infty)=+\infty$.

Before the first lemma, we introduce the notion of essential oscillation of a function. Let $u$ be a measurable function defined on a domain $\mathcal{D}$. The \emph{essential oscillation} of $u$ is defined as
$$
\operatorname*{ess\,osc} (u):=\operatorname*{ess\,sup}_{x\in \mathcal{D}} u(x)-\operatorname*{ess\,inf}_{x\in \mathcal{D}} u(x).
$$

Let $U\in A$. We define the $\nu^\prime$-\emph{essential oscillation} of $U$ as
\begin{equation}
\nu^\prime\operatorname*{-ess\,osc} (U)=\max\limits_{i\in\{1,\ldots,N\}}\operatorname*{ess\,osc}  (\nu^\prime\circ U_i), \label{nuprimeessosc}
\end{equation}
where, for $i\in\{1,\ldots,N\}$,
\begin{eqnarray}
\operatorname*{ess\,osc} (\nu^\prime \circ U_i)&=&\operatorname*{ess\,sup}_{x\in [\widehat{R}_i,1-\widehat{R}_i]} \nu^\prime(U_i(x))-\operatorname*{ess\,inf}_{x\in [\widehat{R}_i,1-\widehat{R}_i]} \nu^\prime(U_i(x))\notag\\
&=&\nu^\prime\left(\operatorname*{ess\,sup}_{x\in [\widehat{R}_i,1-\widehat{R}_i]} U_i(x)\right)-\nu^\prime\left(\operatorname*{ess\,inf}_{x\in [\widehat{R}_i,1-\widehat{R}_i]} U_i(x)\right). \notag\\
\label{diffnuprime}
\end{eqnarray}
The second equality holds since $\nu^\prime$ is strictly increasing and continuous. We set $\operatorname*{ess\,osc} (\nu^\prime \circ U_i)=0$ if $U_i=0$. 

Here is the first lemma.
\begin{lem} \label{Lemma1}
Let $\omega\geq 0$. Let $U\in A$ be such that $\widehat{F}(U)\leq \omega$.
Let
\begin{equation}
p(\omega):=4\max\limits_{i\in\{1,\ldots,N\}}\vert Z_i\vert \sqrt{K\omega}. \label{c}
\end{equation}
If
\begin{equation}
\nu^\prime\operatorname*{-ess\,osc}(U)> p(\omega), \label{cL1}
\end{equation}
then there exists $V\in A$ such that $\Vert V\Vert\leq \Vert U\Vert$ and $\widehat{F}(V)<\widehat{F}(U)$.
\end{lem}

\begin{proof}
	Suppose $\nu^\prime\operatorname*{-ess\,osc} (U)>p(\omega)$. Let $i\in\{1,\ldots,N\}$ be such that
	$$
	\nu^\prime\operatorname*{-ess\,osc} (U)=\operatorname*{ess\,osc} (\nu^\prime \circ U_i).
	$$
	There exist $a,b\geq 0$ with $a<b$ and $\nu^\prime(b)-\nu^\prime(a)$ arbitrarily close (from below) to $\operatorname*{ess\,osc} (\nu^\prime \circ U_i)$ such that the disjoint sets
	$$
	\mathcal{A}=\left\{x\in \left[\widehat{R}_i,1-\widehat{R}_i\right]: U_i(x)\leq a\right\}\text{\ and \ }\mathcal{B}=\left\{x\in \left[\widehat{R}_i,1-\widehat{R}_i\right]: U_i(x)\geq b\right\}
	$$
	have non-zero measures $m(\mathcal{A})$ and $m(\mathcal{B})$.

	Let $\Delta\in X$ be a perturbation of $U$ with components
	\begin{eqnarray}
		&&\Delta_i=\alpha\chi_{\mathcal{A}}-\beta\chi_{\mathcal{B}}\notag \\
		&&\Delta_k=0\text{\ for $k\in\{1,\ldots,N\}$ with $k\neq i$}\notag ,\\
		\label{Delta}
	\end{eqnarray}
	where $\alpha,\beta>0$ and $\chi_{\mathcal{A}}$ and $\chi_{\mathcal{B}}$ are the characteristic functions of $\mathcal{A}$ and $\mathcal{B}$.
	
	Let $V=U+\Delta$. We check the following conditions:
	\begin{itemize}
		\item [1)]$V_k\geq 0$, $k\in\{1,\ldots,N\}$, i.e., $V$ satisfies the non-negativity condition;
		\item [2)] $E(V)=0$, i.e., $V$ satisfies the electroneutrality condition;
		\item [3)] $\Vert V\Vert \leq \Vert U\Vert$;
		\item [4)] $\widehat{F}(V)<\widehat{F}(U)$.
	\end{itemize}
	
	Condition 1). We have $V_k\geq 0$, $k\in\{1,\ldots,N\}$, if $\beta\leq b$.
	
	Condition 2). We have $E(V)=0$ if and only if $[\mathcal{T}(\Delta)](1)=0$ and we have (recall (\ref{Tu1}))
	\begin{eqnarray*}
		[\mathcal{T}(\Delta)](1)=\alpha Z_ic_im(\mathcal{A})-\beta Z_ic_im(\mathcal{B})=0
	\end{eqnarray*} 
	if and only if
	$$
	\beta=k\alpha\text{\ with\ }k=\frac{m(\mathcal{A})}{m(\mathcal{B})}.
	$$
	
	For the next conditions, 3) and 4), we assume
	$$
	\alpha\leq\alpha_0:=\frac{b}{k}\text{\ and\ }\beta=k\alpha,
	$$
	to have conditions 1) and 2) satisfied.

	Condition 3). We have $\Vert V\Vert\leq \Vert U\Vert$ if  $\Vert V_i\Vert_p\leq \Vert U_i\Vert_p$. Observe that
	\begin{eqnarray*}
	\Vert V_i\Vert_p^p-\Vert U_i\Vert_p^p&=&\int\limits_{x\in \mathcal{A}}(U_i(x)+\alpha)^pdx+\int\limits_{x\in \mathcal{B}}(U_i(x)-k\alpha)^pdx\\
	&&-\left(\int\limits_{x\in \mathcal{A}}U_i(x)^pdx+\int\limits_{x\in \mathcal{B}}U_i(x)^pdx\right).
	\end{eqnarray*}
	Therefore, there exists $\alpha_1\in\left(0,\alpha_0\right]$ such that, for $\alpha\in(0,\alpha_1]$, $\Vert V_i\Vert_p<\Vert U_i\Vert_p$ if the function
	$$
	g(\alpha)=\int\limits_{x\in \mathcal{A}}(U_i(x)+\alpha)^pdx+\int\limits_{x\in \mathcal{B}}(U_i(x)-k\alpha)^pdx,\ \alpha\in\left[0,\alpha_0\right],
	$$
	has $g^\prime(0)<0$. Since the function $u\mapsto u^{p-1}$, $u\geq 0$, is increasing,  we obtain
	\begin{eqnarray*}
		g^\prime(0)&=&g^\prime(\alpha)\big |_{\alpha=0}\\
		&=&p\left(\int\limits_{x\in \mathcal{A}}(U_i(x)+\alpha)^{p-1}dx-k\int\limits_{x\in \mathcal{B}}(U_i(x)-k\alpha)^{p-1}dx\right)\Bigg |_{\alpha=0}\\
		&=&p\left(\int\limits_{x\in \mathcal{A}}(U_i(x))^{p-1}dx-k\int\limits_{x\in \mathcal{B}}(U_i(x))^{p-1}dx\right)\\
		&\leq&p\left(a^{p-1}m(\mathcal{A})-kb^{p-1}m(\mathcal{B})\right)\\
		&=&pm(\mathcal{A})\cdot \left(a^{p-1}-b^{p-1}\right)<0.
	\end{eqnarray*}
	We conclude that there exists $\alpha_1\in(0,\alpha_0]$ such that, for $\alpha\in(0,\alpha_1]$, we have $\Vert V\Vert \leq \Vert U \Vert$. 
	
	Condition 4). We write $\widehat{F}=\widehat{F}^{\mathrm{en}}+\widehat{F}^{\mathrm{el}}$ (recall (\ref{FhatU}), (\ref{FenhatU}) and (\ref{FelhatU})). Regarding the entropic part $\widehat{F}^{\mathrm{en}}$, we have
	\begin{eqnarray*}
		\widehat{F}^{\mathrm{en}}(V)-\widehat{F}^{\mathrm{en}}(U)&=&c_i\int\limits_{x\in \mathcal{A}}\nu(U_i(x)+\alpha)dx+c_i\int\limits_{x\in \mathcal{B}}\nu(U_i(x)-k\alpha)dx\\
		&&-\left(c_i\int\limits_{x\in \mathcal{A}}\nu(U_i(x))dx+c_i\int\limits_{x\in \mathcal{B}}\nu(U_i(x))dx\right).
	\end{eqnarray*}
	Regarding the electrostatic part $\widehat{F}^{\mathrm{el}}$, written as $\widehat{F}^{\mathrm{el}}=KL$ as
	in the proof of Proposition \ref{Thconvex}, we have
	\begin{eqnarray*}
		L(V)-L(U)&=&\langle\widehat{\sigma}+\mathcal{T}(V),\widehat{\sigma}+\mathcal{T}(V)\rangle
		-\langle\widehat{\sigma}+\mathcal{T}(U),\widehat{\sigma}+\mathcal{T}(U)\rangle\\
		&=&\langle\widehat{\sigma}+\mathcal{T}(U)+\mathcal{T}(\Delta),\widehat{\sigma}+\mathcal{T}(U)+\mathcal{T}(\Delta)\rangle
		-\langle\widehat{\sigma}+\mathcal{T}(U),\widehat{\sigma}+\mathcal{T}(U)\rangle\\
		&=&2\langle\widehat{\sigma}+\mathcal{T}(U),\mathcal{T}(\Delta)\rangle+\langle \mathcal{T}(\Delta),\mathcal{T}(\Delta)\rangle\\
		&=&2\langle\widehat{\sigma}+\mathcal{T}(U),Z_ic_i\mathcal{L}_i(\alpha\chi_{\mathcal{A}}-k\alpha \chi_{\mathcal{B}})\rangle\\
		&&+\langle Z_ic_i\mathcal{L}_i(\alpha\chi_\mathcal{A}-k\alpha \chi_\mathcal{B}),Z_ic_i\mathcal{L}_i(\alpha\chi_\mathcal{A}-k\alpha \chi_\mathcal{B})\rangle\\
		&=&2Z_ic_i\alpha\langle\widehat{\sigma}+\mathcal{T}(U),\mathcal{L}_i(\chi_\mathcal{A}-k\chi_\mathcal{B})\rangle\\
		&&+Z_i^2c_i^2\alpha^2\langle \mathcal{L}_i(\chi_\mathcal{A}-k \chi_\mathcal{B}),\mathcal{L}_i(\chi_\mathcal{A}-k\chi_\mathcal{B})\rangle.
	\end{eqnarray*}

	Therefore, there exists $\alpha_2\in\left(0,\alpha_1\right]$ such that, for $\alpha\in(0,\alpha_2]$,
	$$
	\widehat{F}(V)-\widehat{F}(U)=\widehat{F}^{\mathrm{en}}(V)-\widehat{F}^{\mathrm{en}}(U)+K(L(V)-L(U))<0
	$$
	if the function
	\begin{eqnarray*}
		&&h(\alpha)=c_i\int\limits_{x\in \mathcal{A}}\nu(U_i(x)+\alpha)dx+c_i\int\limits_{x\in \mathcal{B}}\nu(U_i(x)-k\alpha)dx\\
		&&\quad\quad\quad\   +2KZ_ic_i\langle\widehat{\sigma}+\mathcal{T}(U),\mathcal{L}_i(\chi_\mathcal{A}-k\chi_\mathcal{B})\rangle\alpha\\
		&&\quad\quad\quad\   +KZ_i^2c_i^2\langle \mathcal{L}_i(\chi_\mathcal{A}-k \chi_\mathcal{B}),\mathcal{L}_i(\chi_\mathcal{A}-k\chi_\mathcal{B})\rangle\alpha^2,\\
		&&\quad\quad\quad\  \alpha\in\left[0,\alpha_1\right],
	\end{eqnarray*}
	has $h^\prime(0)<0$. We have
	\begin{eqnarray}
		h^\prime(0)&=&h^\prime(\alpha)\big |_{\alpha=0}\notag \\
		&=&\left(c_i\int\limits_{x\in \mathcal{A}}\nu^\prime(U_i(x)+\alpha)dx-c_ik\int\limits_{x\in \mathcal{B}}\nu^\prime(U_i(x)-k\alpha)dx\right.\notag \\
		&&\quad\quad\left. +2KZ_ic_i\langle\widehat{\sigma}+\mathcal{T}(U),\mathcal{L}_i(\chi_\mathcal{A}-k\chi_\mathcal{B})\rangle\right.\notag\\
		&&\quad\quad  +2KZ_i^2c_i^2
		\langle \mathcal{L}_i(\chi_\mathcal{A}-k\chi_\mathcal{B}),\mathcal{L}_i(\chi_\mathcal{A}-k\chi_\mathcal{B})\rangle\alpha\Bigg)\Bigg|_{\alpha=0}\notag\\
		&=& c_i\int\limits_{x\in \mathcal{A}}\nu^\prime(U_i(x))dx-c_ik\int\limits_{x\in \mathcal{B}}\nu^\prime(U_i(x))dx\notag \\
		&&+2KZ_ic_i\langle\widehat{\sigma}+\mathcal{T}(U),\mathcal{L}_i(\chi_\mathcal{A}-k\chi_\mathcal{B})\rangle. \notag \\ \label{dersum}
	\end{eqnarray}
	The fact that we can differentiate under the integral signs in the second
	equality of \eqref{dersum} follows from the monotone convergence theorem applied
	to the difference quotients, which decrease monotonically to the derivative
	as $\alpha\to0^+$, since $\nu$ is convex. Regarding the first term in (\ref{dersum}), we have
	\begin{eqnarray*}
		&&c_i\int\limits_{x\in \mathcal{A}}\nu^\prime(U_i(x))dx-c_ik\int\limits_{x\in \mathcal{B}}\nu^\prime(U_i(x))dx\\
		&&\leq c_i\nu^\prime(a)m(\mathcal{A})-c_ik\nu^\prime(b)m(\mathcal{B})=c_i\left(\nu^\prime(a)-\nu^\prime(b)\right)m(\mathcal{A}).
	\end{eqnarray*}
	Regarding the second term in (\ref{dersum}), we have
	\begin{eqnarray*}
	&&\vert \langle\widehat{\sigma}+\mathcal{T}(U),\mathcal{L}_i(\chi_\mathcal{A}-k\chi_\mathcal{B})\rangle \vert \\
	&&\leq \Vert \widehat{\sigma}+\mathcal{T}(U)\Vert_2\Vert \mathcal{L}_i(\chi_\mathcal{A}-k\chi_\mathcal{B})\Vert_2
	\leq \Vert \widehat{\sigma}+\mathcal{T}(U)\Vert_2\Vert \mathcal{L}_i(\chi_\mathcal{A}-k\chi_\mathcal{B})\Vert_\infty.
	\end{eqnarray*}
	Since
	$$
	K\Vert \widehat{\sigma}+\mathcal{T}(U)\Vert_2^2\leq \sum\limits_{i=1}^{N}c_{i}\left\langle \nu\circ U
	_i\right\rangle_i+K\Vert \widehat{\sigma}+\mathcal{T}(U)\Vert_2^2=\widehat{F}(U)\leq \omega,
	$$
	we obtain
	$$
	\Vert \widehat{\sigma}+\mathcal{T}(U)\Vert_2\leq \sqrt{\frac{\omega}{K}}.
	$$
	Moreover, by (\ref{operatornorm}) we obtain
	\begin{equation*}
	\Vert \mathcal{L}_i(\chi_\mathcal{A}-k\chi_\mathcal{B})\Vert_\infty\leq \Vert \chi_\mathcal{A}-k\chi_\mathcal{B}\Vert_1\leq \Vert \chi_\mathcal{A}\Vert_1+k\Vert \chi_\mathcal{B}\Vert_1=m(\mathcal{A})+km(\mathcal{B})=2m(\mathcal{A}).
	\end{equation*}
	Thus, in (\ref{dersum}) we have
	$$
	\vert 2KZ_ic_i\langle\widehat{\sigma}+\mathcal{T}(U),\mathcal{L}_i(\chi_\mathcal{A}-k\chi_\mathcal{B})\rangle\vert \leq 4\vert Z_i\vert c_i\sqrt{K\omega}\ m(\mathcal{A}).
	$$
Therefore, we obtain
	\begin{equation*}
		h^\prime(0)\leq -c_i\left(\nu^\prime(b)-\nu^\prime(a)-4\vert Z_i\vert \sqrt{K\omega}\right)m(\mathcal{A}).
	\end{equation*}
	We can conclude that if
	$$
	\nu^\prime\operatorname*{-ess\,osc} (U)=\operatorname*{ess\,osc} (\nu^\prime \circ U_i)> p(\omega),
	$$
	where $p(\omega)$ is defined in (\ref{c}), then we can take  $a$ and $b$ with $\nu^\prime(b)-\nu^\prime(a)$ arbitrarily close to $\operatorname*{ess\,osc} (\nu^\prime\circ U_i)$ so as to have $h^\prime(0)<0$. As a consequence, there exists $\alpha_2\in\left(0,\alpha_1\right]$ such that, for $\alpha\in(0,\alpha_2]$, $\widehat{F}(V)<\widehat{F}(U)$.
	
	In summary, if  (\ref{cL1}) holds, then there exist  $a,b\geq 0$ with $a<b$ such that, by taking $\alpha>0$ sufficiently small and $$\beta=\frac{m(\mathcal{A})}{m(\mathcal{B})}\alpha$$ in (\ref{Delta}), we obtain $V\in A$, $\Vert V\Vert \leq \Vert U\Vert$ and $\widehat{F}(V)<\widehat{F}(U)$ with $V=U+\Delta$.

\end{proof}

Before the second lemma, we introduce the sets $I^+$ and $I^{-}$ and the notion of $\varepsilon$-small co-ions.

In the case $\widehat{\sigma}\neq 0$ (charged planes), let
$$
I^+:=\{i\in\{1,\ldots,N\}: Z_i\widehat{\sigma}>0\}\text{\ and\ }I^-:=\{i\in\{1,\ldots,N\}: Z_i\widehat{\sigma}<0\},
$$
i.e., $I^+$  is the set of indices $i$ such that the $i$-th type of ion is a co-ion and $I^-$  is the set of indices $i$ such that the $i$-th type of ion is a counter-ion. Observe that $I^{+}\cap I^{-}=\emptyset$ and $I^{+}\cup I^{-}=\{1,\ldots,N\}$. Moreover, we have $I^+\neq \emptyset$ and $I^{-}\neq \emptyset$, since (\ref{bulken}) holds. Finally, observe that the electroneutrality condition (\ref{encond}) can be written as
\begin{equation}
	\sum\limits_{i\in I^-}\vert Z_i\vert c_i \Vert U_i\Vert_1=2\vert \widehat{\sigma}\vert+\sum\limits_{i\in I^+}\vert Z_i\vert c_i \Vert U_i\Vert_1. \label{encond2}
\end{equation}
In the limiting case $\widehat{\sigma}=0$ (uncharged planes), we set
$$
I^+:=\{i\in\{1,\ldots,N\}: Z_i>0\}\text{\ and\ }I^-:=\{i\in\{1,\ldots,N\}: Z_i<0\}.
$$
With this definition, $I^{+}$ and $I^{-}$ remain a partition of $\{1,\ldots,N\}$ and (\ref{encond2}) remains valid for $\widehat{\sigma}=0$.

Let $\varepsilon>0$. We say that $U\in A$ \emph{has $\varepsilon$-small co-ions} if $\Vert U_i\Vert_1\leq \varepsilon$ for all $i\in I^+$.

Finally, for $U\in A$, let
$$
\Vert U\Vert_\infty=\max\limits_{i=1,\ldots,N}\Vert U_i\Vert_\infty.
$$
Since $U_i\in L^p(\widehat{R}_i,1-\widehat{R}_i)$, $i\in\{1,\ldots,N\}$, we have
$$
\Vert U_i\Vert_\infty=\operatorname*{ess\,sup}_{x\in [\widehat{R}_i,1-\widehat{R}_i]} U_i(x)\in[0,+\infty].
$$ 

Here is the second lemma. 
\begin{lem} \label{Lemma2}
	Let $\omega\geq 0$ and let $\varepsilon>0$. Let $U\in A$ be such that $\widehat{F}(U)\leq \omega$ and $U$ does not
	have $\varepsilon$-small co-ions. Let
	\begin{equation}
	q(\varepsilon):=\min\left\{\frac{2\vert \widehat{\sigma}\vert+\min\limits_{k\in I^+}\vert Z_k\vert c_k\varepsilon}{\max\limits_{k\in I^-}\vert Z_k\vert c_k\vert I^{-}\vert},\varepsilon \right\} \label{q}
	\end{equation}
	and let
	\begin{equation}
		r(\omega,\varepsilon):=\max\limits_{\substack{i,j\in\{1,\ldots,N\}\\\text{$Z_i$ and $Z_j$ are of opposite sign}}}\left(-\frac{|Z_i|}{|Z_j|}\nu^\prime\left(q(\varepsilon)\right)+4\vert Z_i\vert \sqrt{K\omega}\right). \label{M}
	\end{equation}
	If
	\begin{equation}
\Vert U\Vert_\infty >\left(\nu^\prime\right)^{-1}\left(r(\omega,\varepsilon)\right), \label{cL2}
	\end{equation}
then there exists $V\in A$ such that $\Vert V\Vert \leq \Vert U\Vert$ and  $\widehat{F}(V)<\widehat{F}(U)$.
\end{lem}

\begin{proof}
	Suppose $\Vert U\Vert_{\infty}>\left(\nu^\prime\right)^{-1}\left(r(\omega,\varepsilon)\right)$. Let $i\in\{1,\ldots,N\}$ be such that
	$
	\Vert U\Vert_\infty =\Vert U_i\Vert_\infty.
	$
	Let $M\in\mathbb{R}$ be such that
	$$
	\Vert U\Vert_{\infty}>M>\left(\nu^\prime\right)^{-1}\left(r(\omega,\varepsilon)\right).
	$$
	There exists a subset $\mathcal{A}$ of $[\widehat{R}_i,1-\widehat{R}_i]$ with $m(\mathcal{A})>0$ where we have
	$
	U_i\geq M>0.
	$
	
	If $i\in I^+$, then, by the electroneutrality condition (\ref{encond2}), there exists $j\in I^{-}$ such that
	$$
	\vert Z_j\vert c_j \Vert U_j\Vert_1\geq \frac{2\vert \widehat{\sigma}\vert+\sum\limits_{k\in I^+}\vert Z_k\vert c_k \Vert U_k\Vert_1}{\vert I^{-}\vert} 
	$$
	and, by the fact that $U$ does not have $\varepsilon$-small co-ions, there exists $j_1\in I^+$ such that $\Vert U_{j_1}\Vert_1>\varepsilon$ and then
	$$
	2\vert \widehat{\sigma}\vert+\sum\limits_{k\in I^+}\vert Z_k\vert c_k \Vert U_k\Vert_1\geq 2\vert \widehat{\sigma}\vert +\vert Z_{j_1}\vert c_{j_1} \Vert U_{j_1}\Vert_1>  2\vert \widehat{\sigma}\vert +\vert Z_{j_1}\vert c_{j_1}\varepsilon. 
	$$
	Thus
	$$
	\Vert U_j\Vert_\infty \geq \Vert U_j\Vert_1>\frac{ 2\vert \widehat{\sigma}\vert +\vert Z_{j_1}\vert c_{j_1}\varepsilon}{\vert Z_j\vert c_j\vert I^{-}\vert}\geq \frac{2\vert \widehat{\sigma}\vert+\min\limits_{k\in I^+}\vert Z_k\vert c_k\varepsilon}{\max\limits_{k\in I^-}\vert Z_k\vert c_k\vert I^{-}\vert}.
	$$	
	On the other hand, if $i\in I^-$, then there exists $j\in I^+$ such that
	$$
	\Vert U_j\Vert_\infty\geq  \Vert U_j\Vert_1>\varepsilon,
	$$
	since $U$ does not
	have $\varepsilon$-small co-ions.
	Therefore, for both cases $i\in I^+$ and $i\in I^{-}$, there exists $j\in\{1,\ldots,N\}$ such that $Z_i$ and $Z_j$ have different signs and
	$$
	\Vert U_j\Vert_\infty>q(\varepsilon),
	$$
	where $q(\varepsilon)$ is given in (\ref{q}). Moreover, there exists a subset $\mathcal{B}$ of $[\widehat{R}_j,1-\widehat{R}_j]$ with $m(\mathcal{B})>0$ where we have $U_j> q(\varepsilon)$.
	
	Let $\Delta\in X$ be a perturbation of $U$ with components
	\begin{eqnarray}
		&&\Delta_i=-\alpha\chi_{\mathcal{A}},\ \Delta_j=-\beta\chi_{\mathcal{B}},\notag \\
		&&\Delta_k=0\text{\ for $k=1,\ldots,N$ with $k\neq i$ and $k\neq j$}\notag ,\\
		\label{Delta2}
	\end{eqnarray}
	where $\alpha,\beta>0$ and $\chi_{\mathcal{A}}$ and $\chi_\mathcal{B}$ are the characteristic functions of $\mathcal{A}$ and $\mathcal{B}$.
	
	Let $V=U+\Delta$. We check the following conditions:
	\begin{itemize}
		\item [1)]$V_k\geq 0$, $k\in\{1,\ldots,N\}$, i.e., $V$ satisfies the non-negativity condition;
		\item [2)] $E(V)=0$, i.e., $V$ satisfies the electroneutrality condition;
		\item [3)] $\Vert V\Vert \leq \Vert U\Vert$;
		\item [4)] $\widehat{F}(V)<\widehat{F}(U)$.
	\end{itemize}
	
	Condition 1). We have $V\geq 0$ if $\alpha\leq M$ and $\beta\leq q(\varepsilon)$.
	
	Condition 2). We have $E(V)=0$ if and only if $[\mathcal{T}(\Delta)](1)=0$ and we have
	\begin{eqnarray*}
		[\mathcal{T}(\Delta)](1)=-\alpha Z_ic_im(\mathcal{A})-\beta Z_jc_jm(\mathcal{B})=0
	\end{eqnarray*} 
	if and only if
	$$
	\beta=k\alpha\text{\ with\ }k=-\frac{Z_ic_im(\mathcal{A})}{Z_jc_jm(\mathcal{B})}>0.
	$$
	
	For the next conditions, 3) and 4), we assume
	$$
	\alpha\leq\alpha_0:=\min\left\{M,\frac{q(\varepsilon)}{k}\right\}\text{\ and\ }\beta=k\alpha
	$$
	to have conditions 1) and 2) satisfied.

	Condition 3). We have $\Vert V\Vert\leq \Vert U\Vert$ if  $\Vert V_i\Vert_p\leq \Vert U_i\Vert_p$ and $\Vert V_j\Vert_p\leq \Vert U_j\Vert_p$. Observe that
	$$
	\Vert V_i\Vert_p^p-\Vert U_i\Vert_p^p=\int\limits_{x\in \mathcal{A}}(U_i(x)-\alpha)^pdx-\int\limits_{x\in \mathcal{A}}U_i(x)^pdx\leq 0
	$$
	and
	$$
	\Vert V_j\Vert_p^p-\Vert U_j\Vert_p^p=\int\limits_{x\in \mathcal{B}}(U_j(x)-\beta)^pdx-\int\limits_{x\in \mathcal{B}}U_j(x)^pdx\leq 0
	$$
	and thus we have $\Vert V\Vert \leq \Vert U \Vert$. 
	
	Condition 4). We write $\widehat{F}=\widehat{F}^{\mathrm{en}}+\widehat{F}^{\mathrm{el}}$, with $\widehat{F}^{\mathrm{el}}=KL$, as in the proofs of the previous lemma and Proposition \ref{Thconvex}. 	Regarding the entropic part $\widehat{F}^{\mathrm{en}}$, we have 
	\begin{eqnarray*}
		\widehat{F}^{\mathrm{en}}(V)-\widehat{F}^{\mathrm{en}}(U)&=&c_i\int\limits_{x\in \mathcal{A}}\nu(U_i(x)-\alpha)dx+c_j\int\limits_{x\in \mathcal{B}}\nu(U_j(x)-k\alpha)dx\\
		&&-\left(c_i\int\limits_{x\in \mathcal{A}}\nu(U_i(x))dx+c_j\int\limits_{x\in \mathcal{B}}\nu(U_j(x))dx\right).
	\end{eqnarray*}
	Regarding the electrostatic part $\widehat{F}^{\mathrm{el}}=KL$, we have
	\begin{eqnarray*}
		&&L(V)-L(U)\\
		&&=2\langle\widehat{\sigma}+\mathcal{T}(U),\mathcal{T}(\Delta)\rangle+\langle \mathcal{T}(\Delta),\mathcal{T}(\Delta)\rangle\\
		&&=2\langle\widehat{\sigma}+\mathcal{T}(U),Z_ic_i\mathcal{L}_i(-\alpha\chi_\mathcal{A})+Z_jc_j\mathcal{L}_j(-k\alpha \chi_\mathcal{B})\rangle\\
		&&\quad +\langle Z_ic_i\mathcal{L}_i(-\alpha\chi_\mathcal{A})+Z_jc_j\mathcal{L}_j(-k\alpha \chi_\mathcal{B}),Z_ic_i\mathcal{L}_i(-\alpha\chi_\mathcal{A})+Z_jc_j\mathcal{L}_j(-k\alpha \chi_\mathcal{B})\rangle\\
		&&=-2\alpha\langle\widehat{\sigma}+\mathcal{T}(U),Z_ic_i\mathcal{L}_i(\chi_\mathcal{A})+kZ_jc_j\mathcal{L}_j(\chi_\mathcal{B})\rangle\\
		&&\quad +\alpha^2\langle Z_ic_i\mathcal{L}_i(\chi_\mathcal{A})+kZ_jc_j\mathcal{L}_j(\chi_\mathcal{B}),Z_ic_i\mathcal{L}_i(\chi_\mathcal{A})+kZ_jc_j\mathcal{L}_j(\chi_\mathcal{B})\rangle.
	\end{eqnarray*}

	Therefore, there exists $\alpha_1\in\left(0,\alpha_0\right]$ such that, for $\alpha\in(0,\alpha_1]$,
	$$
	\widehat{F}(V)-\widehat{F}(U)=\widehat{F}^{\mathrm{en}}(V)-\widehat{F}^{\mathrm{en}}(U)+K(L(V)-L(U))<0
	$$
	if the function
	\begin{eqnarray*}
		h(\alpha)&=&c_i\int\limits_{x\in \mathcal{A}}\nu(U_i(x)-\alpha)dx+c_j\int\limits_{x\in \mathcal{B}}\nu(U_j(x)-k\alpha)dx\\ &&-2K\langle\widehat{\sigma}+\mathcal{T}(U),Z_ic_i\mathcal{L}_i(\chi_\mathcal{A})+kZ_jc_j\mathcal{L}_j(\chi_\mathcal{B})\rangle\alpha\\
		&&+K\langle Z_ic_i\mathcal{L}_i(\chi_\mathcal{A})+kZ_jc_j\mathcal{L}_j(\chi_\mathcal{B}),Z_ic_i\mathcal{L}_i(\chi_\mathcal{A})+kZ_jc_j\mathcal{L}_j(\chi_\mathcal{B})\rangle\alpha^2,\\
		&&\alpha\in\left[0,\alpha_0\right],
	\end{eqnarray*}
	has $h^\prime(0)<0$. We have
	\begin{eqnarray}
		h^\prime(0)&=&-c_i\int\limits_{x\in \mathcal{A}}\nu^\prime(U_i(x))dx-c_jk\int\limits_{x\in \mathcal{B}}\nu^\prime(U_j(x))dx\notag \\
		&&-2K\langle\widehat{\sigma}+\mathcal{T}(U),Z_ic_i\mathcal{L}_i(\chi_\mathcal{A})+kZ_jc_j\mathcal{L}_j(\chi_\mathcal{B})\rangle.\notag \\
		\label{dersum2}
	\end{eqnarray}
	Regarding the first term in (\ref{dersum2}), we have
	\begin{eqnarray*}
		&&-c_i\int\limits_{x\in \mathcal{A}}\nu^\prime(U_i(x))dx-c_jk\int\limits_{x\in \mathcal{B}}\nu^\prime(U_j(x))dx\\
		&&\leq -c_i\nu^\prime\left(M\right)m(\mathcal{A})-c_jk\nu^\prime\left(q(\varepsilon)\right)m(\mathcal{B})\\
		&&=-c_i\left(
		\nu'(M)+\frac{|Z_i|}{|Z_j|}\nu'(q(\varepsilon))
		\right)m(\mathcal{A}).
	\end{eqnarray*}
	Regarding the second term in (\ref{dersum2}), we have, as in the proof of the previous lemma,
	\begin{eqnarray*}
	&&\vert\langle\widehat{\sigma}+\mathcal{T}(U),Z_ic_i\mathcal{L}_i(\chi_\mathcal{A})+kZ_jc_j\mathcal{L}_j(\chi_\mathcal{B})\rangle \vert\\
	&&\leq \Vert \widehat{\sigma}+\mathcal{T}(U)\Vert_2\Vert Z_ic_i\mathcal{L}_i(\chi_\mathcal{A})+kZ_jc_j\mathcal{L}_j(\chi_\mathcal{B})\Vert_\infty
	\end{eqnarray*}
	with
	$$
	\Vert \widehat{\sigma}+\mathcal{T}(U)\Vert_2\leq \sqrt{\frac{\omega}{K}}
	$$
	and
	\begin{eqnarray*}
	\Vert Z_ic_i\mathcal{L}_i(\chi_\mathcal{A})+kZ_jc_j\mathcal{L}_j(\chi_\mathcal{B})\Vert_\infty&\leq& \vert Z_i\vert c_i \Vert \mathcal{L}_i(\chi_\mathcal{A})\Vert_\infty +k\vert Z_j\vert c_j \Vert  \mathcal{L}_j(\chi_\mathcal{B})\Vert_\infty\\
	&\leq&\vert Z_i\vert c_i \Vert \chi_\mathcal{A}\Vert_1 +k\vert Z_j\vert c_j \Vert  \chi_\mathcal{B}\Vert_1\\
	&=&\vert Z_i\vert c_i m(\mathcal{A}) +k\vert Z_j\vert c_j m(\mathcal{B})\\
	&=&2\vert Z_i\vert c_im(\mathcal{A}).
	\end{eqnarray*}
	Thus, for the second term in (\ref{dersum2}), we have
	\begin{eqnarray*}
	&&\vert -2K\langle\widehat{\sigma}+\mathcal{T}(U), Z_ic_i\mathcal{L}_i(\chi_\mathcal{A})+kZ_jc_j\mathcal{L}_j(\chi_\mathcal{B})\rangle\vert\\
	&&\leq 4\vert Z_i\vert c_i\sqrt{K\omega}\ m(\mathcal{A}).
	\end{eqnarray*}
Therefore
	\begin{equation*}
		h^\prime(0)\leq -c_i\left(\nu^\prime\left(M\right)+\frac{|Z_i|}{|Z_j|}\nu^\prime\left(q(\varepsilon)\right)-4\vert Z_i\vert \sqrt{K\omega}\right)m(\mathcal{A}).
	\end{equation*}
	We can conclude that if
	$$
	\nu^\prime\left(M\right)>r(\omega,\varepsilon),
	$$
	or, equivalently, 
	$$
	M>\left(\nu^\prime\right)^{-1}\left(r(\omega,\varepsilon)\right),
	$$
	where $r(\omega,\varepsilon)$ is defined in (\ref{M}), then there exists $\alpha_1\in\left(0,\alpha_0\right]$ such that, for $\alpha\in(0,\alpha_1]$, $\widehat{F}(V)<\widehat{F}(U)$.
	
	In summary, if  (\ref{cL2}) holds, then by taking $\alpha>0$ sufficiently small and
	$$
	\beta=-\frac{Z_ic_im(\mathcal{A})}{Z_jc_jm(\mathcal{B})}\alpha
	$$
	in (\ref{Delta2}), we obtain $V\in A$, $\Vert V\Vert\leq \Vert U\Vert $ and $\widehat{F}(V)<\widehat{F}(U)$ with $V=U+\Delta$. 
	
\end{proof}

Here is the third lemma. 
\begin{lem}  \label{Lemma3}
	Let $\varepsilon>0$. Let $U\in A$ be such that $U$ has $\varepsilon$-small co-ions. Let
	$$
	s(\varepsilon):=\max\left\{\frac{\varepsilon}{\min\limits_{k\in I^{+}}(1-2\widehat{R}_k)},\frac{2\vert \widehat{\sigma}\vert+\sum\limits_{k\in I^+}\vert Z_k\vert c_k \varepsilon}{\min\limits_{k\in I^{-}}\vert Z_k\vert c_k(1-2\widehat{R}_k)}\right\}.
	$$
	We have
	\begin{equation}
	\Vert U\Vert_\infty \leq \left(\nu^\prime\right)^{-1}\left(\nu^\prime\left(s(\varepsilon)\right)+\nu^\prime\operatorname*{-ess\,osc}(U)\right). \label{eqnLemma3}
	\end{equation}
\end{lem}

\begin{proof}
	If $\nu'\operatorname{-ess\,osc}(U)=+\infty$, the
	right-hand side of (\ref{eqnLemma3}) is $+\infty$,
	so the conclusion is immediate. We may therefore
	assume that $\nu'\operatorname{-ess\,osc}(U)<+\infty$.
	
	Let $i\in\{1,\ldots,N\}$ be such that $\Vert U\Vert_\infty=\Vert U_i\Vert_\infty$. Let
	$$
	a_i=\operatorname*{ess\,inf}_{x\in [\widehat{R}_i,1-\widehat{R}_i]} U_i(x)\text{\ and\ }b_i=\operatorname*{ess\,sup}_{x\in [\widehat{R}_i,1-\widehat{R}_i]} U_i(x)=\Vert U_i\Vert_\infty.
	$$
	We have (recall (\ref{diffnuprime}))
	\begin{equation}
	\nu^\prime(b_i)=\nu^\prime(a_i)+\operatorname*{ess\,osc}(\nu^\prime\circ U_i). \label{aibi}
	\end{equation}

	Suppose $i\in I^+$. We have
	$$
	\varepsilon\geq \Vert U_i\Vert_1\geq (1-2\widehat{R}_i)a_i
	$$
	and then
	$$
	\nu^\prime(a_i)\leq\nu^\prime\left(\frac{\varepsilon}{1-2\widehat{R}_i}\right),
	$$
	which implies
	\begin{eqnarray}
	\nu^\prime(\Vert U\Vert_\infty)=\nu^\prime(\Vert U_i\Vert_\infty)=\nu^\prime(b_i)&\leq&\nu^\prime\left(\frac{\varepsilon}{1-2\widehat{R}_i}\right)+\operatorname*{ess\,osc}(\nu^\prime\circ U_i)\notag \\
	&\leq&\nu^\prime\left(\frac{\varepsilon}{\min\limits_{k\in I^+}(1-2\widehat{R}_k)}\right)+\nu^\prime\operatorname*{-ess\,osc}(U) \notag\\
	\label{thefirst}
	\end{eqnarray}
	by (\ref{aibi}) and (\ref{nuprimeessosc}).
	
	Suppose $i\in I^-$. By the electroneutrality condition (\ref{encond2}), we have
	\begin{eqnarray}
		\vert Z_i\vert c_i \Vert U_i\Vert_1&\leq& \sum\limits_{k\in I^-}\vert Z_k\vert c_k \Vert U_k\Vert_1=2\vert \widehat{\sigma}\vert+\sum\limits_{k\in I^+}\vert Z_k\vert c_k \Vert U_k\Vert_1\notag \\
		&\leq& 2\vert \widehat{\sigma}\vert+\sum\limits_{k\in I^+}\vert Z_k\vert c_k \varepsilon.\notag\\
		 \label{encond3}
	\end{eqnarray}
	Moreover, we have
	$$
	\vert Z_i\vert c_i \Vert U_i\Vert_1\geq \vert Z_i\vert c_i(1-2\widehat{R}_i)a_i
	$$
	and then, by (\ref{encond3}),
	$$
	\nu^\prime(a_i)\leq\nu^\prime\left(\frac{2\vert \widehat{\sigma}\vert+\sum\limits_{k\in I^+}\vert Z_k\vert c_k \varepsilon}{\vert Z_i\vert c_i(1-2\widehat{R}_i)}\right),
	$$
	which implies
	\begin{eqnarray}
	\nu^\prime(\Vert U\Vert_\infty)=\nu^\prime(\Vert U_i\Vert_\infty)=\nu^\prime(b_i)&\leq&\nu^\prime\left(\frac{2\vert \widehat{\sigma}\vert+\sum\limits_{k\in I^+}\vert Z_k\vert c_k \varepsilon}{\vert Z_i\vert c_i(1-2\widehat{R}_i)}\right)+\operatorname*{ess\,osc}(\nu^\prime \circ U_i)\notag\\
	&\leq&\nu^\prime\left(\frac{2\vert \widehat{\sigma}\vert+\sum\limits_{k\in I^+}\vert Z_k\vert c_k \varepsilon}{\min\limits_{k\in I^{-}}\vert Z_k\vert c_k(1-2\widehat{R}_k)}\right)+\nu^\prime\operatorname*{-ess\,osc}(U). \notag\\
	\label{thesecond}
	\end{eqnarray}
	
	By (\ref{thefirst}) and (\ref{thesecond}), we conclude that
	\begin{eqnarray*}
	\nu^\prime(\Vert U\Vert_\infty)&\leq& \max\left\{\nu^\prime\left(\frac{\varepsilon}{\min\limits_{k\in I^{+}}(1-2\widehat{R}_k)}\right),\nu^\prime\left(\frac{2\vert \widehat{\sigma}\vert+\sum\limits_{k\in I^+}\vert Z_k\vert c_k \varepsilon}{\min\limits_{k\in I^{-}}\vert Z_k\vert c_k(1-2\widehat{R}_k)}\right)\right\}\\
	&&+\nu^\prime\operatorname*{-ess\,osc}(U),
	\end{eqnarray*}
	i.e., (\ref{eqnLemma3}).	
\end{proof}

\subsection{Existence and uniqueness of the minimum}
We can now state the existence and uniqueness result for the
problem MP. 

\begin{Thm}
The problem MP has a unique minimum $\overline{U}$ and $\overline{U}_i\in L^\infty(\widehat{R}_i,1-\widehat{R}_i)$, $i\in\{1,\ldots,N\}$.
\end{Thm}

\begin{proof}
We show that there does not exist an increasing diverging sequence $(r_n)_{n\in\mathbb{N}^+}$ of positive numbers such that $r_n\geq r_0$, $n\in\mathbb{N}^+$, and $\Vert \overline{U}(r_n)\Vert\rightarrow +\infty$, $n\rightarrow \infty$. As a consequence, the alternative b) in Proposition \ref {twoalt} cannot hold. Thus, alternative a) holds, and MP has a unique minimum $\overline{U}$.

Suppose that such a sequence exists. Let $U_0\in A(r_0)$ and then we have $U_0\in A(r_n)$, $n\in\mathbb{N}^+$.

By applying Lemma \ref{Lemma1} with $\omega=\widehat{F}(U_0)$, we have
\begin{equation}
\nu^\prime\operatorname*{-ess\,osc}\left(\overline{U}(r_n)\right)\leq p\left(\widehat{F}(U_0)\right),\ n\in\mathbb{N}^+. \label{essosc<}
\end{equation}
By applying Lemma \ref{Lemma2} with $\omega=\widehat{F}\left(U_0\right)$ and $\varepsilon=1$, we have that $\overline{U}(r_n)$ has $1$-small co-ions, for $n\in\mathbb{N}^+$ sufficiently large. Finally, by Lemma \ref{Lemma3} as applied with $\varepsilon=1$ and (\ref{essosc<}), we obtain a contradiction, since $\Vert \overline{U}(r_n)\Vert \rightarrow +\infty$, $n\rightarrow \infty$.

We have proved that the alternative a) holds.

Now, we prove that $\overline{U}_i\in L^\infty(\widehat{R}_i,1-\widehat{R}_i)$, $i\in\{1,\ldots,N\}$. By applying all three lemmas with $\omega=\widehat{F}\left(U_0\right)$ and $\varepsilon=1$, we obtain $\Vert \overline{U}\Vert_\infty<+\infty$. In fact, by Lemma \ref{Lemma1} we have
\begin{equation}
\nu^\prime\operatorname*{-ess\,osc}(\overline{U})\leq p\left(\widehat{F}(U_0)\right), \label{essoscU<}
\end{equation}
by Lemma \ref{Lemma2} we have
$$
\Vert \overline{U}\Vert_\infty \leq \left(\nu^\prime\right)^{-1}\left(r\left(\widehat{F}(U_0),1\right)\right) 
$$
when $\overline{U}$ does not
have $1$-small co-ions, and by Lemma \ref{Lemma3} and (\ref{essoscU<}) we have
$$
\Vert \overline{U}\Vert_\infty \leq \left(\nu^\prime\right)^{-1}\left(\nu^\prime\left(s(1)\right)+\nu^\prime\operatorname*{-ess\,osc}(\overline{U})\right)\leq \left(\nu^\prime\right)^{-1}\left(\nu^\prime\left(s(1)\right)+p\left(\widehat{F}(U_0)\right)\right) 
$$
when $\overline{U}$ has $1$-small co-ions.
\end{proof}

It is important to observe that the minimum $\overline{U}$ is the same for all
$p\in(1,+\infty)$. In fact, given $p,q\in(1,+\infty)$ with $p<q$, let
$\overline{U}$ be the minimum obtained by using the space $L^p$. Since
$\overline{U}_i\in L^\infty(\widehat{R}_i,1-\widehat{R}_i)$,
$i\in\{1,\ldots,N\}$, we have $\overline{U}\in A_q\subset A_p$, where $A_q$ and $A_p$ denote the
domains corresponding to the spaces $L^q$ and $L^p$, respectively. Therefore, $\overline{U}$ is also the minimum obtained by using the space $L^q$.

\subsection{Back to dimensional quantities}
The problem MP is the dimensionless formulation of the problem of minimizing the free energy $F$ defined in~(\ref{F}) with respect to the ion concentrations $n_i$, $i\in\{1,\ldots,N\}$.

The dimensional concentrations $n_i^{\min}$, $i\in\{1,\ldots,N\}$, minimizing the free energy are given by
$$
n_i^{\min}(x)=\overline{U}_i\left(\frac{x}{D}\right)n_{i0},\ x\in[R_i,D-R_i],
$$
as follows from Subsection \ref{nondim}. Moreover, the minimal free energy $F^{\min}$ is given by
$$
\frac{F^{\min}}{Ak_{B}T}=Dn_0\widehat{F}(\overline{U}).
$$

\subsection{Symmetry}\label{Sym}

This subsection aims to show that the minimum of the problem MP has
components that are symmetric with respect to the point $\frac{1}{2}$, which is the midpoint of the interval on
which they are defined.

Let $\mathcal{S}:X\rightarrow X$ be the \emph{symmetry operator} given by
\begin{equation*}
	[\mathcal{S}(U)]_i(x)=U_i(1-x),\ i\in\{1,\ldots,N\}\text{\ and\ }x\in[\widehat{R}_i,1-\widehat{R}_i].
\end{equation*}
Observe that $[\mathcal{S}(U)]_i$, $i\in\{1,\ldots,N\}$, is the reflection of $U_i$ across the midpoint $\frac{1}{2}$ of the interval $[\widehat{R}_i,1-\widehat{R}_i]$.
Moreover, observe that $\mathcal{S}$ is a bijection on $X$ and the restriction $\mathcal{S}|_{X_+}$ is a bijection on $X_+$. 

\begin{pro} \label{PropTS}
For $U\in X$, we have
\begin{equation}
[\mathcal{T}(\mathcal{S}(U))](x)=[\mathcal{T}(U)](1)-[\mathcal{T}(U)](1-x),\ x\in [0,1]. \label{TS}
\end{equation}
\end{pro}

\begin{proof}
Let $i\in\{1,\ldots,N\}$. For $x\in [0,1]$, we have (recall (\ref{def}), or (\ref{R=0}) for $\widehat{R}_i=0$)
\begin{eqnarray*}
	[\mathcal{L}_{i}([\mathcal{S}(U)]_i)] \left( x\right) &=&
	\int\limits_{0}^{x}\frac{1}{2\widehat{R}_i}\int\limits_{z-\widehat{R}_{i}}^{z+\widehat{R}_{i}}
	[\mathcal{S}(U)]_i\left(s\right)dsdz=\int\limits_{0}^{x}\frac{1}{2\widehat{R}_i}\int\limits_{z-\widehat{R}_{i}}^{z+\widehat{R}_{i}}
	U_i\left(1-s\right)dsdz\\
	&=&\int\limits_{0}^{x}\frac{1}{2\widehat{R}_i}\int\limits_{1-z-\widehat{R}_{i}}^{1-z+\widehat{R}_{i}}
	U_i\left(t\right)dtdz=\int\limits_{1-x}^{1}\frac{1}{2\widehat{R}_i}\int\limits_{w-\widehat{R}_{i}}^{w+\widehat{R}_{i}}
	U_i\left(t\right)dtdw\\
	&=&[\mathcal{L}_{i}(U_i)] \left( 1\right)-[\mathcal{L}_{i}(U_i)] \left( 1-x\right).
\end{eqnarray*}
Now, (\ref{TS}) follows by recalling (\ref{loT}).

\end{proof}

The symmetry operator leaves the nondimensional free energy and the electroneutrality condition invariant.

\begin{pro} \label{invariant}
	For $U\in A$, we have $\widehat{F}(\mathcal{S}(U))=\widehat{F}(U)$ and $E(\mathcal{S}(U))=0$.
\end{pro}

\begin{proof}
Regarding the entropic part $\widehat{F}^{\mathrm{en}}$ in (\ref{FenhatU}) of $\widehat{F}$, we have, for $i\in\{1,\ldots, N\}$, 
\begin{eqnarray}
\langle \nu \circ [\mathcal{S}(U)]_i \rangle_i &=&\int\limits_{\widehat{R}_i}^{1-\widehat{R}_i} \nu([\mathcal{S}(U)]_i(x))dx=\int\limits_{\widehat{R}_i}^{1-\widehat{R}_i} \nu(U_i(1-x))dx\notag\\
&=&\int\limits_{\widehat{R}_i}^{1-\widehat{R}_i} \nu(U_i(x))dx=\langle \nu \circ U_i \rangle_i.\notag\\
 \label{ENP}
\end{eqnarray}
Regarding the electrostatic part $\widehat{F}^{\mathrm{el}}$ in (\ref{FelhatU}) of $\widehat{F}$, by Proposition \ref{PropTS} and the electroneutrality condition
$$
E(U)=2\widehat{\sigma}+[\mathcal{T}(U)](1)=0,
$$
we have
\begin{eqnarray}
\langle \widehat{\sigma}+\mathcal{T}(\mathcal{S}(U)),\widehat{\sigma}+\mathcal{T}(\mathcal{S}(U))\rangle &=&\int\limits_{0}^{1} \left(\widehat{\sigma}+[\mathcal{T}(\mathcal{S}(U))](x)\right)^2dx 
\notag\\
&=&\int\limits_{0}^{1} \left(\widehat{\sigma}-2\widehat{\sigma}-[\mathcal{T}(U)](1-x)\right)^2dx\notag \\
&=&\int\limits_{0}^{1} \left(\widehat{\sigma}+[\mathcal{T}(U)](1-x)\right)^2dx\notag\\ 
&=&\int\limits_{0}^{1} \left(\widehat{\sigma}+[\mathcal{T}(U)](x)\right)^2dx\notag \\
&=&\langle \widehat{\sigma}+\mathcal{T}(U),\widehat{\sigma}+\mathcal{T}(U)\rangle. \notag \\
\label{ELP}
\end{eqnarray}
By (\ref{ENP}) and (\ref{ELP}), we have $\widehat{F}(\mathcal{S}(U))=\widehat{F}(U)$. 

Regarding the electroneutrality condition, we have
$$
E(\mathcal{S}(U))=2\widehat{\sigma}+[\mathcal{T}(\mathcal{S}(U))](1)=2\widehat{\sigma}+[\mathcal{T}(U)](1)=0
$$
by Proposition \ref{PropTS}.
\end{proof}

Let $U\in X$. We say that $U$ is \emph{symmetric} if $\mathcal{S}(U)=U$. Here is the main result of this subsection.

\begin{Thm}
	The minimum $\overline{U}$ of MP is symmetric.
\end{Thm}
\begin{proof}
	By Proposition \ref{invariant}, we have $\mathcal{S}(\overline{U})\in A$ and $\widehat{F}(\mathcal{S}(\overline{U}))=\widehat{F}(\overline{U})$. Therefore, $\mathcal{S}(\overline{U})$ is a minimum of the problem MP. Since the minimum is unique, we have $\mathcal{S}(\overline{U})=\overline{U}$.
\end{proof}

The fact that the minimum $\overline{U}$ is symmetric, i.e., each component $\overline{U}_i$, $i\in\{1,\ldots,N\}$, is unchanged under reflection across the midpoint $\frac{1}{2}$ of the interval $[\widehat{R}_i,1-\widehat{R}_i]$, is not surprising. The underlying physical problem, in which we consider the equilibrium of an ionic solution confined between two parallel planes carrying the same surface charge density, is itself invariant under reflection across the midplane between the planes.

\section{Numerical tests}\label{NTsection}

In this section, we report the results of numerical tests for three different examples of the PB problem considered in the present paper.

The three examples are solved numerically by discretizing the
infinite-dimensional dimensionless problem MP and then computing a minimum of the resulting
finite-dimensional minimization problem in MATLAB, using the interior-point
method implemented in the function \texttt{fmincon}. The evaluation of the objective function $\widehat{F}$
given in (\ref{FhatU}) requires numerical integration, which is performed using the trapezoidal rule.

The tests are intended only to illustrate the qualitative behavior of the solutions of the minimization problem MP. A detailed description and analysis of the discretization, including its convergence properties, will be provided in a companion paper.

In each numerical test, we show the following quantities.
\begin{itemize}
	\item The \emph{concentration profiles}, i.e., the concentrations $n_i^{\mathrm{min}}$, $i\in\{1,\ldots,N\}$, of the different types of ions, as functions of
	the coordinate $x$,  minimizing the free energy $F$ given in (\ref{F}), 
	for a given distance $D$ between the two parallel planes. Due to the symmetry
	property described in Subsection~\ref{Sym}, each concentration profile $n_i^{\mathrm{min}}$ is shown only in the left half-interval
	from $R_i$ to $\frac{D}{2}$; since the model excludes penetration of the ions
	into the planes, the concentration profile starts at $x=R_i$. We use $D=20\ \mathrm{nm}$ for all three tests.
	\item The minimal free energy $F^{\min}$ as a function of $D$.
	\item The pressure $P$ defined in (\ref{pressure}), as a function of $D$. It is computed numerically using a central difference
	approximation to the derivative of $F^{\mathrm{min}}$. 
\end{itemize}
The phenomenon of like-charge attraction between the planes occurs at 
distances $D$ for which the pressure $P(D)$ is negative.

In the tests, the surface charge density of the planes is positive and is determined by the
\emph{surface parameter} $a$ through $\sigma=\frac{e}{a}$. The parameter $a$ represents the surface area carrying one elementary charge. The Bjerrum length $l_B$ in \eqref{eq1} is taken to be
$0.7\ \mathrm{nm}$ in all three examples.

\subsection{Test 1}
In the first test, we consider two types of spherical ions (nanoparticles) of equal size, absolute valence, and bulk concentration, with
the first type negatively charged and the second positively charged. This PB problem was addressed in \cite{S}, where the boundary
value problem for the potential obtained through the variational
procedure was solved numerically. The
physical parameters corresponding to Test 1 are reported in
Table~\ref{tab:test1}.

\begin{table}
	\centering
	\caption{Physical parameters used in Test 1.}
	\label{tab:test1}
	\begin{tabular}{ll}
		\hline
		\textbf{Parameter} & \textbf{Value} \\
		\hline
		$D$ & 20 nm \\
		$a$ & $5\,\mathrm{nm}^2$ \\
		$N$ & 2 \\
		$R_1$ & 1 nm \\
		$Z_1$ & $-12$ \\ 
		$n_{10}$ & $1\cdot 10^{-3}\,\mathrm{nm}^{-3}$\\
		$R_2$ & 1 nm \\
		$Z_2$ & $12$ \\
		$n_{20}$ & $1\cdot10^{-3}\,\mathrm{nm}^{-3}$ \\
		\hline
	\end{tabular}
\end{table}

\begin{figure}
	\begin{center}
		\includegraphics[width=.65\textwidth]{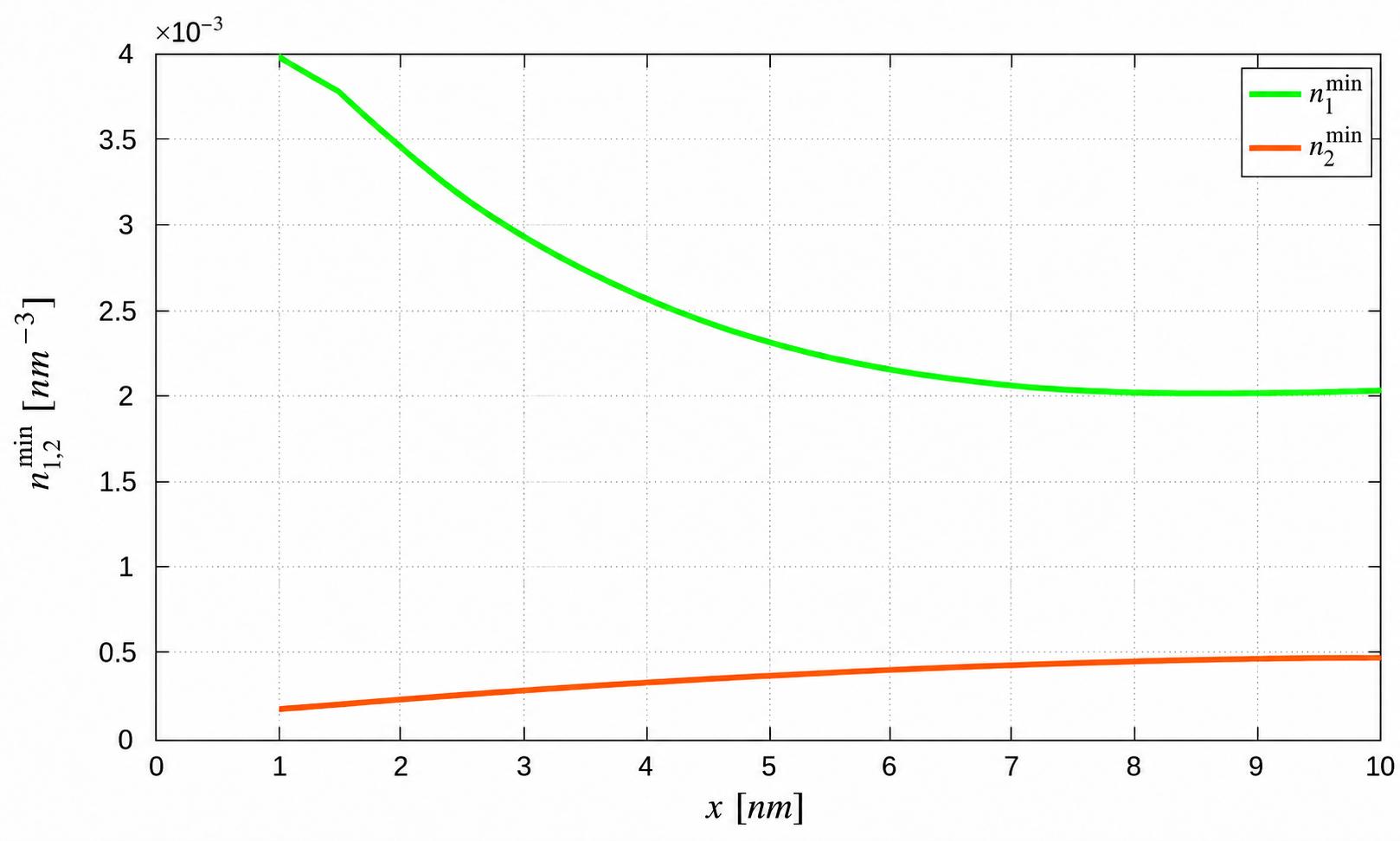}
		\caption{\label{f1} Test 1. Concentration profiles of negatively charged nanoparticles ($n_{1}^{\min}$) and positively charged nanoparticles ($n_{2}^{\min}$). The physical parameters are specified in Table \ref{tab:test1}.}
	\end{center}
\end{figure}

\begin{figure}
	\begin{center}
		\includegraphics[width=.65\textwidth]{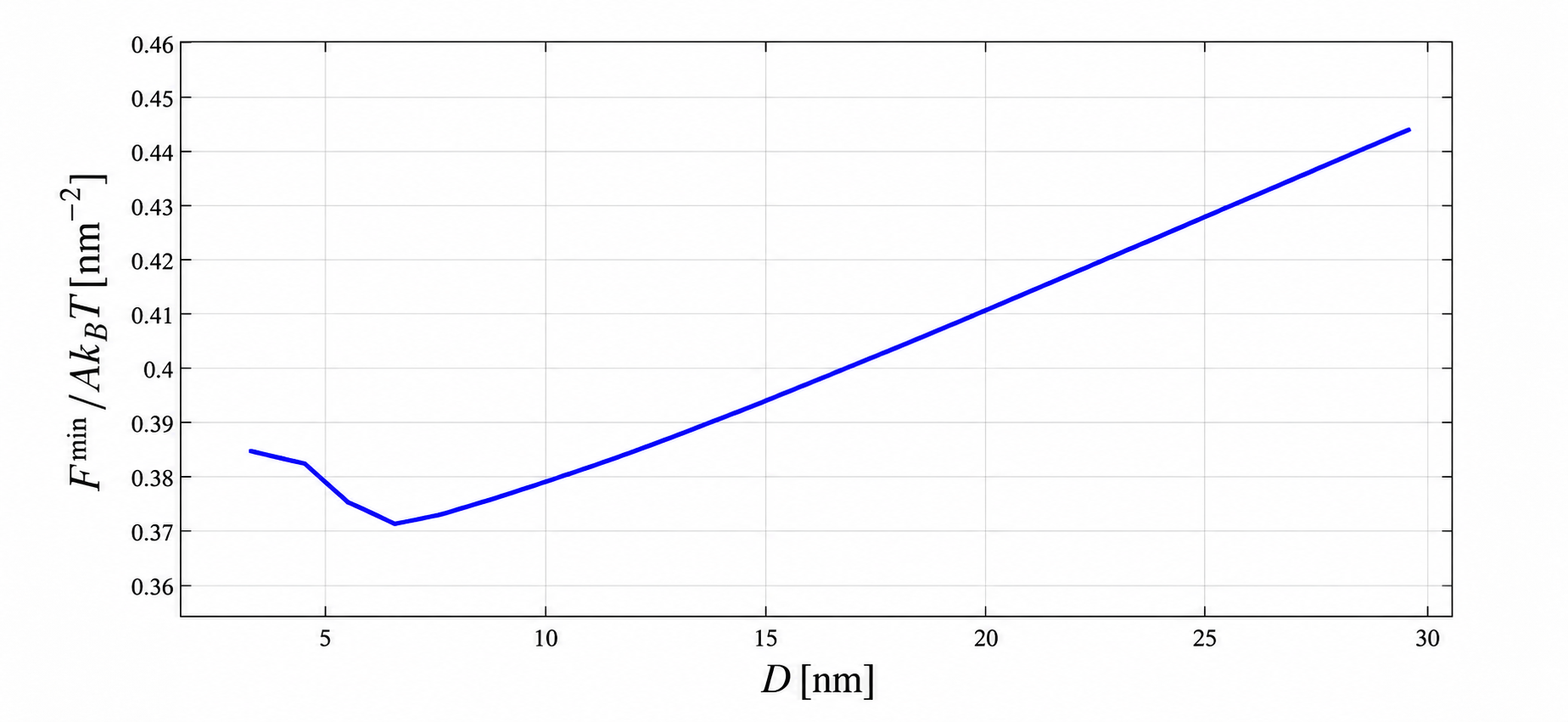}
		\includegraphics[width=.65\textwidth]{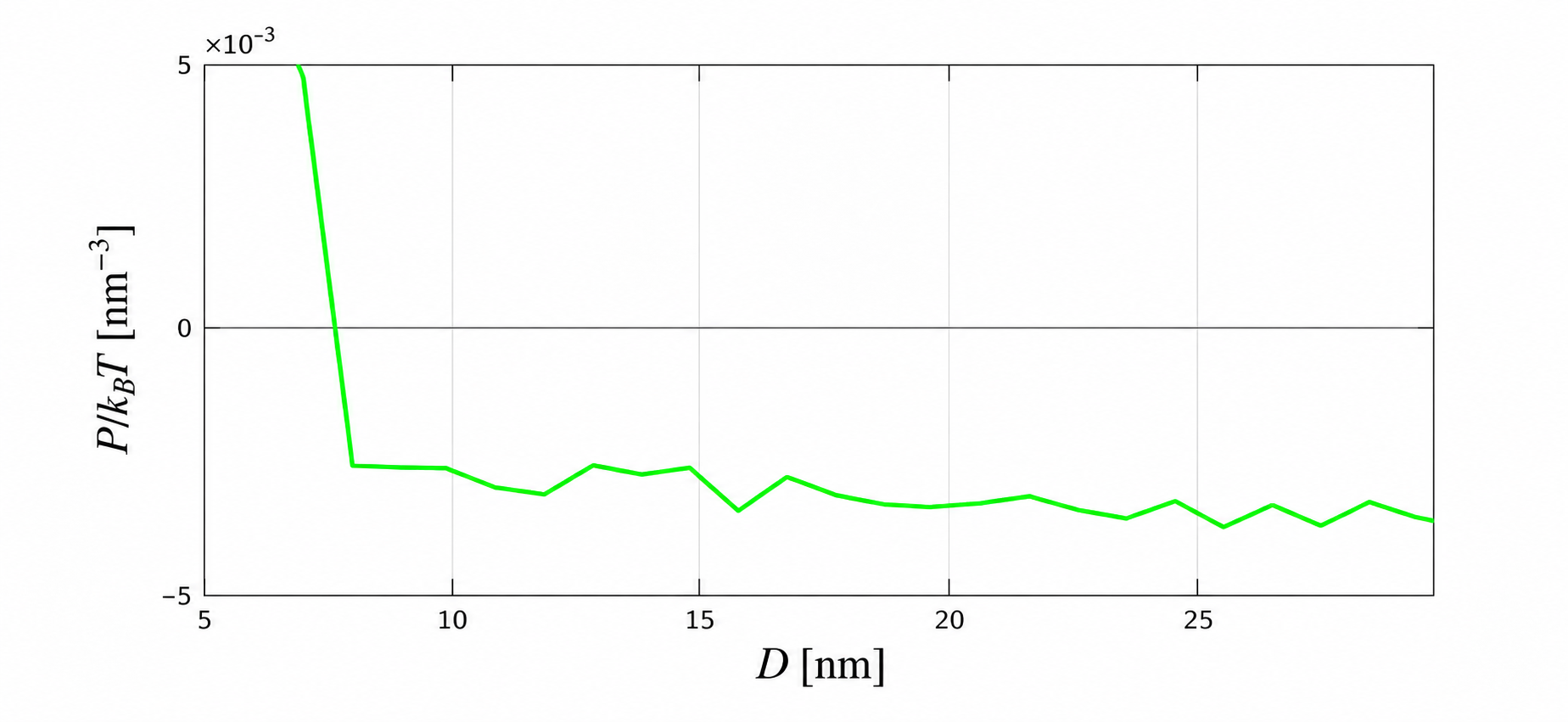}
		\caption{\label{f2} Test 1. Minimal free energy (top) and pressure (bottom) as  functions of $D$. The physical parameters are specified in Table \ref{tab:test1}.}
	\end{center}
\end{figure}

Figure \ref{f1} shows the concentration profiles $n_{1}^{\min}$ and $n_{2}^{\min}$ of the negatively and positively charged nanoparticles, respectively, when $D=20\ \mathrm{nm}$. We see that the concentration $n_{1}^{\min}$ of counter-ions is larger than the concentration $n_{2}^{\min}$ of co-ions. Moreover, near the charged surfaces the concentration of counter-ions is higher. 
Figure \ref{f2} shows that there is like-charge attraction between the planes, and it appears for distances $D$ larger than about $7.5\ \mathrm{nm}$.  

\subsection{\textbf{Test 2}}
In the second test, we consider a mixture of negatively charged nanoparticles  and  positively charged point-like salt ions, with a valence much larger than the absolute valence of the nanoparticles. A related Poisson--Boltzmann problem, with the same types of ions but
with a different entropic term, was considered in \cite{Bohinc2023}.  The
physical parameters corresponding to Test 2 are reported in
Table \ref{tab:test2}. 

\begin{table}[ht]
	\centering
	\caption{Physical parameters used in Test 2.}
	\label{tab:test2}
	\begin{tabular}{ll}
		\hline
		\textbf{Parameter} & \textbf{Value} \\
		\hline
		$D$ & 20 nm \\
		$a$ & $50\,\mathrm{nm}^2$ \\
		$N$ & 2 \\
		$R_1$ & 1 nm \\
		$Z_1$ & $-1$ \\
		$n_{10}$ & $1\cdot10^{-3}\,\mathrm{nm}^{-3}$ \\
		$R_2$ & 0  (point-like ions) \\
		$Z_2$ & $10$ \\
		$n_{20}$ & $1\cdot10^{-4}\,\mathrm{nm}^{-3}$\\		
		\hline
	\end{tabular}
\end{table}

\begin{figure}
	\begin{center}
		\includegraphics[width=.65\textwidth]{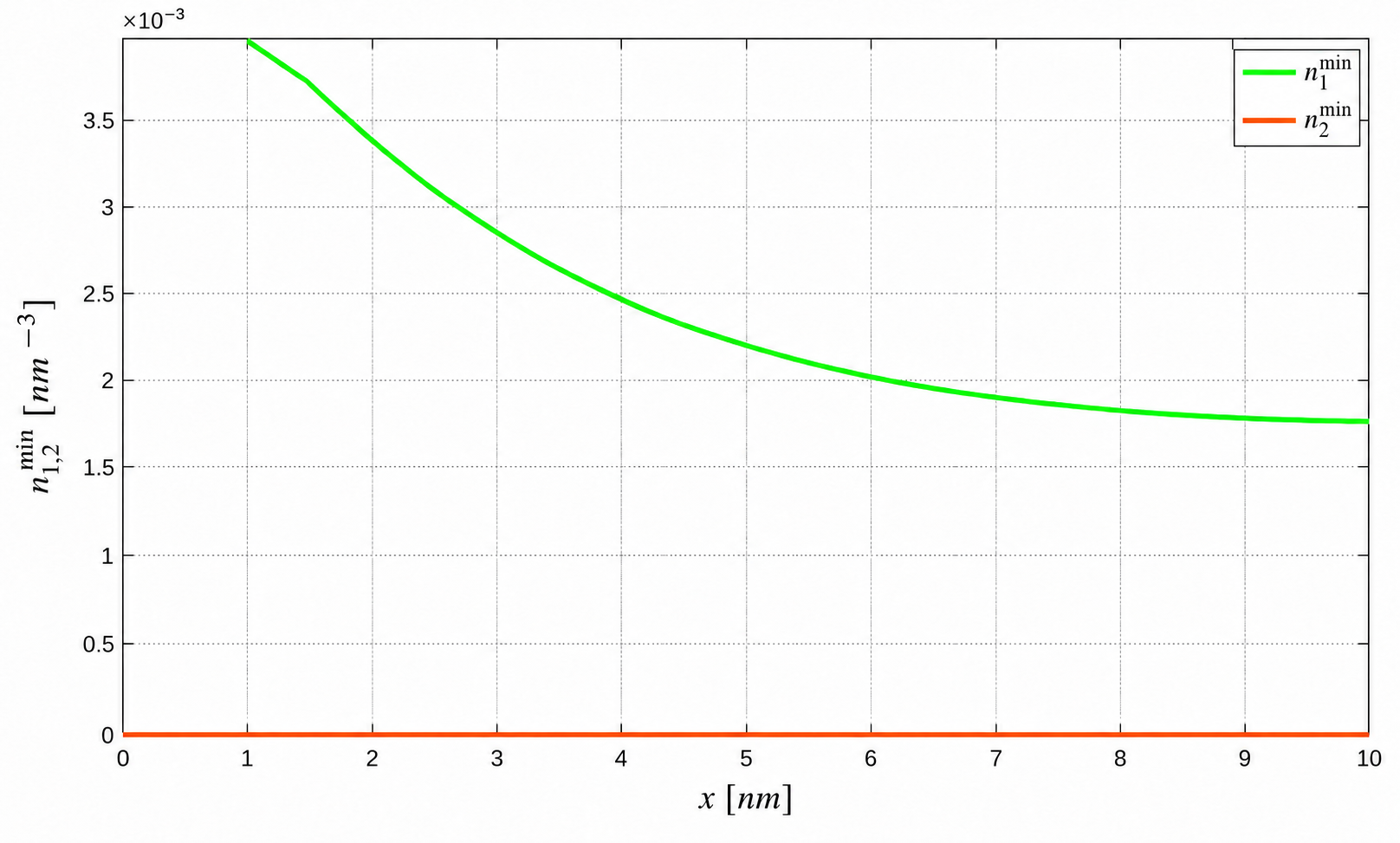}
	
		\caption{ \label{f3} Test 2. Concentration profiles of negatively charged nanoparticles ($n_1^{\mathrm{min}}$) and positively charged salt ions ($n_2^{\mathrm{min}}$). The physical parameters are specified in Table \ref{tab:test2}.}
\end{center}
\end{figure}

\begin{figure}
\begin{center}
	\includegraphics[width=.65\textwidth]{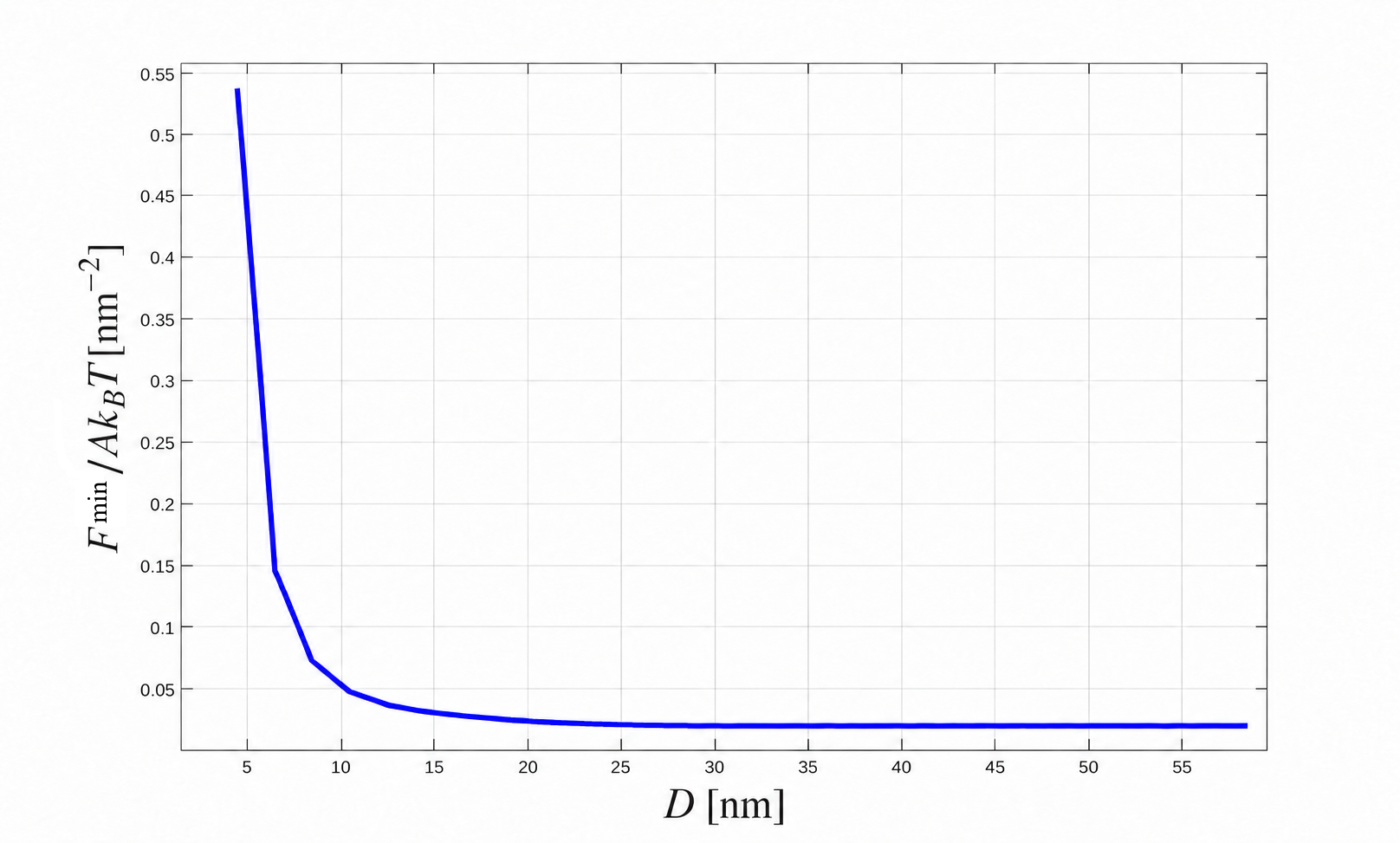}
	\includegraphics[width=.65\textwidth]{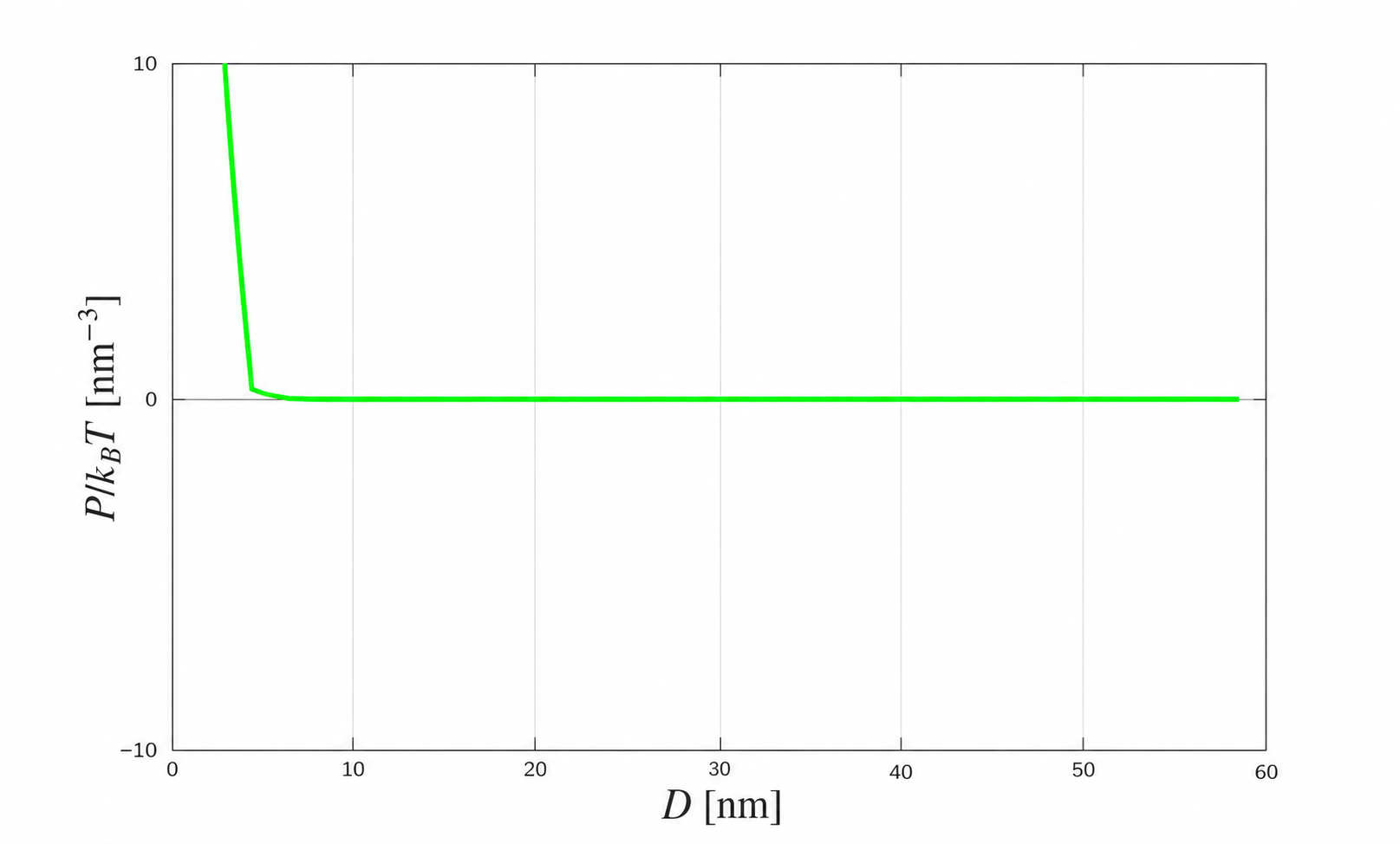}
	
	\caption{ \label{ff4} Test 2. Minimal free energy (top) and pressure (bottom) as  functions of $D$. The physical parameters are specified in Table \ref{tab:test2}.}
\end{center}
\end{figure}

Figure~\ref{f3} depicts the concentration profiles $n_1^{\mathrm{min}}$ and $n_2^{\mathrm{min}}$ of counter-nanoparticles and salt
co-ions, respectively, when $D=20\ \mathrm{nm}$. The
concentration profile $n_1^{\mathrm{min}}$ of the counter-nanoparticles is much larger than the concentration profile $n_2^{\mathrm{min}}$ of salt co-ions. Figure \ref{ff4} indicates that no like-charge attraction occurs between the planes over the range of values of $D$ considered.


\subsection{\textbf{Test 3}}
In the third and last test, we consider both negatively and positively charged nanoparticles of equal size and absolute valence, but with different bulk concentrations, and positively charged point-like salt ions with a valence much larger than the absolute valence of the nanoparticles. This type of PB problem has not been
addressed in the literature. The
physical parameters corresponding to Test 3 are reported in
Table \ref{tab:test3}. 

\begin{table}
	\centering
	\caption{Physical parameters used in Test 3.}
	\label{tab:test3}
	\begin{tabular}{ll}
		\hline
		\textbf{Parameter} & \textbf{Value} \\
		\hline
		$D$ & 20 nm \\
		$a$ & $5\,\mathrm{nm}^2$ \\
		$N$ & 3 \\
		$R_1$ & 1 nm \\
		$Z_1$ & $-1$ \\
		$n_{10}$ & $1\cdot10^{-3}\,\mathrm{nm}^{-3}$ \\
		$R_2$ & 1 nm \\
		$Z_2$ & $1$ \\
		$n_{20}$ & $5\cdot 10^{-4}\,\mathrm{nm}^{-3}$ \\
		$R_3$  & 0  (point-like ions) \\
		$Z_3$ & $10$ \\
		$n_{30}$ & $5\cdot 10^{-5}\ \mathrm{nm}^{-3}$\\
		\hline
	\end{tabular}
\end{table}

\begin{figure}
\begin{center}
	\includegraphics[width=.65\textwidth]{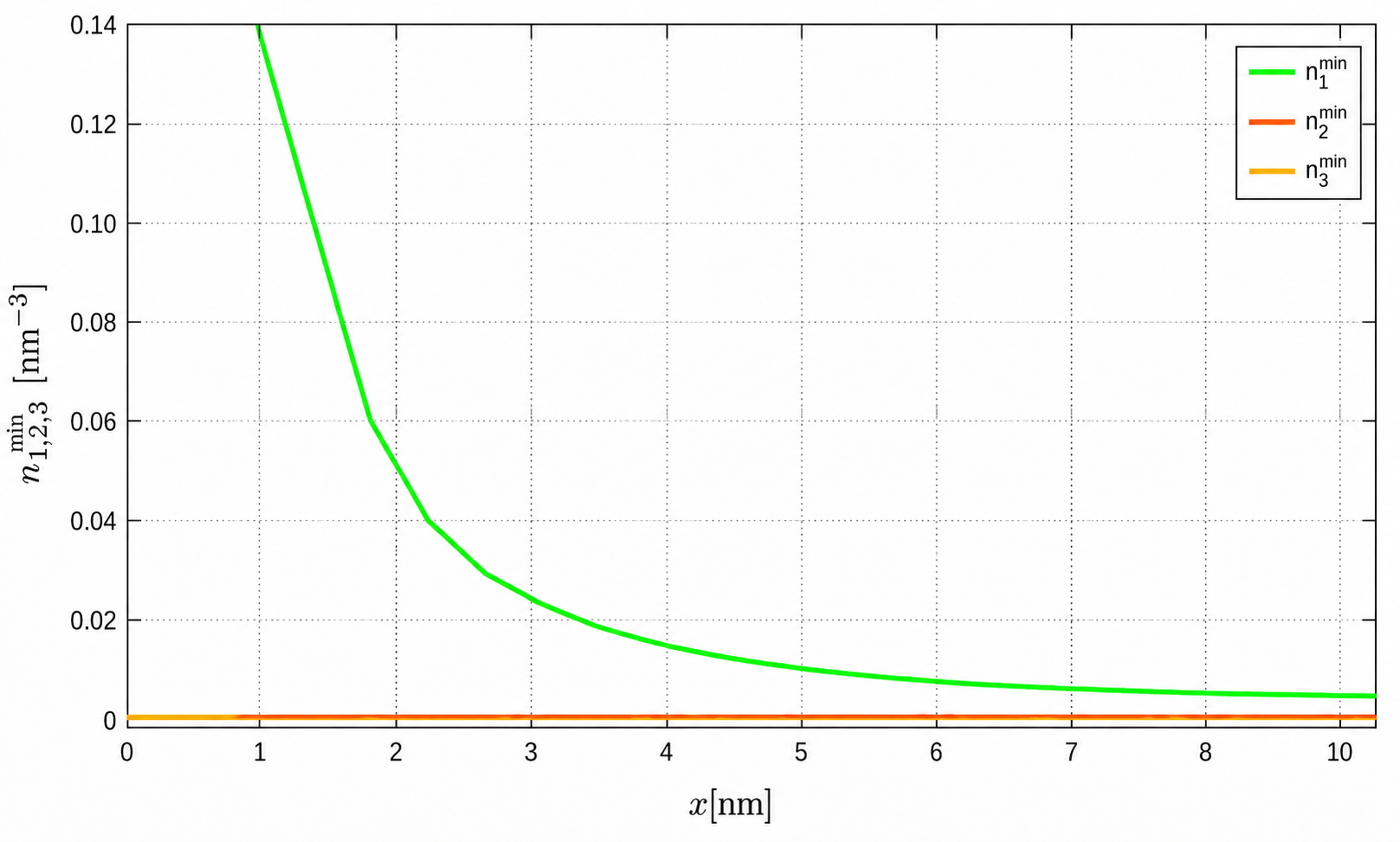}
	\caption{\label{ff5} Test 3. 
		Concentration profiles of the  
		negatively charged nanoparticles ($n_1^{\mathrm{min}}$), positively charged nanoparticles ($n_2^{\mathrm{min}}$) and positively charged salt ions ($n_3^{\mathrm{min}}$). The physical parameters are specified in Table \ref{tab:test3}.}
\end{center}
\end{figure}

\begin{figure}
\begin{center}
	\includegraphics[width=.65\textwidth]{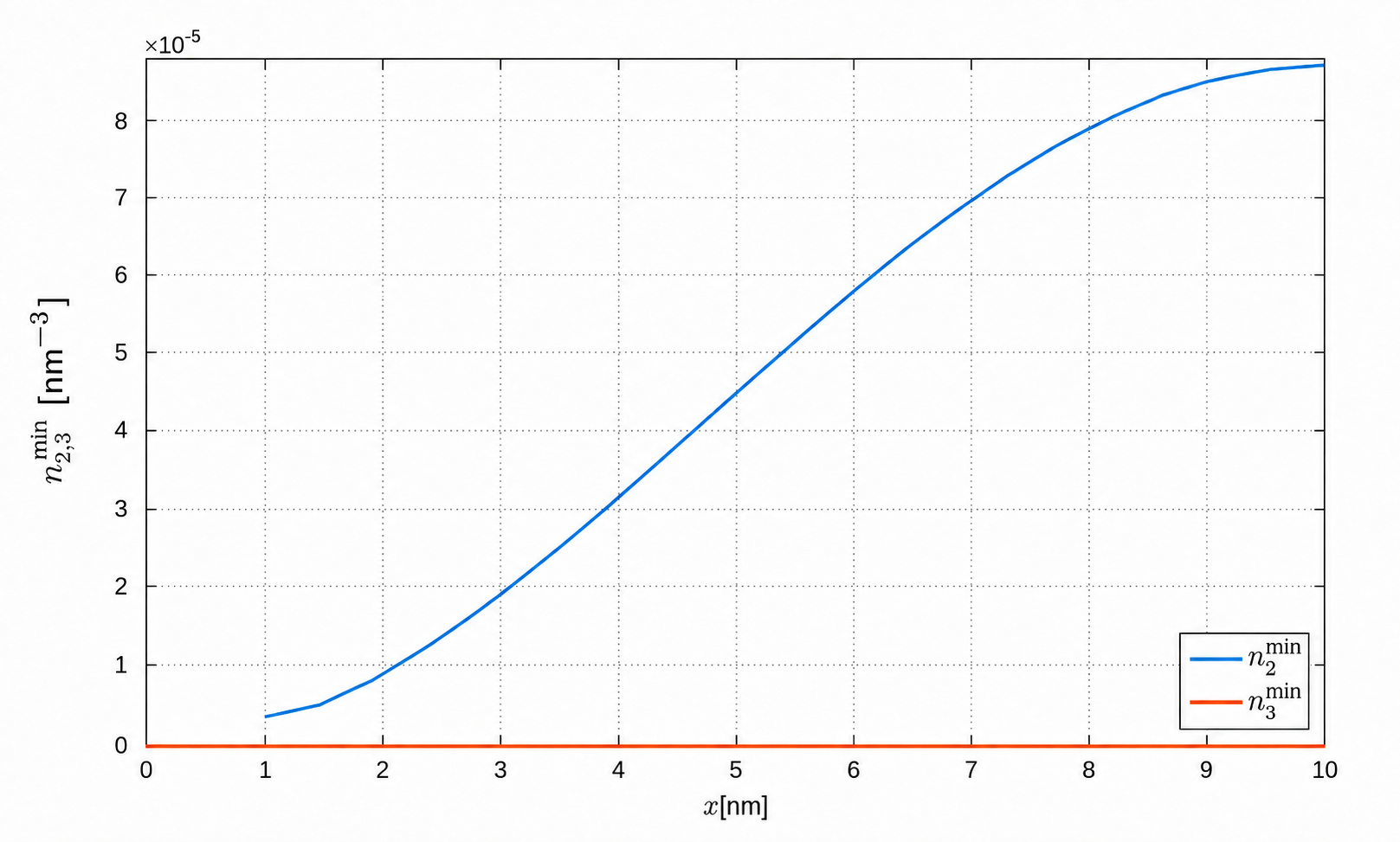}
	\caption{\label{ff6} Test 3. 
		Concentration profiles of the positively charged nanoparticles ($n_2^{\mathrm{min}}$) and positively charged salt ions ($n_3^{\mathrm{min}}$). The physical parameters are specified in Table \ref{tab:test3}.}
\end{center}
\end{figure}

\begin{figure}
\begin{center}
	\includegraphics[width=.65\textwidth]{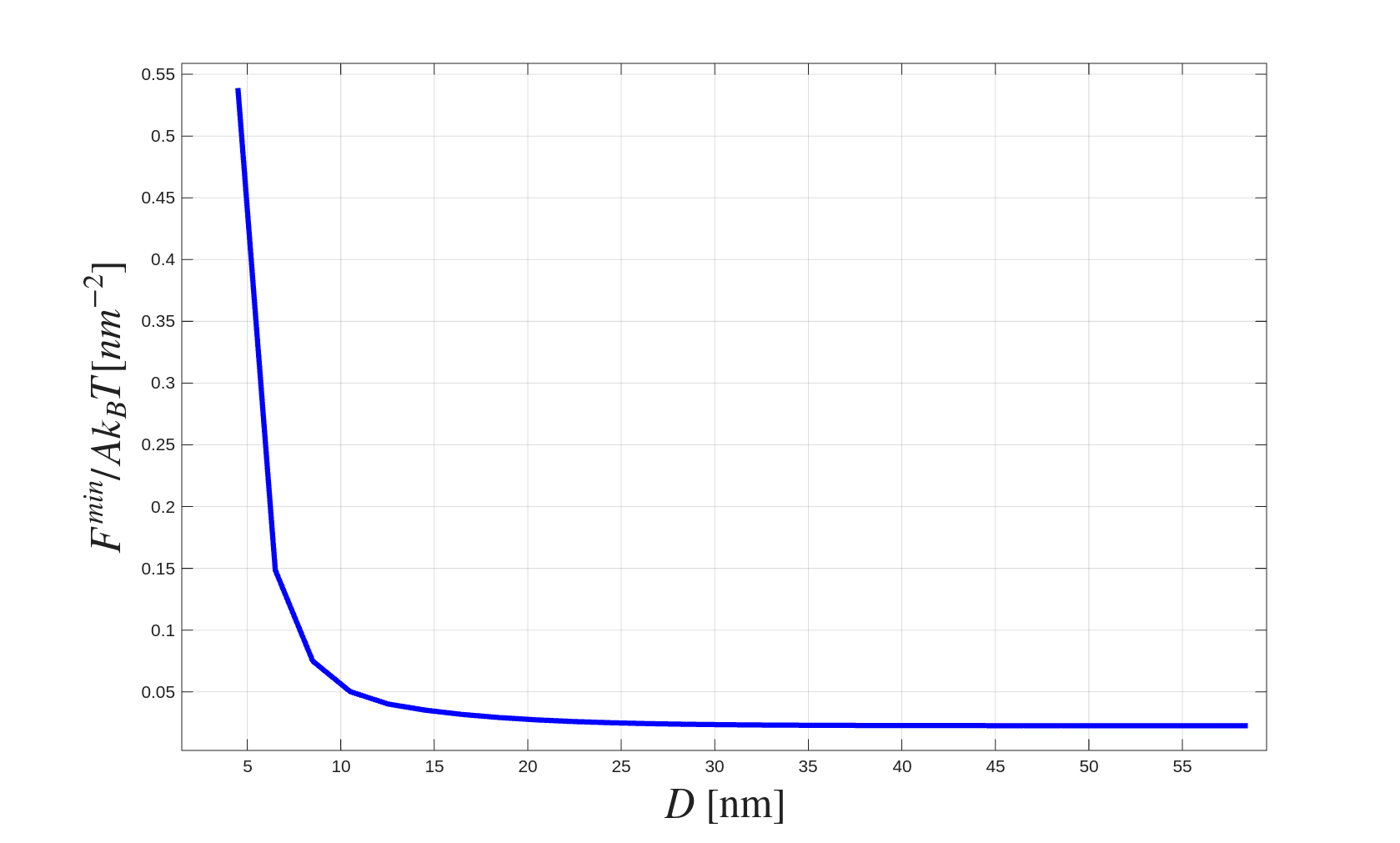}
	\includegraphics[width=.65\textwidth]{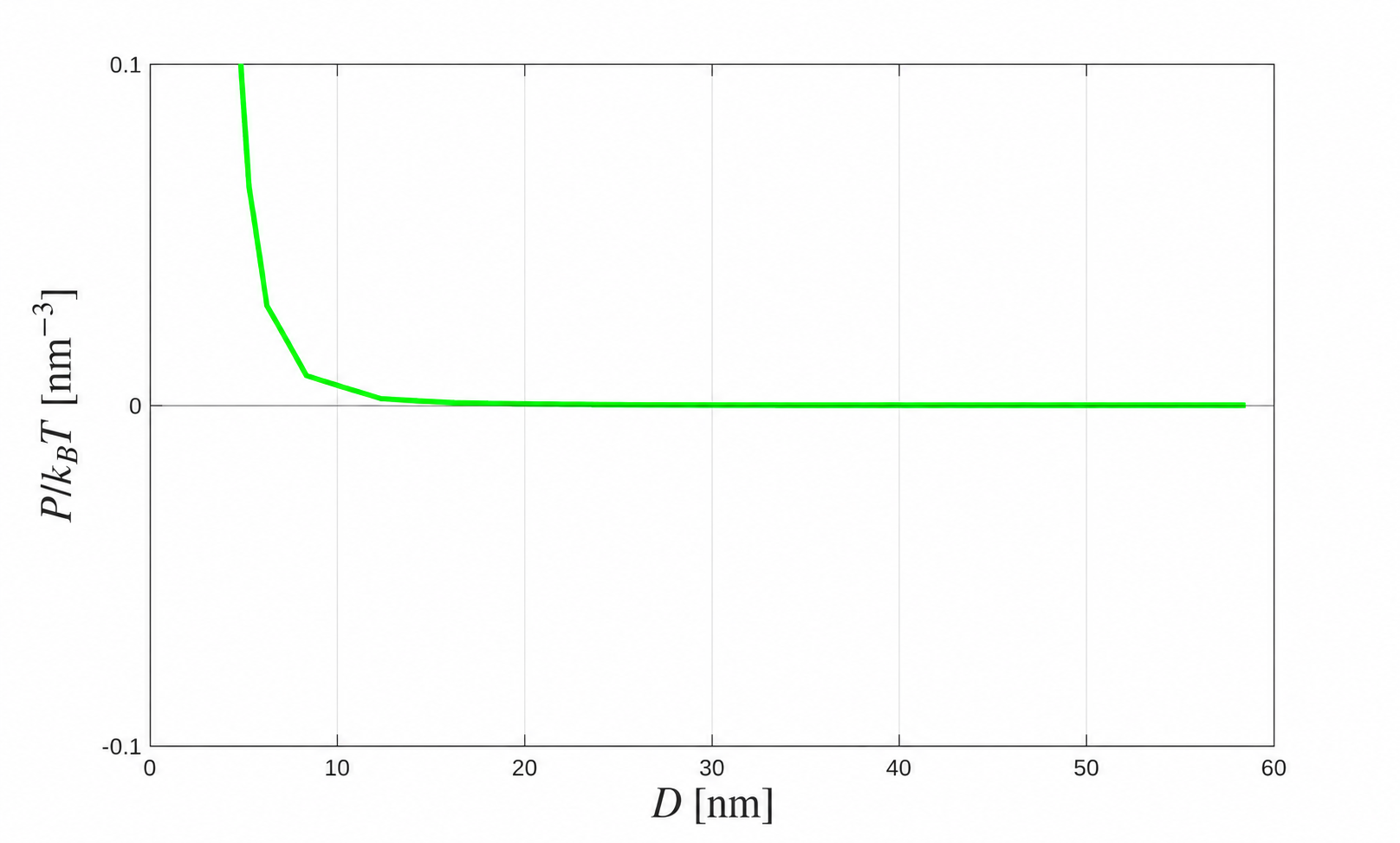}
	\caption{\label{ff7} Test 3. Minimal free energy (top) and pressure (bottom) as functions of $D$. The physical parameters are specified in Table \ref{tab:test3}.}
\end{center}
\end{figure}

Figure \ref{ff5} presents the concentration profiles $n_1^{\mathrm{min}}$, $n_2^{\mathrm{min}}$ and $n_3^{\mathrm{min}}$ of counter-nanoparticles, co-nanoparticles and salt
co-ions, respectively, when $D=20\ \mathrm{nm}$. The
concentration profile $n_1^{\mathrm{min}}$ of the counter-nanoparticles is much larger than the concentration profiles $n_2^{\mathrm{min}}$ and $n_3^{\mathrm{min}}$ of co-nanoparticles and salt co-ions. Figure \ref{ff6} shows $n_2^{\mathrm{min}}$ and $n_3^{\mathrm{min}}$, where $n_2^{\mathrm{min}}$ is much larger than $n_3^{\mathrm{min}}$. Figure \ref{ff7} indicates that no like-charge attraction occurs between the planes over the range of values of $D$ considered. 


\section{Conclusion}

In the present paper, we have considered a Poisson–Boltzmann problem in which the free energy of an ionic solution confined between two charged parallel planes is minimized. The ions in the
solution, which are spherical and carry charges uniformly distributed on their
surfaces, undergo molecular motion modeled as that of an ideal gas and
influenced by the electric field generated by the charged planes and by the ions
themselves. The ions are of different types, characterized by different sphere
radii, valences, and bulk concentrations.

The problem of finding the ion concentrations minimizing the 
free energy is formulated as the dimensionless infinite-dimensional minimization problem
MP in Subsection~\ref{MPsection}. The main result of the present paper is the
existence and uniqueness of a minimum for MP; the proof is given in
Section~\ref{EUsection}. The proof uses a norm-convergence direct method based on the
Banach--Saks property, which is not a standard tool in direct methods. Moreover,
the proof is nonstandard and technically challenging because the functional to be minimized is not coercive in general. The arguments used in the proof also show that the minimum is essentially bounded.

Future work will concern the numerical analysis of the discretization of
the problem MP used to obtain the numerical results in Section~\ref{NTsection}.
Moreover, the present paper opens the way to treating other types of
Poisson--Boltzmann problems. For example, one may consider problems involving more sophisticated
models of the molecular motion of the ions, such as lattice-gas or
hard-sphere models, as well as other types of ions, such as rod-like
ions or Janus particles.

\

\appendix

\section{A Banach--Saks direct method of proof for infinite-dimensional minimization}

Let $X$ be a Banach space, let $A$ be a non-empty subset of $X$ and let $%
F: A\rightarrow \mathbb{R}$ be a functional.

Consider the problem of minimizing $F$, i.e., we are interested in elements $%
\overline{x}\in A$ such that 
\begin{equation*}
	F\left( \overline{x}\right) =\min_{x\in A}F\left( x\right) .
\end{equation*}%
Such an element $\overline{x}$ is called a \emph{minimum} of $F$.

In this appendix, we present a \emph{direct method}, i.e., not based on the derivation of the Euler--Lagrange equation, for proving existence (and uniqueness) of a minimum of $F$, to be used when $X$ is infinite-dimensional. It is based on the Banach--Saks property of the space $X$, uses only norm convergence and avoids weak-topology arguments. For the general theory of the Banach--Saks property, we refer to \cite{Diestel,Partington}.

The topics of Subsections~\ref{Asub1}, \ref{Asub2}, and \ref{Asub4} are
standard and can be found in classical textbooks (see, e.g., \cite{Rindler}).
The main result of Subsection~\ref{Asub3}, Theorem~\ref{existence2}, is derived
by standard arguments, but it has apparently never been explicitly stated in
the literature, although a related minimization result for convex functionals on bounded closed convex subsets of \(L^p\) spaces is given in \cite{Royden} and a closely related argument is used in \cite{ChengNakiboglu}. Theorem~\ref{BSC}
of Subsection~\ref{Asub3} is presented in \cite{ZF}, but the proof given in the
present appendix is shorter and simpler.

\subsection{Minimizing sequences}\label{Asub1}

Since $A$ is non-empty, 
\begin{equation*}
	\inf_{x\in A}F\left( x\right) \in \left[ -\infty ,+\infty \right)
\end{equation*}%
is defined.

\begin{de}
	A sequence $\left( x_{n}\right)_{n\in\mathbb{N}^+}$ in $A$ is called a \emph{minimizing
		sequence} if 
	\begin{equation*}
		\lim_{n\rightarrow \infty }F\left( x_{n}\right) =\inf_{x\in A}F\left(
		x\right) .
	\end{equation*}
\end{de}

Clearly, minimizing sequences exist. In fact:

\begin{itemize}
	\item if $\inf\limits_{x\in A}F\left( x\right) \in \mathbb{R}$, then, for
	any $n\in \mathbb{N}^{+}$, there exists $x_{n}\in A$ such that 
	\begin{equation*}
		\inf_{x\in A}F\left( x\right) \leq F\left( x_{n}\right) <\inf_{x\in
			A}F\left( x\right) +\frac{1}{n};
	\end{equation*}
	
	\item if $\inf\limits_{x\in A}F\left( x\right) =-\infty $, then, for any $%
	n\in \mathbb{N}^{+}$, there exists $x_{n}\in A$ such that 
	\begin{equation*}
		F\left( x_{n}\right) <-n.
	\end{equation*}
\end{itemize}

We are interested in the case in which a bounded minimizing sequence exists.
Clearly, if $A$ is bounded, then every minimizing sequence in $A$ is bounded.
If $A$ is unbounded, a coercivity condition ensures the
boundedness of minimizing sequences.

\begin{de}
	Assume that $A$ is unbounded. $F$ is called \emph{coercive} if for any sequence $\left( x_{n}\right)_{n\in\mathbb{N}^+}$ in $A$ such that $\lim\limits_{n\rightarrow \infty}\Vert
	x_n\Vert=+\infty$, we have 
	\begin{equation*}
		\lim_{n\rightarrow \infty }F\left( x_n\right) =+\infty.
	\end{equation*}
\end{de}

\begin{pro}
	Assume that $A$ is unbounded. If $F$ is coercive, then every minimizing sequence is
	bounded.
\end{pro}

\begin{proof}
	Suppose that $F$ is coercive. Let $\left( x_{n}\right)_{n\in\mathbb{N}^+}$ be a minimizing sequence.
	Suppose that $\left( x_{n}\right)_{n\in\mathbb{N}^+}$ is unbounded. This means that there exists a
	subsequence $\left( x_{n_{k}}\right)_{k\in\mathbb{N}^+} $ such that 
	\begin{equation*}
		\lim_{k\rightarrow \infty }\left\Vert x_{n_{k}}\right\Vert =+\infty .
	\end{equation*}
	Since $F$ is coercive, we have 
	\begin{equation*}
		\lim_{k\rightarrow \infty }F\left( x_{n_{k}}\right) =+\infty .
	\end{equation*}
	This contradicts 
	\begin{equation*}
		\lim_{k\rightarrow \infty }F\left( x_{n_{k}}\right) =\inf_{x\in A}F\left(
		x\right)\in[-\infty,+\infty) .
	\end{equation*}
\end{proof}

\subsection{Existence of a minimum for finite-dimensional minimization problems}\label{Asub2}

Under the assumption that the minimization problem is finite-dimensional, i.e., $X$ is finite-dimensional, the next theorem gives conditions under which a solution to the minimization
problem exists. It requires lower semicontinuity for $F$.

\begin{de}
	Let $x\in A$. $F$ is called \emph{lower semicontinuous} at $x$ if, for any
	sequence $\left( x_{n}\right)_{n\in\mathbb{N}^+} $ in $A$ such that 
	\begin{equation*}
		\lim_{n\rightarrow \infty }x_{n}=x,
	\end{equation*}
	we have 
	\begin{equation*}
		F\left( x\right) \leq \liminf\limits_{n\rightarrow \infty}F\left(
		x_{n}\right) .
	\end{equation*}
	$F$ is called \emph{lower semicontinuous} if it is lower semicontinuous at
	any point $x\in A$.
\end{de}

Observe that $F$ is lower semicontinuous if $F$ is continuous.

\begin{Thm}\label{existence}
	\label{First}Assume that the Banach space $X$ is such that every bounded
	sequence in $X$ has a convergent subsequence. Assume that there exists a
	bounded minimizing sequence. If $A$ is closed and $F$ is lower
	semicontinuous, then there exists a minimum of $F$.
\end{Thm}

\begin{proof}
	Let $\left( x_{n}\right)_{n\in\mathbb{N}^+} $ be a bounded minimizing sequence. Suppose that 
	$A$ is closed and that $F$ is lower semicontinuous. By the assumption on $X$, there
	exists a convergent subsequence $\left( x_{n_{k}}\right)_{k\in\mathbb{N}^+} $ of $\left(
	x_{n}\right)_{n\in\mathbb{N}^+} $. Let $\overline{x}$ be the limit of $\left( x_{n_{k}}\right)_{k\in\mathbb{N}^+} $%
	. Since $A$ is closed, we have $\overline{x}\in A$. Since $F$ is lower
	semicontinuous, we have
	\begin{equation*}
		F\left( \overline{x}\right) \leq \liminf_{k\rightarrow \infty }F\left(
		x_{n_{k}}\right) .
	\end{equation*}%
	Since 
	\begin{equation*}
		\liminf_{k\rightarrow \infty }F\left( x_{n_{k}}\right)
		=\lim\limits_{k\rightarrow \infty}F\left( x_{n_{k}}\right) =\inf_{x\in
			A}F\left( x\right) ,
	\end{equation*}%
	we obtain 
	\begin{equation*}
		F\left( \overline{x}\right) \leq \inf_{x\in A}F\left( x\right)
	\end{equation*}%
	and then, recalling that $\overline{x}\in A$, we can conclude that 
	\begin{equation*}
		F\left( \overline{x}\right) =\min_{x\in A}F\left( x\right) .
	\end{equation*}
\end{proof}

The property in Theorem \ref{existence} that every bounded sequence in $X$ has a convergent subsequence
holds if and only if $X$ is finite-dimensional. This is the content of the 
\emph{Heine--Borel theorem}.

If $A$ is unbounded, the assumption in Theorem~\ref{existence} that there exists a bounded minimizing sequence is satisfied whenever $F$ is coercive.

\subsection{The Banach--Saks property}\label{Asub3}

In this subsection, we present a theorem to be used when the minimization problem is infinite-dimensional, i.e., when $X$ is infinite-dimensional. Specific notions need to be introduced.

\begin{de}
	The non-empty subset $A$ of $X$ is called \emph{convex} if, for all $%
	x_{1},x_{2}\in A$ and $\alpha \in \lbrack 0,1]$, we have 
	\begin{equation*}
		(1-\alpha )x_{1}+\alpha x_{2}\in A.
	\end{equation*}%
	The functional $F$ is called $\emph{convex}$ if $A$ is convex and for all $%
	x_{1},x_{2}\in A$ and $\alpha \in \lbrack 0,1]$, we have 
	\begin{equation*}
		F((1-\alpha )x_{1}+\alpha x_{2})\leq (1-\alpha )F(x_{1})+\alpha F(x_{2}).
	\end{equation*}
\end{de}

\begin{remark}\label{reA}
Suppose that $A$ is convex and $F$ is convex. Then, for every $n\in\mathbb{N}^+$, $%
x_1,\ldots,x_n\in A$ and $\alpha_1,\ldots,\alpha_n\in[0,1]$ with $%
\sum\limits_{i=1}^n \alpha_i=1$, we have 
\begin{equation*}
	\sum\limits_{i=1}^n \alpha_i x_i\in A
\end{equation*}
and 
\begin{equation*}
	F\left(\sum\limits_{i=1}^n \alpha_i x_i\right)\leq \sum\limits_{i=1}^n
	\alpha_i F(x_i).
\end{equation*}
\end{remark}

\begin{de}
	A sequence $\left( x_{n}\right)_{n\in\mathbb{N}^+}$ in $X$ is said to be \emph{Ces\`{a}ro
		convergent} to $\overline{x}\in X$ if 
	\begin{equation*}
		\lim_{N\rightarrow \infty }\frac{1}{N}\sum\limits_{n=1}^{N}x_{n}=\overline{x}
		.
	\end{equation*}
\end{de}

\begin{de}
	The Banach space $X$ is said to have the \emph{Banach--Saks property} if every bounded
	sequence in $X$ has a Ces\`{a}ro convergent subsequence.
\end{de}

For every $p\in(1,+\infty)$, the space $L^p(a,b)$ has the Banach--Saks
property (see the original paper \cite{BS} by Banach and Saks, or the proof given in \cite{Fo}, or use the fact that the Banach--Saks property holds for uniformly convex spaces and $L^p(a,b)$ is uniformly convex). 

In the case where $X$ has the Banach--Saks property, the next theorem gives conditions under which a minimum of $F$ on $A$ exists. 

\begin{Thm}\label{existence2}
	Assume that $X$ has the Banach--Saks property. Assume that there exists a
	bounded minimizing sequence. If $A$ is closed and convex and $F$ is lower
	semicontinuous and convex, then there exists a minimum of $F$.
\end{Thm}

\begin{proof}
	Let $\left( x_{n}\right)_{n\in\mathbb{N}^+}$ be a bounded minimizing sequence. Suppose that  
	$A$ is closed and convex and that $F$ is lower semicontinuous and convex. By the
	assumption on $X$, there exists a Ces\`{a}ro convergent subsequence $\left(
	x_{n_{k}}\right)_{k\in\mathbb{N}^+} $. Let 
	\begin{equation*}
		\overline{x}=\lim_{K\rightarrow \infty }\frac{1}{K}\sum%
		\limits_{k=1}^{K}x_{n_{k}}.
	\end{equation*}%
	Since $A$ is convex, we have, for any $K\in \mathbb{N}^{+}$,  
	\begin{equation*}
		\frac{1}{K}\sum\limits_{k=1}^{K}x_{n_{k}}\in A
	\end{equation*}%
(recall Remark \ref{reA}). Therefore, since $A$ is closed, we have $\overline{x}\in
A$. Since $F$ is convex, we have, for any $K\in \mathbb{N}^{+}$, 
	\begin{equation}
		F\left( \frac{1}{K}\sum\limits_{k=1}^{K}x_{n_{k}}\right) \leq \frac{1}{K}%
		\sum\limits_{k=1}^{K}F\left( x_{n_{k}}\right)  \label{F<=F}
	\end{equation}
	(recall Remark \ref{reA}). Since $F$ is lower semicontinuous, we obtain, by passing to the limit as $%
	K\rightarrow \infty $ in (\ref{F<=F}), 
	\begin{equation*}
		F\left( \overline{x}\right) \leq \liminf\limits_{K\rightarrow \infty
		}F\left( \frac{1}{K}\sum\limits_{k=1}^{K}x_{n_{k}}\right) \leq
		\liminf\limits_{K\rightarrow \infty }\frac{1}{K}\sum\limits_{k=1}^{K}F%
		\left( x_{n_{k}}\right).
	\end{equation*}%
	Since 
	\begin{equation*}
		\lim_{k\rightarrow \infty }F\left( x_{n_{k}}\right) =\inf_{x\in A}F\left(
		x\right),
	\end{equation*}%
	we have 
	\begin{equation*}
		\lim\limits_{K\rightarrow +\infty }\frac{1}{K}\sum\limits_{k=1}^{K}F\left(
		x_{n_{k}}\right) =\inf_{x\in A}F\left( x\right).
	\end{equation*}%
	This follows from the fact that convergent sequences are Ces\`{a}ro
	convergent to the same limit. We conclude that 
	\begin{equation*}
		F\left( \overline{x}\right) \leq \liminf\limits_{K\rightarrow +\infty }\frac{%
			1}{K}\sum\limits_{k=1}^{K}F\left( x_{n_{k}}\right)
		=\lim\limits_{K\rightarrow +\infty }\frac{1}{K}\sum\limits_{k=1}^{K}F\left(
		x_{n_{k}}\right) =\inf_{x\in A}F\left( x\right)
	\end{equation*}%
	and then, since $\overline{x}\in A$, 
	\begin{equation*}
		F\left( \overline{x}\right) =\min_{x\in A}F\left( x\right) .
	\end{equation*}
\end{proof}

As in Theorem \ref{existence}, if $A$ is unbounded, the assumption in Theorem \ref{existence2} that there exists a bounded minimizing sequence is satisfied whenever $F$ is coercive.

Consider the Banach space 
\begin{equation*}
	X=X_{1}\times \cdots\times X_{M},
\end{equation*}%
the Cartesian product of the $M$ Banach spaces $X_{1},\ldots,X_{M}$ with norms $%
\left\Vert \ \cdot\ \right\Vert_{X_1},\ldots,$\ $\left\Vert \ \cdot\
\right\Vert_{X_M}$, respectively. The Banach space $X$ is equipped with the
norm 
\begin{equation*}
	\Vert x\Vert_X =\max\left\{\Vert x_1\Vert_{X_1},\ldots,\Vert x_M
	\Vert_{X_M}\right\},\ x\in X.
\end{equation*}
Regarding the Banach--Saks property of $X$, the following
theorem holds.

\begin{Thm} \label{BSC}
	If the Banach spaces $X_{1},\ldots, X_{M}$ have the Banach--Saks property,
	then $X$ has the Banach--Saks property.
\end{Thm}

The proof of this theorem uses the fact that, in a Banach space with the
Banach--Saks property, every bounded sequence has a Cesàro convergent
subsequence such that each of its subsequences is Cesàro convergent to the same
limit. This result is included in the theorem presented in \cite{EM}.

\begin{proof}
	Let 
	\begin{equation*}
		\left( x_{n}\right)_{n\in\mathbb{N}^+} =\left( \left( x_{n1},\ldots,x_{nM}\right) \right)_{n\in\mathbb{N}^+}
	\end{equation*}%
	be a bounded sequence in $X$. Then, for any $i\in\{1,\ldots,M\}$, the
	sequence $\left( x_{ni}\right) $ is a bounded sequence in $X_{i}$.
	
	Consider the bounded sequence $\left( x_{n1}\right)_{n\in\mathbb{N}^+} $ in $X_{1}$. Since $%
	X_{1}$ has the Banach--Saks property, there exists a subsequence $\left(
	x_{n_{k_{1}}1}\right)_{k_1\in\mathbb{N}^+} $ of $\left( x_{n1}\right)_{n\in\mathbb{N}^+} $ Ces\`{a}ro convergent to $%
	\overline{x}_{1}$ and such that each of its subsequences is Cesàro convergent to $	\overline{x}_{1}$.
	
	Now, consider the bounded sequence $\left( x_{n_{k_{1}}2}\right)_{k_1\in\mathbb{N}^+} $ in $X_{2}$%
	. Since $X_{2}$ has the Banach--Saks property, there exists a subsequence $%
	\left( x_{n_{{k_{1}}_{k_{2}}}2}\right)_{k_2\in\mathbb{N}^+} $ of $\left( x_{n_{k_{1}}2}\right)_{k_1\in\mathbb{N}^+}$
	Ces\`{a}ro convergent to $\overline{x}_{2} $  and such that each of its subsequences is Cesàro convergent to $	\overline{x}_{2}$.
	
	By proceeding in this manner, we see that, for every $i\in\{1,\ldots, M-1\}$, there
	exists a subsequence 
	\begin{equation*}
		\left( x_{n_{{k_{1}}_{\ddots_{{k_{i}}_{k_{i+1}}}}}i+1}\right)_{k_{i+1}\in\mathbb{N}^+},
	\end{equation*}
	of the bounded sequence 
	\begin{equation*}
		\left( x_{n_{{k_{1}}_{\ddots_{k_{i}}}}i+1}\right)_{k_i\in\mathbb{N}^+} 
	\end{equation*}
	in $X_{i+1}$, Ces\`{a}ro convergent to $\overline{x}_{i+1}$  and such that each of its subsequences is Cesàro convergent to $	\overline{x}_{i+1}$.
	
	Thus, for any $i\in \{1,\ldots,M\}$, the sequence 
	\begin{equation*}
		\left(x_{n_{{k_{1}}_{\ddots_{k_{i}}}}i}\right)_{k_i\in\mathbb{N}^+}
	\end{equation*}%
	in $X_{i}$ is Ces\`{a}ro convergent to $\overline{x}_{i}$ and its subsequence
	\begin{equation*}
		\left(x_{n_{{k_{1}}_{\ddots_{k_{M}}}}i}\right)_{k_M\in\mathbb{N}^+}
	\end{equation*}%
	is also Ces\`{a}ro convergent to $\overline{x}_{i}$.
	
	For the subsequence%
	\begin{equation*}
		\left(x_{n_{{k_{1}}_{\ddots_{k_{M}}}}}\right)_{k_M\in\mathbb{N}^+}
	\end{equation*}
	of the bounded sequence $(x_n)_{n\in\mathbb{N}^+}$ in $X$, we have, for any $K\in \mathbb{N}^{+}$, 
	\begin{equation*}
		\frac{1}{K}\sum\limits_{k_{M}=1}^{K}x_{n_{{k_{1}}_{\ddots_{k_{M}}}}}=\left( 
		\frac{1}{K}\sum\limits_{k_{M}=1}^{K}x_{n_{{k_{1}}_{\ddots_{k_{M}}}}1},\ldots,%
		\frac{1}{K}\sum\limits_{k_{M}=1}^{K}x_{n_{{k_{1}}_{\ddots_{k_{M}}}}M}\right)
		.
	\end{equation*}
	Thus%
	\begin{eqnarray*}
		&&\lim\limits_{K\rightarrow \infty }\frac{1}{K}\sum\limits_{k_{M}=1}^{K}x_{n_{{%
					k_{1}}_{\ddots_{k_{M}}}}}\\
		&&=\left( \lim_{K\rightarrow \infty }\frac{1}{K}%
		\sum\limits_{k_{M}=1}^{K}x_{n_{{k_{1}}_{\ddots_{k_{M}}}}1},\ldots,\lim_{K%
			\rightarrow \infty }\frac{1}{K}\sum\limits_{k_{M}=1}^{K}x_{n_{{k_{1}}%
				_{\ddots_{k_{M}}}}M}\right) \\
		&&=\left(\overline{x}_1,\ldots,\overline{x}_M\right).
	\end{eqnarray*}
\end{proof}

\subsection{Uniqueness}\label{Asub4}

The uniqueness of a minimum of $F$ is guaranteed under strict convexity of 
$F$.

\begin{de}
	The functional $F$ is called $\emph{strictly\ convex}$ if $A$ is convex and
	for all $x_{1},x_{2}\in A$ with $x_{1}\neq x_{2}$, we have 
	\begin{equation*}
		F((1-\alpha )x_{1}+\alpha x_{2})<(1-\alpha )F(x_{1})+\alpha F(x_{2})\ \text{%
			for every }\alpha \in \left( 0,1\right) .
	\end{equation*}
	
	$F$ is called $\emph{strictly\ convex}$ \emph{in relaxed form} if $A$ is convex
	and for all $x_{1},x_{2}\in A$ with $x_{1}\neq x_{2}$, we have%
	\begin{equation*}
		F((1-\alpha )x_{1}+\alpha x_{2})<(1-\alpha )F(x_{1})+\alpha F(x_{2})\ \text{%
			for some }\alpha \in \left( 0,1\right) .
	\end{equation*}
\end{de}

Strict convexity implies convexity and, if $F$ is lower semicontinuous or convex, strict
convexity in relaxed form is equivalent to strict convexity.

\begin{Thm}\label{uniqueness}
	If $F$ is strictly convex in relaxed form, then there exists at
	most one minimum of $F$.
\end{Thm}

\begin{proof}
	Suppose that $F$ is strictly convex in relaxed form. Let $\overline{x}_{1},%
	\overline{x}_{2}\in A$, with $\overline{x}_{1}\neq \overline{x}_{2}$, be minima
	of $F$ on $A$ with minimum value $F\left( \overline{x}_{1}\right) =F\left( 
	\overline{x}_{2}\right) $. We have 
	\begin{equation*}
		F((1-\alpha )\overline{x}_{1}+\alpha \overline{x}_{2})<(1-\alpha )F(%
		\overline{x}_{1})+\alpha F(\overline{x}_{2})=F(\overline{x}_{1})
	\end{equation*}%
	for some $\alpha \in \left( 0,1\right) $. This is not possible, since $%
	(1-\alpha )\overline{x}_{1}+\alpha \overline{x}_{2}\in A$ and $F(\overline{x}%
	_{1})$ is the minimum value of $F$ on $A$.
\end{proof}

\end{document}